\documentclass{article}
\usepackage{amsmath}
\usepackage{amsthm}
\usepackage{amsfonts}
\usepackage{amssymb}
\usepackage{biblatex}
\usepackage{tikz-cd}
\usepackage{hyperref}
\usepackage{geometry}

\theoremstyle{definition}
\newtheorem{definition}{Definition}[section]
\newtheorem{condition}[definition]{Condition}
\theoremstyle{plain}
\newtheorem{lemma}[definition]{Lemma}
\newtheorem{theorem}[definition]{Theorem}
\newtheorem*{theorem*}{Theorem}
\newtheorem{proposition}[definition]{Proposition}
\newtheorem{corollary}[definition]{Corollary}
\theoremstyle{remark}
\newtheorem{remark}[definition]{Remark}

\title{A GIT construction of moduli spaces of sheaves of length 2}
\date{\today}
\author{Yikun Qiao}

\begin{document}
\maketitle


\begin{abstract}
	Let $\Bbbk$ be an algebraically closed field of characteristic zero. Let $\mathrm{Sch}/\Bbbk$ denote the category of schemes of finite type over $\Bbbk$. Let $B$ be a connected projective scheme over $\Bbbk$ and let $\mathcal L$ be an ample line bundle on $B$. Let $\tau$ be a Harder-Narasimhan type of length 2, and let $\delta\in\mathbb N$. We say a pure sheaf $\mathcal E$ on $B$ is \emph{$(\tau,\delta)$-stable} if its Harder-Narasimhan filtration $0=\mathcal E_{\leq 0}\subsetneq\mathcal E_{\leq 1}\subsetneq\mathcal E_{\leq 2}=\mathcal E$ is non-splitting, of type $\tau$, with stable subquotients, and $\delta=\dim_\Bbbk\mathrm{Hom}_{\mathcal O_B}(\mathcal E_2,\mathcal E_1)$ for $\mathcal E_i:=\mathcal E_{\leq i}/\mathcal E_{\leq i-1}$. 

	We define a moduli functor $\mathbf M'_{\tau,\delta}$ classifying $(\tau,\delta)$-stable sheaves on $B$ and construct its coarse moduli space by non-reductive geometric invariant theory (GIT). We extend the non-reductive GIT in \cite{BércziGergely2016Pcog} and \cite{BércziGergely2018Gitf} to linear actions on non-reduced schemes, and apply our non-reductive GIT to prove that the sheafification $(\mathbf M'_{\tau,\delta})^\sharp$ on $(\mathrm{Sch}/\Bbbk)_{\acute etale}$ is represented by a quasi-projective scheme. 

	Our methods generalise Jackson's construction of moduli spaces of $(\tau,\delta)$-stable sheaves in \cite{JacksonJoshua2021MSoU} in the category of varieties, to allow non-reduced moduli schemes. 
\end{abstract}

\tableofcontents

\section{Introduction}
{
	Let $\Bbbk$ be an algebraically closed field of characteristic zero. Let $\mathrm{Sch}/\Bbbk$ be the category of schemes of finite type over $\Bbbk$. In this situation, a closed point of a scheme is equivalent to a $\Bbbk$-point. Let $B$ be a connected projective scheme over $\Bbbk$ and let $\mathcal L$ be an ample line bundle on $B$. Let $f:B\to\mathrm{Spec}(\Bbbk)$ be the structure morphism. 

	Let $P_1(t),P_2(t)\in\mathbb Q[t]$ be two numerical polynomials satisfying 
	\begin{equation}
		\frac{P_1(t)}{\alpha_{1,d_1}}>\frac{P_2(t)}{\alpha_{2,d_2}}
	\end{equation}
	where $\alpha_{i,d_i}\in\mathbb N_+$ is the integer such that $\frac{\alpha_{i,d_i}}{d_i!}t^{d_i}$ is the leading monomial of $P_i(t)$, and the total order on $\mathbb Q[t]$ is the lexicographical order. Let $\tau$ be a Harder-Narasimhan type of length 2 associated to Hilbert polynomials $P_1(t),P_2(t)$, in the sense that we say a pure sheaf $\mathcal E$ on $B$ has type $\tau$ if $\mathcal E_i:=\mathcal E_{\leq i}/\mathcal E_{\leq i-1}$ is semistable of Hilbert polynomial $P_i(t)$ for $i=1,2$, where 
	\begin{equation}
		0=\mathcal E_{\leq 0}\subsetneq\mathcal E_{\leq 1}\subsetneq \mathcal E_{\leq 2}=\mathcal E
	\end{equation}
	is the Harder-Narasimhan filtration of $\mathcal E$. Let $\delta\in\mathbb N$, which we choose as another discrete invariant. For a pure sheaf $\mathcal E$ of type $\tau$, we say that the space of \emph{unipotent endomorphisms} of $\mathcal E$ has dimension $\delta$ if $\delta=\dim_\Bbbk \mathrm{Hom}_{\mathcal O_B}(\mathcal E_2,\mathcal E_1)$. 


	Jackson has studied the moduli problem of \emph{$(\tau,\delta)$-stable sheaves} in \cite{JacksonJoshua2021MSoU}. A sheaf $\mathcal E$ on $B$ is called \emph{$(\tau,\delta)$-stable} if 
	\begin{itemize}
		\item $\mathcal E$ has type $\tau$; 
		\item $\mathcal E_i$ is stable for $i=1,2$; 
		\item $\mathcal E\not\cong\mathcal E_1\oplus\mathcal E_2$; 
		\item $\delta=\dim\mathrm{Hom}_{\mathcal O_B}(\mathcal E_2,\mathcal E_1)$. 
	\end{itemize}
	He used non-reductive geometric invariant theory developed in \cite{BércziGergely2016Pcog} and \cite{BércziGergely2018Gitf} to construct moduli spaces of $(\tau,\delta)$-stable sheaves. However, the non-reductive GIT in \cite{BércziGergely2016Pcog} and \cite{BércziGergely2018Gitf} is for actions on \emph{reduced} projective schemes, while the desired moduli spaces may not be necessarily reduced (See \cite{https://doi.org/10.48550/arxiv.0911.2301} Remark 7.8. ). The paper aims to extend results in \cite{JacksonJoshua2021MSoU} to allow potentially non-reduced moduli schemes. 

	For $T\in\mathrm{Sch}/\Bbbk$, we call a coherent sheaf $\mathcal E$ on $B_T:=B\times_\Bbbk T$ a \emph{flat family over $T$ of $(\tau,\delta)$-stable sheaves on $B$} if: 
	\begin{itemize}
		\item there is a length 2 filtration of $\mathcal O_{B_T}$-submodules 
		\begin{equation}
			0=\mathcal E_{\leq 0}\subsetneq\mathcal E_{\leq 1}\subsetneq \mathcal E_{\leq 2}=\mathcal E
		\end{equation}
		such that $\mathcal E_i:=\mathcal E_{\leq i}/\mathcal E_{\leq i-1}$ is a flat family over $T$ of stable sheaves of Hilbert polynomial $P_i(t)$ for $i=1,2$; 
		\item $\mathcal E_{\leq 1}\subsetneq\mathcal E_{\leq 2}$ does not split at any closed $t\in T$; 
		\item for every closed subscheme $T'\hookrightarrow T$, the following $\mathcal O_{T'}$-module is locally free of rank $\delta$
		\begin{equation}
			(B_{T'}\to T')_*\mathcal Hom_{\mathcal O_{B_{T'}}}\big(\mathcal E|_{B_{T'}},\mathcal E_1|_{B_{T'}}\big). 
		\end{equation}
	\end{itemize}
	We define a moduli functor $\mathbf M'_{\tau,\delta}:(\mathrm{Sch}/\Bbbk)^{\mathrm{op}}\to\mathrm{Set}$ which to $T\in\mathrm{Sch}/\Bbbk$ assigns the set 
	\begin{equation}
		\mathbf M'_{\tau,\delta}(T):=\Big\{\begin{matrix}\textrm{isomorphism classes of flat families}\\\textrm{over }T\textrm{ of }(\tau,\delta)\textrm{-stable sheaves on }B\end{matrix}\Big\}. 
	\end{equation}
	The final result of this paper is the following. 
	{
		\begin{theorem*}
			Let $\Bbbk$ be an algebraically closed field of characteristic zero. Let $\mathrm{Sch}/\Bbbk$ denote the category of schemes of finite type over $\Bbbk$. Let $\tau$ be a Harder-Narasimhan type of length 2. Let $\delta\in\mathbb N$. Then there exists a quasi-projective scheme $M_{\tau,\delta}\in\mathrm{Sch}/\Bbbk$ representing $(\mathbf M'_{\tau,\delta})^\sharp$, the sheafification of $\mathbf M'_{\tau,\delta}$ with respect to the étale topology. In particular, the moduli functor $\mathbf M'_{\tau,\delta}$ has a coarse moduli space $M_{\tau,\delta}$. 
		\end{theorem*}
	}

	In Section 2, we first introduce the notion of \emph{flat families of $(\tau,\delta)$-stable sheaves on $B$} parametrised by schemes (Definition \ref{definition: family of (tau,delta)-stable sheaves}), which generalises the notion of $(\tau,\delta)$-stable sheaves. We consider the moduli functor $\mathbf M'_{\tau,\delta}$ sending $T\in\mathrm{Sch}/\Bbbk$ to the set of isomorphism classes of flat families over $T$ of $(\tau,\delta)$-stable sheaves on $B$. Next, we set up the non-reductive GIT problem in the category of schemes. For $m\gg0$ and $V:=\Bbbk^{P(m)}$, there exist locally closed subschemes $Y_\tau^{\mathrm{(s)s}}\subseteq Y_\tau\subseteq \mathrm{Quot}_{V(-m)/B/\Bbbk}^{P(t),\mathcal L}$, where $V(-m):=V\otimes_\Bbbk\mathcal L^{-m}$. The schemes $Y_\tau^{\mathrm{(s)s}}\subseteq Y_\tau$ are determined by their functors of points. For example, a morphism $T\to Y_\tau^{\mathrm{s}}$ is equivalent to an isomorphism class of a flat family over $T$ of quotients of $V(-m):=V\otimes_\Bbbk\mathcal L^{-m}$
	\begin{equation}
		q:V(-m)_T\to\mathcal E,\quad V(-m)_T:=(B_T\to B)^*\big(V(-m)\big)
	\end{equation}
	denoted by $[(\mathcal E,q)]$, such that $\mathcal E$ is a flat family over $T$ of sheaves of type $\tau$ with stable subquotients. Let $\overline{Y_\tau^{\mathrm{s}}}\subseteq \mathrm{Quot}_{V(-m)/B/\Bbbk}^{P(t),\mathcal L}$ denote the scheme theoretic closure. Let $\mathrm P_\tau\subseteq \mathrm{SL}(V)$ be the parabolic subgroup stabilising the partial flat $0=V_{\leq 0}\subseteq V_{\leq 1}\subseteq V_{\leq 2}=V$, where $V_i:=\Bbbk^{P_i(m)}$ and $V_{\leq i}:=V_1\oplus\cdots\oplus V_i$. Recall that $\mathrm{Quot}_{V(-m)/B/\Bbbk}^{P(t),\mathcal L}$ has a very ample line bundle $\mathcal M$. The non-reductive GIT problem we consider is 
	\begin{equation}
		\mathrm P_\tau\curvearrowright\big(\overline{Y_\tau^{\mathrm{s}}},\mathcal M\big|_{\overline{Y_\tau^{\mathrm{s}}}}\big). 
	\end{equation}
	This is an enhancement of the non-reductive GIT problem considered in \cite{JacksonJoshua2021MSoU}, which was essentially the reduction of the linearisation above. 

	In Section 3, we develop the non-reductive GIT for groups of the form $\hat U:=U\rtimes\mathbb G_m$, where $U:=(\mathbb G_a)^r$ is additive and $\mathbb G_m\to \mathrm{Aut}(U)$ induces a representation of $\mathbb G_m$ on $\mathfrak u:=\mathrm{Lie}(U)$ with a unique weight which is positive. Non-reductive GIT for linear $\hat U$-actions on reduced projective schemes has been considered in \cite{BércziGergely2016Pcog}, \cite{BércziGergely2018Gitf} and \cite{hoskins2021quotients}. We follow their ideas and generalise to linear $\hat U$-actions on non-reduced projective schemes. Let $X$ be a projective scheme and let $L$ be an ample line bundle on $X$. Let $\hat U$ act linearly on $X$ with respect to $L$.

	There is a morphism $\phi:\Omega_{X/\Bbbk}\to \mathfrak u^*\otimes_\Bbbk\mathcal O_X$ and we call $\phi$ the \emph{infinitesimal action of $\mathfrak u$ on $X$}. There are important Condition \hyperref[condition: UU]{UU} and Condition \hyperref[condition: WUU]{WUU}, generalising \cite{hoskins2021quotients} Assumption 4.40 and 4.42 for $\hat U\curvearrowright(X,L)$. Condition \hyperref[condition: UU]{UU} assumes that $X^0_{\min}\ne\emptyset$ and that $\mathrm{coker}(\phi)|_{X^0_{\min}}$ is locally free of constant rank. Condition \hyperref[condition: WUU]{WUU} assumes roughly that $\mathrm{coker}(\phi)$ is locally free of the lowest possible rank on an open subset which intersects with $Z_{\min}:=(X^0_{\min})^{\mathbb G_m}$, where $Z_{\min}$ is defined as the \emph{fixed point subscheme} in the sense of \cite{FogartyJohn1973FPS}. These conditions reduce to \cite{hoskins2021quotients} Assumption 4.40 and 4.42 when $X$ is reduced. 

	When Condition \hyperref[condition: UU]{UU} holds, the open subscheme $X^0_{\min}$ has a quasi-projective universal geometric quotient by $U$ (Theorem \ref{theorem: quotient for projective with UU}), and $X^0_{\min}\setminus UZ_{\min}$ has a projective universal geometric quotient by $\hat U$ (Theorem \ref{theorem: projective geometric quotient by U-hat with UU}). Moreover the quotient $X^0_{\min}/U$ represents the associated quotient sheaf on $(\mathrm{Sch}/\Bbbk)_{\acute etale}$ (Corollary \ref{corollary: geometric U-quotients represent quotient sheaves with UU}). When Condition \hyperref[condition: WUU]{WUU} holds, we construct a $\hat U$-equivariant blowing up $b:\widetilde X\to X$ such that Condition \hyperref[condition: UU]{UU} holds for $\hat U\curvearrowright \big(\widetilde X,\widetilde L\big)$, where $\widetilde L:=\mathcal O_{\widetilde{X}}(-E)\otimes_{\mathcal O_{\widetilde X}}b^*L^a$ for $a\gg0$ and $E\subseteq \widetilde X$ is the exceptional divisor. 

	In general, we only assume $X^0_{\min}\ne\emptyset$. We consider the \emph{stratification by unipotent stabilisers} $\coprod_{\delta\in\mathbb N}S_\delta(X^0_{\min})\to X^0_{\min}$, where $S_\delta(X^0_{\min})\subseteq X^0_{\min}$ is the locally closed subscheme which is universal to the property that $\mathrm{coker}(\phi)|_{S_\delta(X^0_{\min})}$ is locally free of rank $\delta$. If $S_\delta(X^0_{\min})\cap Z_{\min}\ne\emptyset$, then Condition \hyperref[condition: WUU]{WUU} holds for $\hat U\curvearrowright\big(\overline{S_\delta(X^0_{\min})},L\big|_{\overline{S_\delta(X^0_{\min})}}\big)$, where $\overline{S_\delta(X^0_{\min})}\subseteq X$ is the scheme theoretic closure. We can apply our results with Condition \hyperref[condition: WUU]{WUU} to get, for example, Theorem \ref{theorem: geometric U and U-hat quotients of strata} and Theorem \ref{theorem: geometric quotients of strata by parabolic}. 

	In Section 4, we apply our non-reductive GIT for the moduli problem $\mathbf M'_{\tau,\delta}$. We have defined $Y_\tau^{\mathrm{(s)s}}\subseteq Y_\tau$ in Section 2. Let $Y_{\tau,\delta}^{\mathrm{s}}:=S_\delta(Y_\tau^{\mathrm{s}})\setminus \mathrm U_\tau Z_{\min}(\overline{Y_\tau^{\mathrm{s}}})$, which is open in $S_\delta(Y_\tau^{\mathrm{s}})$. Similar to $Y_\tau^{\mathrm{(s)s}}\subseteq Y_\tau$, we can descibe $Y_{\tau,\delta}^{\mathrm{s}}\subseteq S_\delta(Y_\tau^{\mathrm{s}})$ by their functors of points. For a morphism $T\to\mathrm{Quot}_{V(-m)/B/\Bbbk}^{P(t),\mathcal L}$, it factors through $Y_{\tau,\delta}^{\mathrm{s}}$ if the associated class of families $[(\mathcal E,q)]$ over $T$ satisfies that $\mathcal E$ is a flat family over $T$ of $(\tau,\delta)$-stable sheaves (Lemma \ref{lemma: moduli descriptions of S_alpha(Y_tau^s) and Y_(tau,alpha)^s}). This links the moduli problem $\mathbf M'_{\tau,\delta}$ with the GIT problem in Section 2. By Theorem \ref{theorem: geometric quotients of strata by parabolic}, the action $\mathrm P_\tau\curvearrowright Y_{\tau,\delta}^{\mathrm{s}}$ has a quasi-projective universal geometric quotient $Y_{\tau,\delta}^{\mathrm{s}}/\mathrm P_\tau$. The quotient $Y_{\tau,\delta}^{\mathrm{s}}\to Y_{\tau,\delta}^{\mathrm{s}}/\mathrm P_\tau$ is the composition of a $\mathrm U_\tau$-quotient and a $\mathrm P_\tau/\mathrm U_\tau$-quotient. In both stages, the GIT quotients represent the associated quotient sheaves on $(\mathrm{Sch}/\Bbbk)_{\acute etale}$. This observation together with a standard categorical argument (for example \cite{HuybrechtsDaniel2010Tgom} Lemma 4.3.1) shows that $(\mathbf M'_{\tau,\delta})^\sharp\cong Y_{\tau,\delta}^{\mathrm{s}}/\mathrm P_\tau$. 
}

\section{The moduli problem and the GIT problem}
{
	Recall that $\mathrm{Sch}/\Bbbk$ is the category of schemes of finite type over $\Bbbk$, and $B$ is a projective scheme over $\Bbbk$ with an ample line bundle $\mathcal L$. Recall $P_i(t),P(t)\in\mathbb Q[t]$ are prescribed Hilbert polynomials. 

	For an integer $n\in\mathbb Z$ and a vector space $V$, let $V(-n):=V\otimes_\Bbbk\mathcal L^{-n}$ denote the coherent sheaf on $B$ twisted from the constant sheaf by $\mathcal L^{-n}$. For $T\in\mathrm{Sch}/\Bbbk$ and $\mathcal F$ a coherent sheaf on $B$, let $B_T:=B\times_\Bbbk T$ denote the base change of $B$ and let $\mathcal F_T:=(B_T\to B)^*\mathcal F$ denote the base change of $\mathcal F$. Consider $\mathrm{Quot}_{V(-n)/B/\Bbbk}^{P,\mathcal L}$, the Quot scheme parametrising quotients $q:V(-n)\to \mathcal E$ such that the $\mathcal E$ has Hilbert polynomial $P(t)\in\mathbb Q[t]$ with respect to $\mathcal L$. 

	\subsection{The moduli problem}\label{subsection: the moduli problem}
	{
		Let $T\in\mathrm{Sch}/\Bbbk$. Recall that $\tau=(P_1,P_2)$ is a Harder-Narasimhan type, and we say a coherent sheaf $\mathcal E$ on $B$ has type $\tau$ if and only if $\mathcal E$ is pure and there exists a length 2 filtration 
		\begin{equation}
			0=\mathcal E_{\leq 0}\subsetneq \mathcal E_{\leq 1}\subsetneq \mathcal E_{\leq 2}=\mathcal E
		\end{equation}
		such that $\mathcal E_i:=\mathcal E_{\leq i}/\mathcal E_{\leq i-1}$ is semistable of Hilbert polynomial $P_i(t)$ for $i=1,2$. 

		{
			We extend the notion of type $\tau$ Harder-Narasimhan filtration for flat families over $T\in\mathrm{Sch}/\Bbbk$. 
			\begin{definition}\label{definition: type tau HNF}
				For $T\in\mathrm{Sch}/\Bbbk$, let $\mathcal E$ be a coherent $\mathcal O_{B_T}$-module flat over $T$. A \emph{type $\tau$ Harder-Narasimhan filtration} of $\mathcal E$ is a filtration of coherent $\mathcal O_{B_T}$-modules 
				\begin{equation}
					0=\mathcal E_{\leq 0}\subsetneq \mathcal E_{\leq 1}\subsetneq \mathcal E_{\leq 2}=\mathcal E
				\end{equation}
				such that $\mathcal E_i:=\mathcal E_{\leq i}/\mathcal E_{\leq i-1}$ is a flat family over $T$ of semistable sheaves of Hilbert polynomial $P_i(t)$ for $i=1,2$. 

				We say $\mathcal E$ is a \emph{flat family over $T$ of type $\tau$ sheaves on $B$} if it admits a type $\tau$ Harder-Narasimhan filtration. 
			\end{definition}

			\begin{lemma}\label{lemma: HNF is unique when exists and HNF commutes with base change}
				For $T\in\mathrm{Sch}/\Bbbk$ and a coherent $\mathcal O_{B_T}$-module $\mathcal E$. If a type $\tau$ Harder-Narasimhan filtration exits, then the filtration is unique. 

				For a morphism $f:S\to T$ in $\mathrm{Sch}/\Bbbk$ and a coherent $\mathcal O_{B_T}$-module $\mathcal E$ with type $\tau$ Harder-Narasimhan filtration $\{\mathcal E_{\leq i}\}_{i=0}^2$, the pullback $(1_B\times f)^*\mathcal E$ has a type $\tau$ Harder-Narasimhan filtration, which is the pullback of $\{\mathcal E_{\leq i}\}_{i=0}^2$. 
			\end{lemma}
			\begin{proof}
				This is a special case of \cite{https://doi.org/10.48550/arxiv.0909.0891} Theorem 5, where the filtration $\{\mathcal E_{\leq i}\}_{i=0}^2$ is called the \emph{relative} Harder-Narasimhan filtration of $\mathcal E$. 
			\end{proof}
		}

		{
			\begin{definition}\label{definition: family of (tau,delta)-stable sheaves}
				For $\delta\in\mathbb N$ and $T\in\mathrm{Sch}/\Bbbk$, a \emph{flat family over $T$ of $(\tau,\delta)$-stable sheaves on $B$} is a flat family over $T$ of type $\tau$ sheaves $\mathcal E$, which has a Harder-Narasimhan filtration 
				\begin{equation}
					0=\mathcal E_{\leq 0}\subsetneq \mathcal E_{\leq 1}\subsetneq \mathcal E_{\leq 2}=\mathcal E
				\end{equation}
				satisfying: 
				\begin{itemize}
					\item $\mathcal E_i=\mathcal E_{\leq i}/\mathcal E_{\leq i-1}$ is a flat family over $T$ of stable sheaves of Hilbert polynomial $P_i(t)\in\mathbb Q[t]$, for $i=1,2$; 
					\item $\mathcal E_{\leq 1}\subseteq \mathcal E_{\leq 2}$ does not split at any closed point $t\in T$; 
					\item for every closed subscheme $T'\subseteq T$, the quasi-coherent $\mathcal O_{T'}$-module
					\begin{equation}
						(B_{T'}\to T')_*\mathcal Hom_{\mathcal O_{B_{T'}}}\big(\mathcal E_2|_{B_{T'}},\mathcal E_1|_{B_{T'}}\big)
					\end{equation}
					is locally free of rank $\delta$. 
				\end{itemize}
			\end{definition}

			\begin{remark}
				Let $T=\mathrm{Spec}(\Bbbk)$ and let $d\in\mathbb N$. A flat family over $T$ of type $\tau$ sheaves is exactly a pure $\mathcal O_B$-module $\mathcal E$ of Harder-Narasimhan type $\tau$. A flat family over $T$ of $(\tau,\delta)$-stable sheaves is exact a pure $\mathcal O_B$-module $\mathcal E$ of type $\tau$ such that $\mathcal E_1,\mathcal E_2$ are stable and $\mathcal E\not\cong \mathcal E_1\oplus\mathcal E_2$ and $\mathrm{Hom}_{\mathcal O_B}(\mathcal E_2,\mathcal E_1)$ has dimension $\delta$. 
			\end{remark}
		}

		The moduli functor $\mathbf M'_{\tau,\delta}$ is defined as follows
		\begin{equation}
			\mathbf M'_{\tau,\delta}:(\mathrm{Sch}/\Bbbk)^{\mathrm{op}}\to \mathrm{Set}
		\end{equation}
		which sends $T\in\mathrm{Sch}/\Bbbk$ to the set 
		\begin{equation}
			\mathbf M'_{\tau,\delta}(T):=\Big\{\begin{matrix}\textrm{isomorphism classes of flat families}\\\textrm{over }T\textrm{ of }(\tau,\delta)\textrm{-stable sheaves on }B\end{matrix}\Big\}
		\end{equation}
		and for a morphism $g:T'\to T$, the map 
		\begin{equation}
			\mathbf M'_{\tau,\delta}(g):\mathbf M'_{\tau,\delta}(T)\to \mathbf M'_{\tau,\delta}(T')
		\end{equation}
		sends the isomorphism class of $\mathcal E$ to the isomorphism class of the class of $\mathcal E|_{B_{T'}}$, where $\mathcal E|_{\mathcal B_{T'}}:=(1_B\times g)^*\mathcal E$. Let $\mathrm P_\tau$ act trivially on $\mathbf M'_{\tau,\delta}$. 

		Consider the étale site $(\mathrm{Sch}/\Bbbk)_{\acute etale}$. We will show that the sheafification $(\mathbf M'_{\tau,\delta})^\sharp$ is representable (Theorem \ref{theorem: sheafification of M'_(tau,delta) is representable}). The representative is the coarse moduli space of $\mathbf M'_{\tau,\delta}$. The moduli space can be constructed by non-reductive GIT, which is quasi-projective. 
	}
	\subsection{Simpson's and Jackson's moduli spaces}
	{
		We introduce Simpson's moduli spaces of stable and semistable sheaves (\cite{SimpsonCarlosT.1994Moro}) and Jackson's moduli spaces of sheaves of Harder-Narasimhan length 2 (\cite{JacksonJoshua2021MSoU}). These moduli spaces were constructed by reductive GIT and non-reductive GIT respectively. These moduli spaces were coarse moduli spaces of naturally defined moduli functors, so they are independent of the way we set up the GIT problems. Moduli spaces of stable sheaves have natural completions, the corresponding moduli spaces of semistable sheaves. In \cite{JacksonJoshua2021MSoU}, completions of the moduli spaces of sheaves of Harder-Narasimhan length 2 were described as non-reductive quotients of blow-ups of the associated GIT problem. However, it is unclear now whether the completion is independent of the GIT problem, and there is no natural moduli functor associated to the completion yet. 

		\subsubsection{Moduli of stable and semistable sheaves}
		{
			Let $B$ be a connected projective scheme over $\Bbbk$ with an ample line bundle $\mathcal L$. Let $P(t)\in\mathbb Q[t]$ be a numerical polynomial. For $j\in\mathbb Z$ and a coherent sheaf $\mathcal F$ on $B$, we denote $\mathcal F(j):=\mathcal F\otimes_{\mathcal O_B}\mathcal L^j$. Moduli spaces of semistable pure sheaves on $B$ of fixed Hilbert polynomials were constructed by Simpson in \cite{SimpsonCarlosT.1994Moro}. A more recent reference is \cite{HuybrechtsDaniel2010Tgom}. 
			
			Let $\mathrm{Sch}/\Bbbk$ be the category of schemes of finite type over $\Bbbk$. For $T\in\mathrm{Sch}/\Bbbk$, a \emph{family over $T$ of (semi)stable sheaves of Hilbert polynomial $P(t)$} is a coherent sheaf $\mathcal F$ on $B_T:=B\times_{\Bbbk}T$ such that $\mathcal F$ is flat over $T$ and $\mathcal F_s$ is a (semi)stable sheaf on $B$ of Hilbert polynomial $P(t)$ for every closed point $s\in T$, where $\mathcal F_s$ is the pullback of $\mathcal F$ along $B\cong B\times_{\Bbbk}\mathrm{Spec}(\kappa(s))\hookrightarrow B_T$. Two families $\mathcal F,\mathcal F'$ on $T$ are called equivalent if there exists a line bundle $\mathcal A$ on $T$ such that $\mathcal F\cong \mathcal F'\otimes_{\mathcal O_{B_T}}(f_T)^*\mathcal A$, where $f_T:B_T\to T$ is the pullback of $f:B\to \mathrm{Spec}(\Bbbk)$. With the notions of families and their equivalence relations, we consider the moduli functors 
			\begin{equation}
				\begin{split}
					\mathbf{M}_{B/\Bbbk}^{P(t),\mathcal L,\mathrm{ss}}:&(\mathrm{Sch}/\Bbbk)^{\mathrm{op}}\to\mathrm{Set}\\
					\mathbf{M}_{B/\Bbbk}^{P(t),\mathcal L,\mathrm{s}}:&(\mathrm{Sch}/\Bbbk)^{\mathrm{op}}\to\mathrm{Set}
				\end{split}
			\end{equation}
			which associate to $T\in\mathrm{Sch}/\Bbbk$ the sets 
			\begin{equation}
				\begin{split}
					\mathbf{M}_{B/\Bbbk}^{P(t),\mathcal L,\mathrm{(s)s}}(T):=&\bigg\{\begin{matrix}\textrm{families on }T\textrm{ of (semi)stable sheaves}\\\textrm{of Hilbert polynomial }P(t)\end{matrix}\bigg\}\bigg/\textrm{equivalence relation}
				\end{split}
			\end{equation}
			and the functors on morphisms are naturally pullbacks of coherent sheaves. These functors are both corepresentable and $\mathbf{M}_{B/\Bbbk}^{P(t),\mathcal L,\mathrm{s}}$ has a coarse moduli space.

			By \cite{SimpsonCarlosT.1994Moro} Theorem 1.1 and \cite{HuybrechtsDaniel2010Tgom} Theorem 3.3.7, the collection of semistable sheaves with Hilbert polynomial $P(t)$ is bounded. By \cite{HuybrechtsDaniel2010Tgom} Lemma 1.7.6, their regularities are uniformly bounded above. For $m\gg0$, all semistable sheaves $\mathcal F$ of Hilbert polynomial $P(t)$ satisfy: 
			\begin{itemize}
				\item $\dim H^0(B,\mathcal F(m))=P(m)$; 
				\item $H^0(B,\mathcal F(m))\otimes_{\Bbbk}\mathcal L^{-m}\to \mathcal F$ is surjective. 
			\end{itemize}
			Let $V=\Bbbk^{P(m)}$. There is a non-canonical linear isomorphism $V\cong H^0(B,\mathcal F(m))$. With one such isomorphism, we have a surjective morphism 
			\begin{equation}
				q:V(-m)\to \mathcal F
			\end{equation}
			which represents a $\Bbbk$-point in $\mathrm{Quot}_{V(-m)/B/\Bbbk}^{P(t),\mathcal L}$. 
			
			Let $R\subseteq\mathrm{Quot}_{V(-m)/B/\Bbbk}^{P(t),\mathcal L}$ be the open subscheme whose closed points are classes $[(\mathcal F,q)]$ such that: 
			\begin{itemize}
				\item[(1)] $\mathcal F$ is semistable; 
				\item[(2)] the induced map $V\cong H^0(B,V)\to H^0(B,\mathcal F(m))$ is an isomorphism. 
			\end{itemize}
			Let $R^{\mathrm{s}}\subseteq R$ be the open subset of classes $[(\mathcal F,q)]$ such that $\mathcal F$ is stable. Condition $(2)$ above is open since the function $[(\mathcal F,q)]\mapsto\dim H^0(B,\ker(q)(m))$ is upper semi-continuous and the condition that $\dim H^0(B,\ker(q)(m))=0$ is equivalent to the condition $(2)$. Condition $(1)$ above is open by \cite{HuybrechtsDaniel2010Tgom} Proposition 2.3.1. This justifies the claim that $R\subseteq\mathrm{Quot}_{V(-m)/B/\Bbbk}^{P(t),\mathcal L}$ is open, and similarly $R^{\mathrm{s}}\subseteq R$ is open. 

			Recall that for $M\gg m\gg0$, there is a closed immersion 
			\begin{equation}
				\mathrm{Quot}_{V(-m)/B/\Bbbk}^{P(t),\mathcal L}\hookrightarrow \mathrm{Grass}\big(V\otimes_\Bbbk H,P(M)\big)
			\end{equation}
			where $H:=H^0(B,\mathcal L^{M-m})$. This is the classifying morphism of the quotient on $\mathrm{Quot}_{V(-m)/B/\Bbbk}^{P(t),\mathcal L}$ 
			\begin{equation}
				V\otimes_\Bbbk H^0(B,\mathcal L^{M-m})\otimes_\Bbbk\mathcal O_{\mathrm{Quot}}=f_{\mathrm{Quot},*}\big(V_{\mathrm{Quot}}(M-m)\big)\to f_{\mathrm{Quot},*}\mathcal U(M)
			\end{equation}
			where $V(-m)_{\mathrm{Quot}}\to \mathcal U$ is the universal quotient on $B_{\mathrm{Quot}}=B\times_\Bbbk\mathrm{Quot}_{V(-m)/B/\Bbbk}^{P(t),\mathcal L}$, and $f_{\mathrm{Quot}}:B_{\mathrm{Quot}}\to \mathrm{Quot}_{V(-m)/B/\Bbbk}^{P(t),\mathcal L}$. 

			The Pl\"ucker immersion is a closed immersion 
			\begin{equation}
				\mathrm{Grass}\big(V\otimes_\Bbbk H,P(M)\big)\hookrightarrow \mathbb P^N:=\mathbb P\bigg(\Big(\bigwedge^{P(M)}\big(V\otimes_\Bbbk H\big)\Big)^*\bigg). 
			\end{equation}
			The very ample line bundle $\mathcal O_{\mathbb P^N}(1)$ pulls back to a very ample line bundle on $\mathrm{Quot}_{V(-m)/B/\Bbbk}^{P(t),\mathcal L}$, denoted by $\mathcal M$. Then we have 
			\begin{equation}\label{equation: ample line bundle M on Quot}
				\mathcal M=\det\bigg(f_{\mathrm{Quot},*}\Big(\mathcal U\otimes_{\mathcal O_{B_{\mathrm{Quot}}}}\mathcal L_{\mathrm{Quot}}^{\otimes M}\Big)\bigg). 
			\end{equation}

			There is an action of $\mathrm{GL}(V)$ on $\mathrm{Quot}_{V(-m)/B/\Bbbk}^{P(t),\mathcal L}$, which is determined by maps 
			\begin{equation}
				h_{\mathrm{GL}(V)}(T)\times h_{\mathrm{Quot}}(T)\to h_{\mathrm{Quot}}(T),\quad T\in\mathrm{Sch}/\Bbbk. 
			\end{equation}
			Let $(g,q)\in h_{\mathrm{GL}(V)}(T)\times h_{\mathrm{Quot}}(T)$. The morphism $g:T\to \mathrm{GL}(V)$ corresponds to an invertible morphism $g:V\otimes_\Bbbk\mathcal O_T\to V\otimes_\Bbbk\mathcal O_T$ on $T$. Let $q:V(-m)_T\to \mathcal F$ denote the family over $T$ of quotients. We set the image of $(g,q)$ to be the following family of quotients 
			\begin{equation}
				\begin{tikzcd}
					V(-m)_T\ar[rr,"g^{-1}(-m)"]&&V(-m)_T\ar[r,"q"]&\mathcal F. 
				\end{tikzcd}
			\end{equation}
			It is easy to check the above maps are natural in $T\in\mathrm{Sch}/\Bbbk$ and each defines a group action. They define an algebraic group action $\mathrm{GL}(V)\curvearrowright \mathrm{Quot}_{V(-m)/B/\Bbbk}^{P(t),\mathcal L}$. Moreover, the action by $\mathrm{GL}(V)$ is linear with respect to $\mathcal M$. See for example \cite{HuybrechtsDaniel2010Tgom} Section 4.3, where the authors considered right actions. 
			
			We consider the linear action of the subgroup $\mathrm{SL}(V)$. It is easy to see $R$ is invariant under $\mathrm{SL}(V)$. Let $\overline{R}\hookrightarrow \mathrm{Quot}_{V(-m)/B/\Bbbk}^{P(t),\mathcal L}$ be the scheme theoretic image of $R\to \mathrm{Quot}_{V(-m)/B/\Bbbk}^{P(t),\mathcal L}$, which is a subscheme on the closure of $R$ by \href{https://stacks.math.columbia.edu/tag/01R8}{Lemma 01R8}. The linearisation on $\mathcal M$ pulls back to a linearisation $\mathcal M|_{\overline{R}}$. We can apply reductive GIT for the linear action of $\mathrm{SL}(V)$ on $\overline{R}$ with respect to $\mathcal M|_{\overline{R}}$. Note that $\mathcal M$ depends on $M\gg m\gg0$. 
			{
				\begin{theorem}[\cite{HuybrechtsDaniel2010Tgom} Theorem 4.3.3, Theorem 4.3.4]\label{theorem: simpson's construction of moduli of semistable sheaves}
					For $M\gg m\gg0$, we have 
					\begin{itemize}
						\item GIT (semi)stable points are isomorphism classes $[(\mathcal F,q)]\in\overline{R}$ such that $\mathcal F$ is \\(semi)stable and $q:V(-m)\to\mathcal E$ induces an isomorphism $V\cong H^0(B,\mathcal E(m))$, i.e. 
						\begin{equation}
							\overline{R}^{\mathrm{SL}(V),\mathcal M,\mathrm{ss}}=R,\quad \overline{R}^{\mathrm{SL}(V),\mathcal M,\mathrm{s}}=R^{\mathrm{s}}; 
						\end{equation}
						\item Two semistable points $[(\mathcal F,q)],[(\mathcal F',q')]\in R$ are (GIT) $S$-equivalent if and only if $\mathrm{gr}(\mathcal F)\cong \mathrm{gr}(\mathcal F')$, where $\mathrm{gr}(\mathcal F)$ denotes the associated graded sheaf of the Jordan-H\"older filtration of $\mathcal F$; 
						\item A semistable point $[(\mathcal F,q)]\in R$ has closed $\mathrm{SL}(V)$-orbit in $R$ if and only if $\mathcal F$ is polystable, i.e. $\mathcal F$ is a direct sum of stable sheaves; 
						\item $\overline{R}/\!/\mathrm{SL}(V)$ corepresents $\mathbf{M}_{B/\Bbbk}^{P(t),\mathcal L,\mathrm{ss}}$ and $R^{\mathrm{s}}/\mathrm{SL}(V)$ corepresents $\mathbf{M}_{B/\Bbbk}^{P(t),\mathcal L,\mathrm{s}}$. Moreover $R^{\mathrm{s}}/\mathrm{SL}(V)$ is a coarse moduli space of $\mathbf{M}_{B/\Bbbk}^{P(t),\mathcal L,\mathrm{s}}$. 
					\end{itemize}
				\end{theorem}
				
				\begin{remark}
					The coarse moduli space of $\mathbf{M}_{B/\Bbbk}^{P(t),\mathcal L,\mathrm{s}}$ has a natural compactification, which is the corepresentative of $\mathbf{M}_{B/\Bbbk}^{P(t),\mathcal L,\mathrm{ss}}$, which is independent of the GIT problem. 
				\end{remark}
				
				\begin{remark}
					If there exists a semistable sheaf $\mathcal F$ of Hilbert polynomial $P(t)$ that is not polystable, then $\mathbf{M}_{B/\Bbbk}^{P(t),\mathcal L,\mathrm{ss}}$ does not have a coarse moduli space, since $\eta:\mathbf{M}_{B/\Bbbk}^{P(t),\mathcal L,\mathrm{ss}}\to h_{\overline{R}/\!/\mathrm{SL}(V)}$ is not bijective on $\mathrm{Spec}(\Bbbk)$. 
				\end{remark}
			}
		}
		
		\subsubsection{Moduli of unstable sheaves}\label{subsubsection of moduli of unstable sheaves}
		{	
			Let $\tau$ be a Harder-Narasimhan type of length 2. Let $P_i(t)$ for $1\leq i\leq 2$ be the Hilbert polynomials associated to $\tau$ in the sense that if $\mathcal E$ is a pure sheaf of type $\tau$ on $B$, then $P_i(t)$ is the Hilbert polynomial of $\mathrm{HN}_i(\mathcal E)/\mathrm{HN}_{i-1}(\mathcal E)$ for $1\leq i\leq 2$. We have that $P(t):=P_1(t)+P_2(t)$ is the Hilbert polynomial of $\mathcal E$. Let $V_i:=\Bbbk^{P_i(m)}$ for $1\leq i\leq 2$ and let $V:=V_1\oplus V_2$. Denote $V_{\leq i}:=V_1\oplus\cdots\oplus V_i$. Let $\mathrm{P}_\tau\subseteq\mathrm{SL}(V)$ be the parabolic subgroup of $g\in\mathrm{SL}(V)$ stabilising the partial flag $0=V_{\leq 0}\subseteq V_{\leq 1}\subseteq V_{\leq 2}=V$. 
			
			It turns out that it is useful to refine the classification problem according to the dimension of the space of unipotent endomorphisms. In \cite{JacksonJoshua2021MSoU} Definition 7.12 and \cite{HoskinsVictoria2021QbPG} Definition 4.18, the so-called \emph{$\tau$-stability} of pure sheaves was defined. In \cite{HoskinsVictoria2021QbPG} Definition 6.1, the so-called \emph{refined Harder-Narasimhan type} was defined, which encodes information on the dimension of the space of unipotent endomorphisms. 
			
			Theorem 1.1 in \cite{SimpsonCarlosT.1994Moro} shows that sheaves of fixed Harder-Narasimhan type are bounded. In the proof of \cite{HoskinsVictoria2012Qous} Proposition 6.8, a locally closed subscheme $Y_{\tau}^{\mathrm{ss}}\subseteq \mathrm{Quot}_{V(-m)/B/\Bbbk}^{P(t),\mathcal L}$ was described as a parameter space of sheaves of type $\tau$. Closed points in $Y_\tau^{\mathrm{ss}}$ are exactly isomorphism classes $[(\mathcal E,q)]\in\mathrm{Quot}_{V(-m)/B/\Bbbk}^{P(t),\mathcal L}$ such that each $\mathcal E_i$ is semistable with Hilbert polynomial $P_i(t)$ for $1\leq i\leq 2$, where $\mathcal E_i:=q\big(V_{\leq i}(-m)\big)\big/q\big(V_{\leq i-1}(-m)\big)$. In the same manner, one can define locally closed subschemes $Y_\tau^{\mathrm{s}}$ and $Y_\tau$ such that $Y_\tau^{\mathrm{s}}\subseteq Y_\tau^{\mathrm{ss}}\subseteq Y_\tau\subseteq\mathrm{Quot}_{V(-m)/B/\Bbbk}^{P(t),\mathcal L}$. Closed points in $Y_\tau\subseteq\mathrm{Quot}_{V(-m)/B/\Bbbk}^{P(t),\mathcal L}$ are classes $[(\mathcal E,q)]$ such that $\mathcal E_i$ has Hilbert polynomial $P_i(t)$, and closed points in $Y_\tau^{\mathrm{s}}$ are classes $[(\mathcal E,q)]\in Y_\tau^{\mathrm{ss}}$ such that $\mathcal E_i$ is stable. The locally closed subsets $Y_\tau^{\mathrm{(s)s}}\subseteq Y_\tau$ were studied in \cite{HoskinsVictoria2012Qous}, \cite{JacksonJoshuaJames2018Msou}, \cite{JacksonJoshua2021MSoU} and \cite{HoskinsVictoria2021QbPG}, where moduli spaces of \emph{$\tau$-stable} sheaves were constructed as parabolic GIT quotients of $Y_\tau^{\mathrm{s}}$. 
			
			We will recover the locally closed subschemes $Y_\tau^{\mathrm{(s)s}}\subseteq Y_\tau\subseteq \mathrm{Quot}_{V(-m)/B/\Bbbk}^{P(t),\mathcal L}$ and actions of $\mathrm{P}_\tau$ on them from a more functorial perspective in this subsection. From this new perspective, we do not only see closed points but also view $T$-valued points of $Y^{\mathrm{(s)s}}\subseteq Y_\tau$ as families over $T$ of quotients.

			For $Y_\tau$, we describe the associated functor of points $\mathbf Y_\tau\subseteq\mathbf{Quot}_{V(-m)/B/\Bbbk}^{P(t),\mathcal L}$ as a \emph{sub-functor represented by locally closed subschemes}. Let $T\in\mathrm{Sch}/\Bbbk$ and let $[(\mathcal E,q)]\in\mathbf{Quot}_{V(-m)/B/\Bbbk}^{P(t),\mathcal L}(T)$ be an isomorphism class of a quotient 
			\begin{equation}
				\begin{tikzcd}
					V(-m)_T\ar[r,"q"]&\mathcal E\ar[r]&0. 
				\end{tikzcd}
			\end{equation}
			Define the following notation for $1\leq i\leq 2$
			\begin{equation}\label{equation: notations of E_i, E_<=i}
				\begin{split}
					\mathcal E_{\leq i}:=&q\big(V_{\leq i}(-m)_T\big)\subseteq \mathcal E\\
					\mathcal E_i:=&\mathcal E_{\leq i}/\mathcal E_{\leq i-1}. 
				\end{split}
			\end{equation}
			There is a morphism $q_i:V_i(-m)_T\to \mathcal E_i$ such that the following diagram commutes
			\begin{equation}\label{equation: definition of q_i: V_i(-m)_T --> E_i}
				\begin{tikzcd}
					0\ar[r]&V_{\leq i-1}(-m)_T\ar[r]\ar[d,"q_{\leq i-1}"]&V_{\leq i}(-m)_T\ar[r]\ar[d,"q_{\leq i}"]&V_i(-m)_T\ar[r]\ar[d,dashed,"q_i"]&0\\
					0\ar[r]&\mathcal E_{\leq i-1}\ar[r]&\mathcal E_{\leq i}\ar[r]&\mathcal E_i\ar[r]&0
				\end{tikzcd}
			\end{equation}
			where $q_{\leq i}:=q|_{V_{\leq i}(-m)_T}$ and then the Five lemma shows the surjectivity of $q_i$. 

			The subset $\mathbf{Y}_\tau(T)$ consists of $[(\mathcal E,q)]$ such that each $\mathcal E_i$ is flat over $T$ with Hilbert polynomial $P_i(t)$ for $1\leq i\leq 2$. By \href{https://stacks.math.columbia.edu/tag/00HM}{Lemma 00HM} and the fact that $\mathcal E$ is flat over $T$, we have that each $\mathcal E_{\leq i}$ is flat over $T$ if each $\mathcal E_i$ is flat over $T$. Since flatness and Hilbert polynomials are preserved under pullbacks, we have that $\mathbf{Y}_\tau$ defines a sub-functor of $\mathbf{Quot}_{V(-m)/B/\Bbbk}^{P(t),\mathcal L}$.

			The following theorem: \emph{flattening stratification} is important in the construction of Quot schemes. Here we apply it to show that $\mathbf{Y}_\tau\subseteq \mathbf{Quot}_{V(-m)/B/\Bbbk}^{P(t),\mathcal L}$ is represented by locally closed subschemes. 
			
			{
				\begin{theorem}[Theorem 2.1.5, \cite{HuybrechtsDaniel2010Tgom}]\label{theorem: flattening stratification}
					Let $X\to T$ be a projective morphism of schemes of finite type over $\Bbbk$. Let $\mathcal L$ be a line bundle on $X$ which is ample relative to $X\to T$. Let $\mathcal F$ be a coherent sheaf on $X$. Then: 
					\begin{itemize}
						\item The set of Hilbert polynomials is finite 
						\begin{equation}
							\mathcal P:=\bigg\{P(t)\in\mathbb Q[t]:\begin{matrix}P(t)\textrm{ is the Hilbert polynomial of }\mathcal F\\\textrm{at some closed point of }T\end{matrix}\bigg\}; 
						\end{equation}
						\item For each $P(t)\in\mathcal P$ the following functor is represented by a locally closed subscheme $T_P\subseteq T$ 
						\begin{equation}
							\mathbf{Flat}_{\mathcal F/X/T}^{P(t),\mathcal L}:(\mathrm{Sch}/T)^{\mathrm{op}}\to\mathrm{Set},\quad T'\mapsto \begin{cases}\{\mathrm{pt}\},&\begin{matrix}\mathcal F_{T'}\textrm{ is flat over }T'\\\textrm{ with Hilbert polynomial }P(t)\end{matrix}\\
							\emptyset,&\textrm{otherwise}\end{cases}
						\end{equation}
						where $\mathrm{Sch}/T$ is the category of schemes of finite type over $T$, and $\mathcal F_{T'}$ is the pullback of $\mathcal F$ along $T'\to T$; 
						\item The morphism $\coprod_{P\in\mathcal P}T_P\to T$ is bijective on points. 
					\end{itemize}
				\end{theorem}
			}
			
			{
				\begin{lemma}\label{lemma: representability of Y_tau subset Quot}
					$\mathbf Y_\tau\subseteq\mathbf{Quot}_{V(-m)/B/\Bbbk}^{P(t),\mathcal L}$ is representable by locally closed subschemes. 
				\end{lemma}
				\begin{proof}
					Let $T\in\mathrm{Sch}/\Bbbk$ and consider the following diagram in $\mathrm{PSh}(\mathrm{Sch}/\Bbbk)$
					\begin{equation}
						\begin{tikzcd}
							&h_T\ar[d,"\xi"]\\
							\mathbf{Y}_\tau\ar[r,"\subseteq"]&\mathbf{Quot}_{V(-m)/B/\Bbbk}^{P(t),\mathcal L}
						\end{tikzcd}
					\end{equation}
					where $\xi:h_T\to \mathbf{Quot}_{V(-m)/B/\Bbbk}^{P(t),\mathcal L}$ is identified with a class $[(\mathcal E,q)]\in \mathbf{Quot}_{V(-m)/B/\Bbbk}^{P(t),\mathcal L}(T)$ by the Yoneda lemma. 
					
					Recall $\mathcal E_{\leq i}:=q\big(V_{\leq i}(-m)_T\big)$ and $\mathcal E_i:=\mathcal E_{\leq i}/\mathcal E_{\leq i-1}$. For each $1\leq i\leq 2$ consider the functor $\mathbf{Flat}_{\mathcal E_i/B_T/T}^{P_i(t),\mathcal L_T}$ defined in Theorem \ref{theorem: flattening stratification}, which is represented by a locally closed subscheme $T_{\xi,P_i}\subseteq T$. Let 
					\begin{equation}
						T_\xi:=T_{\xi,P_1}\times_T T_{\xi,P_2}
					\end{equation}
					which is a locally closed subscheme of $T$. 
					
					For $T'\to T$ in $\mathrm{Sch}/\Bbbk$, denote $\mathcal E|_{T'}$ the pullback of $\mathcal E$ along $B_{T'}\to B_T$, and similar for $\mathcal E_i,\mathcal E_{\leq i}$. We have
					\begin{equation}
						\begin{split}
							&T'\to T \textrm{ factors through }T'\to T_\xi\to T\\
							\textrm{if and only if}\quad&T'\to T\textrm{ factors though }T_{\xi,P_i}\textrm{ for each }1\leq i\leq 2\\
							\textrm{if and only if}\quad&\begin{cases}\mathcal E_i|_{T'}\textrm{ is flat over }T'\\ \textrm{with Hilbert polynomial }P_i(t)\end{cases}\textrm{ for }1\leq i\leq 2\\
							\textrm{if and only if}\quad&(T'\to T)^*\xi\in \mathbf{Y}_\tau (T'). 
						\end{split}
					\end{equation}
					The above equivalences show the following diagram is Cartesian 
					\begin{equation}
						\begin{tikzcd}
							h_{T_\xi}\ar[r,"\subseteq"]\ar[d]&h_T\ar[d,"\xi"]\\
							\mathbf{Y}_\tau\ar[r,"\subseteq"]&\mathbf{Quot}_{V(-m)/B/\Bbbk}^{P(t),\mathcal L}
						\end{tikzcd}
					\end{equation}
					which proves the representability of $\mathbf{Y}_\tau\subseteq\mathbf{Quot}_{V(-m)/B/\Bbbk}^{P(t),\mathcal L}$ by locally closed subschemes. 
				\end{proof}
			}

			Since $\mathbf{Quot}_{V(-m)/B/\Bbbk}^{P(t),\mathcal L}$ is represented by $\mathrm{Quot}_{V(-m)/B/\Bbbk}^{P(t),\mathcal L}$, by Lemma \ref{lemma: representability of Y_tau subset Quot}, there is a locally closed subscheme $Y_\tau\subseteq \mathrm{Quot}_{V(-m)/B/\Bbbk}^{P(t),\mathcal L}$ representing $\mathbf Y_\tau$. 

			Consider $\prod_{i=1}^2\mathrm{Quot}_{V_i(-m)/B/\Bbbk}^{P_i(t),\mathcal L}$, where the fibre product is over $\Bbbk$. There is a morphism $\prod_{i=1}^2\mathrm{Quot}_{V_i(-m)/B/\Bbbk}^{P_i(t),\mathcal L}\to \mathrm{Quot}_{V(-m)/B/\Bbbk}^{P(t),\mathcal L}$. From their functors of points, we easily see that this morphism is a monomorphism, and it is proper by \href{https://stacks.math.columbia.edu/tag/01W6}{Lemma 01W6}, thus a closed immersion by \href{https://stacks.math.columbia.edu/tag/04XV}{Lemma 04XV}. The closed immersion $\prod_{i=1}^2\mathrm{Quot}_{V_i(-m)/B/\Bbbk}^{P_i(t),\mathcal L}\hookrightarrow\mathrm{Quot}_{V(-m)/B/\Bbbk}^{P(t),\mathcal L}$ factors through $Y_\tau$, which is easily seen from the associated functors of points. There is a retraction $p_\tau:Y_\tau\to \prod_{i=1}^2\mathrm{Quot}_{V_i(-m)/B/\Bbbk}^{P_i(t),\mathcal L}$ such that the following diagram commutes
			\begin{equation}\label{equation: diagram of retraction p_tau on Y_tau}
				\begin{tikzcd}
					\prod_{i=1}^2\mathrm{Quot}_{V_i(-m)/B/\Bbbk}^{P_i(t),\mathcal L}\ar[r,"\subseteq"]\ar[rd,"="]&Y_\tau\ar[d,"p_\tau"]\\&\prod_{i=1}^2\mathrm{Quot}_{V_i(-m)/B/\Bbbk}^{P_i(t),\mathcal L}
				\end{tikzcd}
			\end{equation}
			where $p_\tau$ sends an isomorphism class $[(\mathcal E,q)]$ of families over $T\in\mathrm{Sch}/\Bbbk$ to $\Big[\Big(\bigoplus_{i=1}^2 \mathcal E_i,\bigoplus_{i=1}^2q_i\Big)\Big]$, where $\mathcal E_i$ and $q_i$ are defined in equation \eqref{equation: notations of E_i, E_<=i} and equation \eqref{equation: definition of q_i: V_i(-m)_T --> E_i}. 

			Consider $(\mathbb G_m)^2\to\mathrm{GL}(V)$ such that $(t_1,t_2)\mapsto t_11_{V_1}\oplus t_21_{V_2}\in\mathrm{GL}(V)$. Then $(\mathbb G_m)^2$ acts on $Y_\tau$. With this $(\mathbb G_m)^2$-action, we have that $\prod_{i=1}^2\mathrm{Quot}_{V_i(-m)/B/\Bbbk}^{P_i(t),\mathcal L}$ is the \emph{fixed point subscheme} $(Y_\tau)^{(\mathbb G_m)^2}\hookrightarrow Y_\tau$ in the sense of \cite{FogartyJohn1973FPS}. It is easy to see if we consider $\lambda:\mathbb G_m\to(\mathbb G_m)^2$ given by $t\mapsto (t^{w_1},t^{w_2})$ such that $w_1>w_2$, then the fixed point subscheme of the action of $\lambda$ is also $\prod_{i=1}^2\mathrm{Quot}_{V_i(-m)/B/\Bbbk}^{P_i(t),\mathcal L}$
			\begin{equation}
				\prod_{i=1}^2\mathrm{Quot}_{V_i(-m)/B/\Bbbk}^{P_i(t),\mathcal L}\cong (Y_\tau)^{(\mathbb G_m)^2}=(Y_\tau)^{\lambda}. 
			\end{equation}
			
			We then describe the open subschemes $Y_\tau^{\mathrm{s}}\subseteq Y_\tau^{\mathrm{ss}}\subseteq Y_\tau$. For $1\leq i\leq 2$, let $R_i^{\mathrm{(s)s}}\subseteq \mathrm{Quot}_{V_i(-m)/B/\Bbbk}^{P_i(t),\mathcal L}$ be the open subscheme whose closed points are $[(\mathcal E_i,q_i)]$ such that $\mathcal E_i$ is a (semi)stable pure sheaf on $B$ and $q:V_i(-m)\to\mathcal E_i$ induces an isomorphism $V_i\cong H^0(B,\mathcal E_i(m))$ (cf. Theorem \ref{theorem: simpson's construction of moduli of semistable sheaves}). We define $Y_\tau^{\mathrm{(s)s}}$ as the fibre product indicated in the diagram below 
			\begin{equation}\label{equation: defining diagram of Y_tau^(s)s}
				\begin{tikzcd}
					Y_\tau^{\mathrm{(s)s}}\ar[rd,phantom,very near start,"\lrcorner"]\ar[r,"\subseteq"]\ar[d,"p_\tau^{\mathrm{(s)s}}"]&Y_\tau\ar[d,"p_\tau"]\\
					\prod_{i=1}^2R_i^{\mathrm{(s)s}}\ar[r,"\subseteq"]&\prod_{i=1}^2\mathrm{Quot}_{V_i(-m)/B/\Bbbk}^{P_i(t),\mathcal L}. 
				\end{tikzcd}
			\end{equation}
			We have that $Y_\tau^{(s)s}\subseteq Y_\tau$ is open. So $p_\tau^{\mathrm{(s)s}}:Y_\tau^{\mathrm{(s)s}}\to \prod_{i=1}^2R_i^{\mathrm{(s)s}}$ are also retractions as base changes of $p_\tau$. In particular, $p_\tau$ and $p_\tau^{\mathrm{(s)s}}$ are surjective. 

			There is a $\mathrm P_\tau$-action on $Y_\tau$ induced from the $\mathrm{SL}(V)$-action on $\mathrm{Quot}_{V(-m)/B/\Bbbk}^{P(t),\mathcal L}$. Similarly there are $\mathrm P_\tau$-actions on $Y_\tau^{\mathrm{s}}\subseteq Y_\tau^{\mathrm{ss}}$. These $\mathrm P_\tau$-actions are linearised with respect to the restrictions of $\mathcal M$ on $\mathrm{Quot}_{V(-m)/B/\Bbbk}^{P(t),\mathcal L}$, where $\mathcal M$ is the ample line bundle described in \eqref{equation: ample line bundle M on Quot}. 
			
			In \cite{JacksonJoshuaJames2018Msou} and \cite{JacksonJoshua2021MSoU}, Jackson considered the induced reduced subscheme structure on $Y_\tau^{\mathrm{s}}$ for a Harder-Narasimhan type $\tau$ of length 2, and stratify this reduced subscheme further according to the dimension of the stabiliser in the unipotent radical of $\mathrm P_\tau$. Non-reductive GIT could be applied to each stratum and the following theorem was proved. 

			{
				\begin{definition}[cf. \cite{JacksonJoshua2021MSoU} Definition 7.1.2, \cite{HoskinsVictoria2021QbPG} Definition 6.1]\label{definition: Jackson's (tau,alpha)-stability}
					Let $B$ be a connected projective scheme over $\Bbbk$ with an ample line bundle $\mathcal L$. Let $\tau$ be a Harder-Narasimhan type of length 2. Let $\delta\in \mathbb N$. A pure sheaf $\mathcal E$ on $B$ is called \emph{$(\tau,\delta)$-stable} if: 
					\begin{itemize}
						\item $\mathcal E$ has Harder-Narasimhan type $\tau$; 
						\item if $\{\mathcal E_{\leq i}\}_{i=0}^2$ is the Harder-Narasimhan filtration, then $\mathcal E_i:=\mathcal E_{\leq i}/\mathcal E_{\leq i-1}$ is stable for $i=1,2$; 
						\item $\mathcal E\not\cong \mathcal E_1\oplus\mathcal E_2$; 
						\item $\dim \mathrm{Hom}_{\mathcal O_B}(\mathcal E_2,\mathcal E_1)=\delta$. 
					\end{itemize}
				\end{definition}

				\begin{theorem}[\cite{JacksonJoshua2021MSoU} Theorem 7.1.3]\label{theorem: Jackson's moduli of length 2}
					Let $B$ be a connected projective scheme over $\Bbbk$ with an ample line bundle $\mathcal L$. Let $\tau$ be a Harder-Narasimhan type of length 2 and let $\delta\in \mathbb N$. Then there exists a quasi-projective moduli space $M_{\tau,\delta}$ of $(\tau,\delta)$-stable sheaves on $B$, in the sense that closed points in $M_{\tau,\delta}$ are exactly isomorphism classes of $(\tau,\delta)$-stable sheaves. 
				\end{theorem}
				
				\begin{remark}
					The moduli spaces $M_{\tau,\delta}$ above are reduced schemes since they were constructed by the classic U-hat theorem (cf. Theorem 2.20, \cite{BércziGergely2016Pcog}), which only produces reduced quotient schemes. 
				\end{remark}
			}

			The non-reductive GIT problem considered in \cite{JacksonJoshua2021MSoU} to prove Theorem \ref{theorem: Jackson's moduli of length 2} was for the linearisation 
			\begin{equation}
				\mathrm P_\tau\curvearrowright \big(\overline{Y_\tau^{\mathrm{s}}},\;\mathcal M\big|_{\overline{Y_\tau^{\mathrm{s}}}}\big) 
			\end{equation}
			where $\overline{Y_\tau^{\mathrm{s}}}\subseteq \mathrm{Quot}_{V(-m)/B/\Bbbk}^{P(t),\mathcal L}$ is the scheme theoretic closure of $Y_\tau^{\mathrm{s}}$ in $\mathrm{Quot}_{V(-m)/B/\Bbbk}^{P(t),\mathcal L}$. In Appendix \ref{appendix: on the boundaries of Y_tau and Y_tau^s}, we study the boundary of $Y_\tau^{\mathrm{s}}\subseteq\overline{Y_\tau^{\mathrm{s}}}$ and Proposition \ref{proposition: Y_tau^s closure and its Z_min^K-s} shows that when applying our main theorem in non-reductive GIT (Theorem \ref{theorem: geometric U and U-hat quotients of strata}). 

			For each $\delta\in\mathbb N$, there exists a locally closed subset $Y_{\tau,\delta}^{\mathrm{s}}$ (See equation \eqref{equation: definition of Y_(tau,delta)^s} for its definition), whose closed points are $[(\mathcal E,q)]$ such that $\mathcal E$ is $(\tau,\alpha)$-stable. Each of $Y_{\tau,\delta}^{\mathrm{s}}$ has a quasi-projective universal geometric quotient by $\mathrm P_\tau$
			\begin{equation}
				Y_{\tau,\delta}^{\mathrm{s}}\to Y_{\tau,\delta}^{\mathrm{s}}/\mathrm P_\tau. 
			\end{equation}
			Note that in \cite{JacksonJoshua2021MSoU}, only closed points of $Y_{\tau,\delta}^{\mathrm{s}}\subseteq Y_\tau^{\mathrm{s}}\subseteq \overline{Y_\tau^{\mathrm{s}}}$ were described, and scheme structures were implicitly the reduced ones. However by \cite{https://doi.org/10.48550/arxiv.0911.2301} Remark 7.8, even for rank 2 unstable bundles over curves, the moduli space need not be reduced. 

			{
				One key fact peculiar to length two is the following, which fails for higher lengths. 
				\begin{lemma}[\cite{JacksonJoshua2021MSoU} Corollary 7.2.6]\label{lemma: dim Stab_U constant along p_tau for length 2 sheaves}
					For any closed point $x\in Y_\tau^{\mathrm{ss}}$, we have 
					\begin{equation}
						\dim\mathrm{Stab}_{\mathrm U_\tau}(x)=\dim\mathrm{Stab}_{\mathrm U_\tau}(p_\tau(x))
					\end{equation}
					where $p_\tau:Y_\tau\to \prod_{i=1}^2\mathrm{Quot}_{V(-m)/B/\Bbbk}^{P_i(t),\mathcal L}$ is the retraction and $\mathrm U_\tau\subseteq\mathrm P_\tau$ is the unipotent radical. 
				\end{lemma}
			}

			We will consider the non-reduced version of the GIT problem above. This requires: 
			\begin{itemize}
				\item an extension of non-reductive GIT developed in \cite{BércziGergely2016Pcog}, \cite{BércziGergely2018Gitf} to linear actions on non-reduced schemes; 
				\item a scheme theoretic description of $Y_{\tau,\delta}^{\mathrm{s}}$. 
			\end{itemize}
			We will address these in the following two sections respectively. 
		}
	}
}

\section{Non-reductive GIT for graded additive groups}
{

	Let $U$ be an additive group, i.e. $U\cong(\mathbb G_a)^r$ for $r:=\dim U$. Let $\lambda:\mathbb G_m\to \mathrm{Aut}(U)$ be a 1-parameter subgroup acting on $U$ by conjugation such that the induced representation of $\mathbb G_m$ on $\mathfrak u:=\mathrm{Lie}(U)$ has a unique positive weight. Let $\hat U:=U\rtimes_\lambda\mathbb G_m$ be the semi-direct product of $U$ and the 1PS $\lambda$. We will consider group actions by $\hat U$. 


	\subsection{Some definitions for group actions}
	{
		All schemes are Noetherian over $\Bbbk$. Let $G$ be an affine algebraic group. Let $X$ be a scheme acted on by $G$ via the action morphism $\sigma:G\times X\to X$. Let $\mathfrak g^*$ be the dual of the Lie algebra of $G$. 

		\subsubsection{Infinitesimal stabilisers of scheme-valued points}
		{
			{
				\begin{definition}\label{definition: infinitesimal stabiliser of T-points}
					Let $\sigma:G\times X\to X$ be an action. Let $t:T\to X$ be a morphism in $\mathrm{Sch}/\Bbbk$. The \emph{stabiliser} of the $T$-point $t$, denoted by $\mathrm{Stab}_G(t)$ or $\mathrm{Stab}_G(T)$, is defined as the following fibre product 
					\begin{equation}
						\begin{tikzcd}
							\mathrm{Stab}_G(t)\ar[rrr]\ar[d]\ar[rrrd,phantom,very near start,"\lrcorner"]&&&T\ar[d,"{(t,1_T)}"]\\
							G\times T\ar[rrr,"{\big(\sigma\circ(\mathrm{pr}_0,t\circ\mathrm{pr}_1),\;\mathrm{pr}_1\big)}"]&&&X\times T. 
						\end{tikzcd}
					\end{equation}
					The \emph{infinitesimal stabiliser} of the $T$-point $t$, denoted by $\mathrm{Stab}_{\mathfrak g}(t)$ or $\mathrm{Stab}_{\mathfrak g}(T)$, is defined as the normal sheaf of the closed immersion $e'':T\hookrightarrow \mathrm{Stab}_G(t)$
					\begin{equation}
						\begin{tikzcd}
							T\ar[rrrrd,"1_T",bend left =15]\ar[rd,"e''"]\ar[rdd,swap,bend right =20,"e'"]\\
							&\mathrm{Stab}_G(t)\ar[rrr]\ar[d]\ar[rrrd,phantom,very near start,"\lrcorner"]&&&T\ar[d,"{(t,1_T)}"]\\
							&G\times T\ar[rrr,"{\big(\sigma\circ(\mathrm{pr}_0,t\circ\mathrm{pr}_1),\;\mathrm{pr}_1\big)}"]&&&X\times T
						\end{tikzcd}
					\end{equation}
					induced by the closed immersion $e':T\hookrightarrow G\times T$ at the identity $e:\mathrm{Spec}(\Bbbk)\hookrightarrow G$
					\begin{equation}
						\mathrm{Stab}_{\mathfrak g}(t):=\mathcal N_{e''}\equiv\mathcal Hom_{\mathcal O_T}\big(\mathcal C_{T/\mathrm{Stab}_G(t)},\mathcal O_T\big). 
					\end{equation}
				\end{definition}
			}

			Consider the following diagram for a morphism $t:T\to X$ in $\mathrm{Sch}/\Bbbk$ 
			\begin{equation}
				\begin{tikzcd}
					T\ar[r,"e''"]\ar[d,equal]&\mathrm{Stab}_G(t)\ar[d]\\
					T\ar[r,"e'"]\ar[d,equal]&G\times T\ar[d,"{\big(\sigma\circ(\mathrm{pr_G},t\circ\mathrm{pr}_T),\mathrm{pr}_T\big)}"]\\
					T\ar[r,"{(t,1_T)}"]&X\times T
				\end{tikzcd}
			\end{equation}
			where horizontal morphisms are all immersions. By \href{https://stacks.math.columbia.edu/tag/01R3}{Lemma 01R3}, there is a chain of conormal sheaves on $T$
			\begin{equation}
				\begin{tikzcd}
					\mathcal C_{(t,1_T)}\ar[r]&\mathfrak g^*\otimes_\Bbbk\mathcal O_T\ar[r]&\mathcal C_{T/\mathrm{Stab}_G(t)}\ar[r]&0. 
				\end{tikzcd}
			\end{equation}
			which is exact since $\mathrm{Stab}_G(t)\to G\times T$ is an immersion and $\mathrm{Stab}_G(t)\cong (G\times T)\times_{X\times T}T$. The following diagram is Cartesian 
			\begin{equation}
				\begin{tikzcd}
					T\ar[rd,phantom,very near start,"\lrcorner"]\ar[r,"{(t,1_T)}"]\ar[d,"t"]&X\times T\ar[d,"1_X\times t"]\\
					X\ar[r,"\delta"]&X\times X
				\end{tikzcd}
			\end{equation}
			where $\delta:X\to X\times X$ is the diagonal morphism, and horizontal morphisms are immersions. By \href{https://stacks.math.columbia.edu/tag/01R4}{Lemma 01R4}, there is a surjective morphism of conormal sheaves 
			\begin{equation}
				\begin{tikzcd}
					t^*\Omega_{X/\Bbbk}\ar[r]&\mathcal C_{(t,1_T)}\ar[r]&0. 
				\end{tikzcd}
			\end{equation}
			Then the following sequence is exact 
			\begin{equation}\label{equation: cokernel sequence of t^*(phi)}
				\begin{tikzcd}
					t^*\Omega_{X/\Bbbk}\ar[r]&\mathfrak g^*\otimes_\Bbbk\mathcal O_T\ar[r]&\mathcal C_{T/\mathrm{Stab}_G(t)}\ar[r]&0. 
				\end{tikzcd}
			\end{equation}

			{
				\begin{definition}\label{definition: infinitesimal action of g on X}
					Let $\sigma:G\times X\to X$ be an action. Choose $t:T\to X$ to be the identity morphism $1_X:X\to X$. We call $\phi:\Omega_{X/\Bbbk}\to \mathfrak g^*\otimes_\Bbbk\mathcal O_X$ constructed above the \emph{infinitesimal action of $\mathfrak g$ on $X$}. 
				\end{definition}
			}

			{
				\begin{lemma}\label{lemma: coker of infinitesimal stabiliser restricts to invariant subschemes}
					Let $f:X_1\to X_2$ be an immersion of schemes. Let $G$ act on $X_1,X_2$ such that $f$ is $G$-invariant. Let $\phi_i:\Omega_{X_i/\Bbbk}\to\mathfrak g^*\otimes_\Bbbk\mathcal O_{X_i}$ be the infinitesimal action of $\mathfrak g$ on $X_i$ for $i=1,2$. Then we have $\mathrm{coker}(\phi_1)\cong f^*\mathrm{coker}(\phi_2)$ fitting into the diagram
					\begin{equation}
						\begin{tikzcd}
							f^*\Omega_{X_2/\Bbbk}\ar[r,"f^*\phi_2"]\ar[d]&\mathfrak g^*\otimes_\Bbbk f^*\mathcal O_{X_2}\ar[r]\ar[d,"\cong"]&f^*\mathrm{coker}(\phi_2)\ar[r]\ar[d,"\cong"]&0\\
							\Omega_{X_1/\Bbbk}\ar[r,"\phi_1"]&\mathfrak g^*\otimes_\Bbbk\mathcal O_{X_1}\ar[r]&\mathrm{coker}(\phi_1)\ar[r]&0. 
						\end{tikzcd}
					\end{equation}
				\end{lemma}
				\begin{proof}
					The morphism $f$ is $G$-invariant, so the left square commutes. The left square induces a morphism $f^*\mathrm{coker}(\phi_2)\to \mathrm{coker}(\phi_1)$. Since $f$ is an immersion, we have that $f^*\Omega_{X_2/\Bbbk}\to \Omega_{X_1/\Bbbk}$ is surjective. Therefore $f^*\mathrm{coker}(\phi_2)\to \mathrm{coker}(\phi_1)$ is an isomorphism by the Five lemma. 
				\end{proof}
			}

			{
				\begin{lemma}\label{lemma: Stab and base change}
					Let $G$ act on $X$. In the category $\mathrm{Sch}/X$, we have 
					\begin{itemize}
						\item[(1)] The formation of $\mathrm{Stab}_G(t)$ commutes with arbitrary base change, in the sense that if $T'\to T$ is a morphism in $\mathrm{Sch}/X$ between $t':T'\to X$ and $t:T\to X$, then the following are Cartesian 
						\begin{equation}
							\begin{tikzcd}
								\mathrm{Stab}_G(t')\ar[r]\ar[d]\ar[rd,phantom,very near start, "\lrcorner"]&T'\ar[d]\\
								\mathrm{Stab}_G(t)\ar[r]&T
							\end{tikzcd},\quad
							\begin{tikzcd}
								\mathrm{Stab}_G(t')\ar[r]\ar[d]\ar[rd,phantom,very near start, "\lrcorner"]&\mathrm{Stab}_G(t)\ar[d]\\
								G\times T'\ar[r]&G\times T
							\end{tikzcd},\quad 
							\begin{tikzcd}
								T'\ar[r]\ar[d]\ar[rd,phantom,very near start,"\lrcorner"]&\mathrm{Stab}_G(t')\ar[d]\\
								T\ar[r]&\mathrm{Stab}_G(t). 
							\end{tikzcd}
						\end{equation}

						\item[(2)] The formation of $\mathrm{Stab}_{\mathfrak g}(t)$ commutes with flat base change, in the sense that if $T'\to T$ is a flat morphism in $\mathrm{Sch}/X$ between $t':T'\to X$ and $t:T\to X$, then 
						\begin{equation}
							(T'\to T)^*\mathrm{Stab}_{\mathfrak g}(t)\cong \mathrm{Stab}_{\mathfrak g}(t'). 
						\end{equation}

						\item[(3)] The formation of sequence \eqref{equation: cokernel sequence of t^*(phi)} commutes with arbitrary base change in the sense that if $T'\to T$ is a morphism in $\mathrm{Sch}/X$ between $t':T'\to X$ and $t:T\to X$, then the sequence 
						\begin{equation}
							\begin{tikzcd}
								t^*\Omega_{X/\Bbbk}\ar[r]&\mathfrak g^*\otimes_\Bbbk\mathcal O_T\ar[r]&\mathcal C_{T/\mathrm{Stab}_G(t)}\ar[r]&0
							\end{tikzcd}
						\end{equation}
						pulls back along $T'\to T$ to the sequence 
						\begin{equation}
							\begin{tikzcd}
								(t')^*\Omega_{X/\Bbbk}\ar[r]&\mathfrak g^*\otimes_\Bbbk\mathcal O_{T'}\ar[r]&\mathcal C_{T'/\mathrm{Stab}_G(t')}\ar[r]&0. 
							\end{tikzcd}
						\end{equation}
						Therefore sequence \eqref{equation: cokernel sequence of t^*(phi)} is the cokernel sequence of $t^*\phi$ for $\phi:\Omega_{X/\Bbbk}\to \mathfrak g^*\otimes_\Bbbk\mathcal O_X$. 
					\end{itemize}
				\end{lemma}
				\begin{proof}
					It is routine to check $(1)$ and $(3)$. Note that in $(2)$ we need flatness since $\mathrm{Stab}_{\mathfrak g}(t)=\mathcal Hom_{\mathcal O_T}\big(\mathcal C_{T/\mathrm{Stab}_G(t)},\mathcal O_T\big)$ and the formation of dual sheaves commutes with flat base change by \href{https://stacks.math.columbia.edu/tag/0C6I}{Lemma 0C6I}. 
				\end{proof}

				\begin{lemma}\label{lemma: Stab_g(t') criterion for local freedom of coker(phi)|_T}
					Let $G$ act on $X$. Let $\phi:\Omega_{X/\Bbbk}\to \mathfrak g^*\otimes_\Bbbk\mathcal O_X$ be the infinitesimal action. Let $t:T\to X$ be a $T$-point. Let $\delta\in\mathbb N$. The following are equivalent: 
					\begin{itemize}
						\item[(1)] $t^*\mathrm{coker}(\phi)$ is locally free of rank $\delta$ on $T$; 
						\item[(2)] $\mathrm{Stab}_{\mathfrak g}(t')$ is locally free of rank $\delta$ on $T'$ for all $T'\to T$, where $t'$ is the composition $t':T'\to T\xrightarrow{t}X$; 
						\item[(3)] the same in $(2)$ with $T'\to T$ ranging over all closed subschemes of $T$. 
					\end{itemize}
				\end{lemma}
				\begin{proof}
					Assume $(1)$. Then $(t')^*\mathrm{coker}(\phi)$ is locally free of rank $\delta$. Recall $\mathrm{Stab}_{\mathfrak g}(t')=\mathcal Hom_{\mathcal O_{T'}}\big((t')^*\mathrm{coker}(\phi),\mathcal O_{T'}\big)$, which is locally free of rank $\delta$ on $T'$. This proves $(1)\implies (2)$. 

					Obviously $(2)\implies(3)$. 

					Assume $(3)$. For $k\in\mathbb N$, let $Z_k\subseteq T$ denote the closed subscheme associated to the $k$th Fitting ideal of $t^*\mathrm{coker}(\phi)$. For $S_k:=Z_{k-1}\setminus Z_k$, we have $\mathrm{coker}(\phi)|_{S_k}:=(S_k\to T)^*t^*\mathrm{coker}(\phi)$ is locally free of rank $k$ by \href{https://stacks.math.columbia.edu/tag/05P8}{Lemma 05P8}. Then $\mathrm{Stab}_{\mathfrak g}(S_k)=\mathcal Hom_{\mathcal O_{S_k}}\big(\mathrm{coker}(\phi)|_{S_k},\mathcal O_{S_k}\big)$ is locally free of rank $k$ on $S_k$. Condition $(3)$ implies that $\mathrm{Stab}_{\mathfrak g}(Z_{k-1})$ is locally free of rank $\delta$ on $Z_{k-1}$. Since $S_k\to Z_{k-1}$ is an open immersion, we have $\mathrm{Stab}_{\mathfrak g}(Z_{k-1})|_{S_k}\cong\mathrm{Stab}_{\mathfrak g}(S_k)$ by Lemma \ref{lemma: Stab and base change}. Therefore if $k\ne\delta$, then $S_k=\emptyset$. Then $Z_{\delta-1}=T$ and $Z_\delta=\emptyset$. 

					Condition $(3)$ implies that $\mathrm{Stab}_{\mathfrak g}(T)$ is locally free of rank $\delta$ on $T$. \href{https://stacks.math.columbia.edu/tag/05P8}{Lemma 05P8} implies that $t^*\mathrm{coker}(\phi)$ is locally generated by $\leq \delta$ sections. Let $W\subseteq T$ be an affine open subset such that $\mathrm{Stab}_{\mathfrak g}(T)|_W$ is free of rank $\delta$ and such that there exists an exact sequence of $\mathcal O_W$-modules for some $N\in\mathbb N$
					\begin{equation}
						\begin{tikzcd}
							\mathcal O_W^{\oplus N}\ar[r,"\eta"]&\mathcal O_W^{\oplus\delta}\ar[r,"\zeta"]&\mathrm{coker}(\phi)|_W\ar[r]&0
						\end{tikzcd}
					\end{equation}
					where $\mathrm{coker}(\phi)|_W:=t^*\mathrm{coker}(\phi)|_W$. The sequence of dual $\mathcal O_W$-modules is exact 
					\begin{equation}
						\begin{tikzcd}
							0\ar[r]&\mathrm{Stab}_{\mathfrak g}(W)\ar[r,"\zeta^*"]&\mathcal O_W^{\oplus\delta}\ar[r,"\eta^*"]&\mathcal O_W^{\oplus N}. 
						\end{tikzcd}
					\end{equation}
					We have that $\mathrm{Stab}_{\mathfrak g}(W)$ is free of rank $\delta$ since $\mathrm{Stab}_{\mathfrak g}(W)\cong\mathrm{Stab}_{\mathfrak g}(T)|_W$ by Lemma \ref{lemma: Stab and base change}. It is well known that $\det(\zeta^*)\in\mathcal O_W(W)$ is a non-zerodivisor when $\zeta^*$ is an injective map of finite free modules of the same rank $\delta$. Moreover, $\mathrm{coker}(\zeta^*)$ is annihilated by $\det(\zeta^*)$ and $\mathrm{coker}(\zeta^*)\cong\mathrm{im}(\eta^*)\subseteq \mathcal O_W^{\oplus N}$ is contained in a torsion free module. Therefore $\mathrm{coker}(\zeta^*)\cong \mathrm{im}(\eta^*)=0$, i.e. $\eta=0$, which implies that $\zeta:\mathcal O_W^{\oplus\delta}\to \mathrm{coker}(\phi)|_W$ is an isomorphism. This proves $(1)$. 
				\end{proof}
			}

			For the infinitesimal action $\phi:\Omega_{X/\Bbbk}\to \mathfrak g^*\otimes_\Bbbk \mathcal O_X$, let $\mathrm{Fit}_k(\phi)\subseteq \mathcal O_X$ denote the $k$th Fitting ideal of $\mathrm{coker}(\phi)$. 

			{
				\begin{lemma}\label{lemma: equivalence of dimStab>k and Fit_k in m_x}
					Let $G$ act on $X$. Let $\phi:\Omega_{X/\Bbbk}\to \mathfrak g^*\otimes_\Bbbk\mathcal O_X$ be the infinitesimal action. Let $\mathbb K$ be a field. Let $x:\mathrm{Spec}(\mathbb K)\to X$ be a $\mathbb K$-point. Let $k\in\mathbb Z$. The following are equivalent: 
					\begin{itemize}
						\item[(1)] $\dim_{\mathbb K}\mathrm{Stab}_{\mathfrak g}(x)>k$; 
						\item[(2)] $\mathrm{Fit}_k(\phi)\subseteq \ker(x^\sharp )$, where $x^\sharp :\mathcal O_X\to x_*(\mathbb K)$ is the morphism of sheaves of rings associated to $x$. 
					\end{itemize}
				\end{lemma}
				\begin{proof}
					Let $r:=\dim_\Bbbk\mathfrak g$. When $k\geq r$ or $k\leq -1$, both $(1)$ and $(2)$ are false or both are true. We only need to consider when $0\leq k<r$. 

					Let $\mathrm{Spec}(A)\subseteq X$ be an open neighbourhood of $x\in X$. Let $x_A:\mathrm{Spec}(\mathbb K)\to \mathrm{Spec}(A)$ be the morphism when $x$ is viewed as a $\mathbb K$-point of $\mathrm{Spec}(A)$. On $\mathrm{Spec}(A)$, the infinitesimal action $\phi:\Omega_{X/\Bbbk}\to\mathfrak g^*\otimes_\Bbbk\mathcal O_X$ corresponds to a map of $A$-modules 
					\begin{equation}
						\phi_A:\Omega_{A/\Bbbk}\to\mathfrak g^*\otimes_\Bbbk A. 
					\end{equation}
					Choose a surjective map 
					\begin{equation}
						\alpha:A^{\oplus J}\to \Omega_{A/\Bbbk}
					\end{equation}
					for an index set $J$. Choose a basis $(u_1,\cdots,u_r)$ of $\mathfrak g^*$. Let $M=(m_{ij})_{\substack{1\leq i\leq r\\j\in J}}$ be the matrix with entries $m_{ij}\in A$ such that 
					\begin{equation}
						(\phi_A\circ\alpha)(e_j)=\sum_{i=1}^ru_i\otimes m_{ij},\quad \textrm{for all }j\in J. 
					\end{equation}
					Then the ideal $\mathrm{Fit}_k(\phi_A)\subseteq A$ is generated by $(r-k)\times(r-k)$-minors of $M$.

					Let $x_A^\sharp :A\to \mathbb K$ be the ring map corresponding to the point $x_A:\mathrm{Spec}(\mathbb K)\to\mathrm{Spec}(A)$. The base change of $\phi_A\circ\alpha$ along $x_A^\sharp :A\to \mathbb K$ is 
					\begin{equation}
						x_A^*(\phi_A\circ\alpha):\mathbb K^{\oplus J}\to \mathfrak g^*\otimes_\Bbbk\mathbb K
					\end{equation}
					where $\Bbbk\to \mathbb K$ corresponds to $\mathrm{Spec}(\mathbb K)\xrightarrow{x}X\to \mathrm{Spec}(\Bbbk)$. Then the linear map $x_A^*(\phi_A\circ\alpha)$ is represented by the matrix $x_A^\sharp (M):=\big(x_A^\sharp (m_{ij})\big)$ with entries in $\mathbb K$. Note that $\xi\in\mathfrak g\otimes_\Bbbk\mathbb K$ stabilises $x$ if and only if $\mathrm{im}\big(x_A^*(\phi_A\circ\alpha)\big)\subseteq\xi^\perp$. 

					We have 
					\begin{equation}
						\begin{split}
							(1)\quad\iff&\quad \textrm{there exist independent }\xi_0,\cdots,\xi_k\in\mathfrak g\otimes_\Bbbk\mathbb K\textrm{ stabilising }x\\
							\iff&\quad \textrm{there exist independent }\xi_0,\cdots,\xi_k\in\mathfrak g\otimes_\Bbbk\mathbb K\\
							&\qquad\textrm{such that }\mathrm{im}\big(x_A^*(\phi_A\circ\alpha)\big)\subseteq (\xi_0,\cdots,\xi_k)^\perp\\
							\iff&\quad\dim_{\mathbb K} \mathrm{im}\big(x_A^*(\phi_A\circ\alpha)\big)<r-k\\
							\iff&\quad\mathrm{rank}(x_A^\sharp M)<r-k\\
							\iff&\quad \textrm{any }(r-k)\times(r-k)\textrm{-minor of }x_A^\sharp M\textrm{ vanishes}\\
							\iff&\quad \textrm{any }(r-k)\times(r-k)\textrm{-minor of }M\textrm{ is in }\ker(x_A^\sharp )\\
							\iff&\quad \mathrm{Fit}_k(\phi_A)\subseteq\ker(x_A^\sharp )\\
							\iff&\quad (2)
						\end{split}
					\end{equation}
				\end{proof}

				\begin{remark}
					If $\mathrm{coker}(\phi)$ is locally free of rank $k$ on $X$, then $\mathrm{Fit}_{k-1}(\phi)=0$ and $\mathrm{Fit}_k(\phi)=\mathcal O_X$. By Lemma \ref{lemma: equivalence of dimStab>k and Fit_k in m_x}, this implies $\dim_{\mathbb K}\mathrm{Stab}_{\mathfrak g}(x)=k$ for any field $\mathbb K$ and $\mathbb K$-point $x\in X$. Conversely, if $\dim_{\mathbb K}\mathrm{Stab}_{\mathfrak g}(x)=k$ for any $\mathbb K$ and $\mathbb K$-point $x\in X$, then we only have $\mathrm{Fit}_{k-1}(\phi)\subseteq\mathrm{nil}(X)$ and $\mathrm{Fit}_k(\phi)=\mathcal O_X$, where $\mathrm{nil}(X)\subseteq \mathcal O_X$ is the nilpotent radical. 
				\end{remark}
			}
		}

		\subsubsection{Scheme-valued infinitesimal transformations}
		{
			Let $T$ be a scheme over $\Bbbk$. Recall that a $T$-valued point of $X$ is a morphism $t:T\to X$. We will define \emph{$T$-valued infinitesimal transformations in $G$} and when a $T$-valued infinitesimal transformation stabilises a $T$-valued point. 
			
			Let $\Bbbk[\epsilon]:=\Bbbk[t]/t^2$ be the algebra of dual numbers. For $T\in\mathrm{Sch}/\Bbbk$, let $T[\epsilon]:=T\times_{\Bbbk}\mathrm{Spec}(\Bbbk[\epsilon])$. The underlying topological spaces of $T$ and $T[\epsilon]$ are the same. We can think of $T[\epsilon]$ as a scheme on the topological space $|T|$ with the structure sheaf 
			\begin{equation}
				\mathcal O_{T[\epsilon]}:=\mathcal O_T\otimes_\Bbbk \Bbbk[\epsilon]\cong\mathcal O_T\oplus\epsilon\mathcal O_T. 
			\end{equation}
			There is a natural closed immersion $T\hookrightarrow T[\epsilon]$, corresponding to $\mathcal O_{T[\epsilon]}\xrightarrow{\epsilon=0}\mathcal O_T$.

			{
				\begin{definition}\label{definition: T-valued infinitesimal transformation in G}
					Let $G$ be an affine algebraic group over $\Bbbk$. A \emph{$T$-valued infinitesimal transformation in $G$} is a morphism $\rho:T[\epsilon]\to G$ satisfying the diagram 
					\begin{equation}
						\begin{tikzcd}
							T\ar[r,hook]\ar[d]&T[\epsilon]\ar[d,"\rho"]\\
							\mathrm{Spec}(\Bbbk)\ar[r,hook,"e"]&G. 
						\end{tikzcd}
					\end{equation}
				\end{definition}
				
				\begin{definition}\label{definition: T-infinitesimal transformations stabilising T-points}
					Let $G$ act on $X$. Let $t:T\to X$ be a $T$-valued point of $X$, and let $\rho:T[\epsilon]\to G$ be a $T$-valued infinitesimal transformation in $G$. We say that $\rho$ \emph{stabilises} $t$ if the following diagram commutes 
					\begin{equation}
						\begin{tikzcd}
							T[\epsilon]\ar[r,"{(\rho,t\circ\mathrm{pr})}"]\ar[d,"\mathrm{pr}"]&G\times X\ar[d,"\sigma"]\\
							T\ar[r,"t"]&X. 
					\end{tikzcd}
					\end{equation}
				\end{definition}
				
				\begin{remark}\label{remark: vectors in Lie algebra stabilising a k-point}
					If we take $T=\mathrm{Spec}(\Bbbk)$, then $T$-valued points of $X$ are $\Bbbk$-points, and $T$-valued infinitesimal transformations in $G$ are vectors in the Lie algebra $\mathfrak g$. For a $T$-valued point $x\in X$ and a $T$-valued infinitesimal transformation $\xi\in\mathfrak g$, the definition of when $\xi$ stabilises $x$ in Definition \ref{definition: T-infinitesimal transformations stabilising T-points} is equivalent to the requirement that the pullback of $\xi:\Omega_{X/\Bbbk}\to\mathcal O_X$ along $x:\mathrm{Spec}(\Bbbk)\hookrightarrow X$ is zero. 
				\end{remark}
			}
			
			{
				\begin{lemma}\label{lemma: factorisation of T[epsilon] iff morphism of conormal is zero}
					Consider the following diagram in $\mathrm{Sch}/\Bbbk$ such that $i:Z\to Y$ is an immersion and $j:T\hookrightarrow T[\epsilon]$ is the closed immersion via $\epsilon=0$
					\begin{equation}
						\begin{tikzcd}
							T\ar[r,hook,"j"]\ar[d,"g"]&T[\epsilon]\ar[d,"\rho"]\ar[ld,dashed]\\
							Z\ar[r,"i"]&Y.
						\end{tikzcd}
					\end{equation}
					Then $\rho:T[\epsilon]\to Y$ factors through $i:Z\to Y$ if and only if the morphism $g^*\mathcal C_{Z/Y}\to \epsilon\mathcal O_T$ in \href{https://stacks.math.columbia.edu/tag/01R3}{Lemma 01R3} is zero, where $\epsilon\mathcal O_T\cong\mathcal C_{T/T[\epsilon]}$. 
				\end{lemma}
				\begin{proof}
					Assume without loss of generality that $T,Z,Y$ are affine schemes, and $i:Z\hookrightarrow Y$ is a closed immersion. Let $T=\mathrm{Spec}(A)$, $Y=\mathrm{Spec}(R)$ and $Z=\mathrm{Spec}(R/I)$. Let $\varphi:R/I\to A$ correspond to $g:T\to Z$, and let $\psi:R\to A[\epsilon]$ corresponds to $\rho:T[\epsilon]\to Y$. 
					
					Since $\rho\circ j=i\circ g$, the right square in the following diagram commutes 
					\begin{equation}
						\begin{tikzcd}
							0\ar[r]&I\ar[r]\ar[d,"\psi|_I"]&R\ar[r]\ar[d,"\psi"]&R/I\ar[r]\ar[d,"\varphi"]&0\\
							0\ar[r]&\epsilon A\ar[r]&A[\epsilon]\ar[r]&A\ar[r]&0
						\end{tikzcd}
					\end{equation}
					and thus $\psi|_I:I\to \epsilon A$ exists so that the left square commutes. 
		
					The morphism $g^*\mathcal C_{Z/Y}\to \epsilon\mathcal O_T$ corresponds to the following 
					\begin{equation}
						\eta:I/I^2\otimes_{R/I,\varphi}A\to \epsilon A,\quad \overline h\otimes 1\mapsto \psi|_I(h). 
					\end{equation}
					Then we have $\eta=0$ if and only if $\psi(I)=0$. 
					
					The morphism $\rho:T[\epsilon]\to Y$ factors as $T[\epsilon]\to Z\to Y$ if and only if the ring map $\psi:R\to A[\epsilon]$ factors as $R\to R/I\to A[\epsilon]$, which is equivalent to $\psi(I)=0$, and thus equivalent to $\eta=0$. This completes the proof. 
				\end{proof}
			}
			{
				\begin{lemma}\label{lemma: bijections between infinitesimal transformations stabilising a point with Stab_g(t)}
					Let $\sigma:G\times X\to X$ be an action. Let $t:T\to X$ be a $T$-valued point. Then we have the following commutative diagram of $\Gamma(T,\mathcal O_T)$-modules
					\begin{equation}
						\begin{tikzcd}
							\{\rho:\rho\textrm{ is a }T\textrm{-valued infinitesimal transformation}\}\ar[r,"\cong"{sloped}]&\Gamma(T,\mathfrak g\otimes_\Bbbk\mathcal O_T)\\
							\Big\{\begin{matrix}\rho:\;\rho\textrm{ is a }T\textrm{-valued infinitesimal transformation}\\\textrm{and }\rho\textrm{ stabilises }t\end{matrix}\Big\}\ar[r,"\cong"{sloped}]\ar[u,"\subseteq"{sloped}]&\Gamma\big(T,\mathrm{Stab}_{\mathfrak g}(t)\big)\ar[u,"\subseteq"{sloped}]
						\end{tikzcd}
					\end{equation}
					where both horizontal maps are bijections and the top map sends $\rho:T[\epsilon]\to G$ to the morphism of conormal sheaves $C_\rho:\mathfrak g^*\otimes_\Bbbk\mathcal O_T\to \epsilon\mathcal O_T$ by applying \href{https://stacks.math.columbia.edu/tag/01R3}{Lemma 01R3} to the diagram 
					\begin{equation}
						\begin{tikzcd}
							T\ar[r,hook]\ar[d,equal]&T[\epsilon]\ar[d,"\rho"]\\
							T\ar[r,hook,"{(e,1_T)}"]&G\times T. 
						\end{tikzcd}
					\end{equation}
					Moreover, the maps are natural for $T\in\mathrm{Sch}/X$. 
				\end{lemma}
				\begin{proof}
					It is easy to see the maps are $\Gamma(T,\mathcal O_T)$-linear and are natural for $T\in\mathrm{Sch}/X$. 
		
					We construct the inverse of the top map. Let $\eta\in\Gamma(T,\mathfrak g\otimes_\Bbbk\mathcal O_T)$. We can think $\eta$ as a map of $\Gamma(T,\mathcal O_T)$-modules 
					\begin{equation}
						\eta:\mathfrak g^*\otimes_\Bbbk \Gamma(T,\mathcal O_T)\to \Gamma(T,\mathcal O_T).
					\end{equation}
					Consider the ring map $\rho^*:\mathcal O(G)\to\Gamma(T,\mathcal O_T)[\epsilon]$ such that for $h\in\mathcal O(G)$
					\begin{equation}
						\rho^*(h)=e^*(h)\oplus \epsilon\;\eta\big(\overline{h-e^*(h)}\otimes 1\big)\in\Gamma(T,\mathcal O_T)\oplus\epsilon\Gamma(T,\mathcal O_T)\cong\Gamma(T[\epsilon],\mathcal O_{T[\epsilon]})
					\end{equation}
					where $h-e^*(h)\in\mathfrak m_e:=\ker(e^*)$ represents $\overline{h-e^*(h)}\in\mathfrak m_e/\mathfrak m_e^2\cong \mathfrak g^*$. The ring map $\rho^*$ corresponds to a morphism $\rho:T[\epsilon]\to G$. It is easy to check that the map $\eta\mapsto \rho$ provides the required inverse. This establishes the top bijection. 
					
					Let $\rho:T[\epsilon]\to G$ be a $T$-valued infinitesimal transformation. The following diagram commutes
					\begin{equation}\label{equation: 3 parallelograms}
						\begin{tikzcd}
							&T\ar[r,hook]\ar[d,"t"]&T[\epsilon]\ar[d,"{(\rho,t\circ\mathrm{pr}_T)}"]\\
							&X\ar[r,hook,"\gamma"]\ar[dl]\ar[dr,equal]&G\times X\ar[dl]\ar[dr,"{(\sigma,\mathrm{pr}_X)}"]\\
							\mathrm{Spec}(\Bbbk)\ar[r,hook,"e"]&G&X\ar[r,"\delta"]&X\times X
						\end{tikzcd}
					\end{equation}
					where $\gamma:X\to G\times X$ is the pullback of $e:\mathrm{Spec}(\Bbbk)\hookrightarrow G$ in the bottom left parallelogram, $\delta:X\to X\times X$ is the diagonal morphism, $\mathrm{pr}_T:T[\epsilon]\to T$ is the projection and $\mathrm{pr}_X:G\times X\to X$ is the projection. 

					Apply \href{https://stacks.math.columbia.edu/tag/01R3}{Lemma 01R3} to the top square and the composition of the top square and bottom left parallelogram to obtain two morphisms between conormal sheaves 
					\begin{equation}
						\begin{split}
							\mathcal C_{(\rho,t\circ\mathrm{pr}_T)}:&\mathfrak g^*\otimes_\Bbbk\mathcal O_T\to\epsilon\mathcal O_T\\
							\mathcal C_\rho:&\mathfrak g^*\otimes_\Bbbk\mathcal O_T\to\epsilon\mathcal O_T. 
						\end{split}
					\end{equation}
					such that $\mathcal C_\rho=\mathcal C_{(\rho,t\circ\mathrm{pr}_T)}\circ t^*\mathcal C_{\mathrm{pr}_G}$, where $\mathcal C_{\mathrm{pr}_G}:\mathfrak g^*\otimes_\Bbbk\mathcal O_X\to \mathfrak g^*\otimes_\Bbbk\mathcal O_X$ is the morphism of conormal sheaves associated to the bottom right parallelogram. Since the bottom left parallelogram is Cartesian with southwest arrows flat, by \href{https://stacks.math.columbia.edu/tag/01R4}{Lemma 01R4}, we have $\mathcal C_{\mathrm{pr}_G}=1_{\mathfrak g^*\otimes_\Bbbk\mathcal O_X}$. Then $\mathcal C_\rho=\mathcal C_{(\rho,t\circ\mathrm{pr}_T)}$. 
		
					We have that $\rho$ stabilises $t:T\to X$ if and only if the dashed arrow below exists in the following diagram 
					\begin{equation}\label{equation: diagram of top bottom composition}
						\begin{tikzcd}
							T\ar[r,hook]\ar[d,"f"]&T[\epsilon]\ar[dl,dashed]\ar[d,"{\big(\sigma\circ(\rho,t\circ\mathrm{pr}_T),\;t\circ\mathrm{pr}_T\big)}"]\\
							X\ar[r,"\delta"]&X\times X
						\end{tikzcd}
					\end{equation}
					where the square \eqref{equation: diagram of top bottom composition} is the composition of the top and bottom right parallelogram in diagram \eqref{equation: 3 parallelograms}. Note that $\phi:\Omega_{X/\Bbbk}\to\mathfrak g^*\otimes_\Bbbk \mathcal O_X$ is the morphism of conormal sheaves associated to the bottom right parallelogram. By Lemma \ref{lemma: factorisation of T[epsilon] iff morphism of conormal is zero}, we have that the dashed arrow above exists if and only if $t^*\Omega_{X/\Bbbk}\to \epsilon\mathcal O_T$ is zero, i.e. the following composition is zero
					\begin{equation}
						\begin{tikzcd}
							t^*\Omega_{X/\Bbbk}\ar[r,"t^*\phi"]&\mathfrak g^*\otimes_\Bbbk \mathcal O_T\ar[rr,"\mathcal C_{(\rho,t\circ\mathrm{pr}_T)}"]&&\epsilon\mathcal O_T. 
						\end{tikzcd}
					\end{equation}

					Therefore $\rho$ stabilises $t$ if and only if $\mathcal C_\rho\circ t^*\phi=\mathcal C_{(\rho,t\circ\mathrm{pr}_T)}\circ t^*\phi=0$. Recall that $\rho\mapsto \mathcal C_\rho$ defines a bijection 
					\begin{equation}
						\{\rho:\rho\textrm{ is a }T\textrm{-valued infinitesimal transformation}\}\cong\mathrm{Hom}_{\mathcal O_T}(\mathfrak g^*\otimes_\Bbbk \mathcal O_T,\epsilon\mathcal O_T).
					\end{equation}
					Then the images of those $\rho$ which stabilise $t$ constitute the kernel of 
					\begin{equation}
						\mathrm{Hom}_{\mathcal O_T}(\mathfrak g^*\otimes_\Bbbk \mathcal O_T,\mathcal O_T)\to \mathrm{Hom}_{\mathcal O_T}(t^*\Omega_{X/\Bbbk},\mathcal O_T)
					\end{equation}
					which is $\mathrm{Hom}_{\mathcal O_T}(t^*\mathrm{coker}(\phi),\mathcal O_T)\cong\Gamma\big(T,\mathrm{Stab}_{\mathfrak g}(t)\big)$. This establishes the bottom bijection in the lemma. 
				\end{proof}
			}
		}
	}

	\subsection{Some conditions for the linearisation}
	{
		From now on, all schemes are of finite type over $\Bbbk$. Let $X$ be a projective scheme and let $L$ be an ample line bundle on $X$. Let $\hat U$ act on $X$ linearly with respect to $L$. In \cite{HoskinsVictoria2021QbPG}, conditions called \emph{Upstairs Unipotent Stabiliser} (\cite{HoskinsVictoria2021QbPG} Assumption 4.40) and \emph{Weak Upstairs Unipotent Stabiliser} (\cite{HoskinsVictoria2021QbPG} Assumption 4.42) were considered for the linearisation $\hat U\curvearrowright(X,L)$, when $X$ is reduced. We will generalise these conditions for potentially non-reduced $X$. Firstly, a non-emptiness condition is required, which always holds when $X$ is reduced but can fail in non-reduced cases. 

		\subsubsection{A non-emptiness condition}
		{
			We will define an open subscheme $X^0_{\min}\subseteq X$ and a closed subscheme $Z_{\min}\hookrightarrow X$. We need a condition to ensure that $X^0_{\min}\ne\emptyset$ and $Z_{\min}\ne\emptyset$. 
			
			We first consider when $\hat U$ acts on projective spaces. Let $V$ be a finite dimensional representation of $\hat U$. For any $w\in\mathbb Z$ denote 
			\begin{equation}
				V_{\lambda=w}:=\{v\in V:\lambda(t).v=t^wv,\;\textrm{for all closed point }t\in\mathbb G_m\}. 
			\end{equation}
			We call $w\in\mathbb Z$ a \emph{weight} if $V_{\lambda=w}\ne0$, and we call $V_{\lambda=w}\subseteq V$ the weight subspace of the weight $w$ with respect to $\lambda:\mathbb G_m\to \hat U$. Let $V_{\lambda=\min}\subseteq V$ denote the weight subspace whose weight is the minimum. If $V\ne0$, then $V_{\min}\ne0$. Let $w_{\min}$ denote the minimal weight. For $v\in V$, decompose according to $\mathbb G_m$-weights 
			\begin{equation}
				v=v_{\lambda=\min}+v_{\lambda>\min},\quad v_{\lambda=\min}\in V_{\lambda=\min},\;v_{\lambda>\min}\in \bigoplus_{w>w_{\min}}V_{\lambda=w}. 
			\end{equation}

			Define locally closed subschemes $\mathbb P(V)^0_{\min}\subseteq \mathbb P(V)$ and $Z_{\min}\big(\mathbb P(V)\big)\hookrightarrow \mathbb P(V)$ whose sets of closed points are 
			\begin{equation}
				\begin{split}
					\mathbb P(V)^0_{\min}\overset{\mathrm{cl}}{=\!=}&\big\{[v]\in\mathbb P(V):v_{\lambda=\min}\ne0\big\}\\
					Z_{\min}\big(\mathbb P(V)\big)\overset{\mathrm{cl}}{=\!=}&\big\{[v]\in\mathbb P(V):v_{\lambda>\min}=0\big\}
				\end{split}
			\end{equation}
			where $\overset{\mathrm{cl}}{=\!=}$ means they have the same closed points. It is easy to see that $\mathbb P(V)^0_{\min}$ is an open subscheme and $Z_{\min}\big(\mathbb P(V)\big)$ is a closed subscheme, and $Z_{\min}\big(\mathbb P(V)\big)\hookrightarrow \mathbb P(V)^0_{\min}$. 

			The non-emptiness condition for projective spaces is vacuous, since if $\mathbb P(V)\ne\emptyset$, then $\mathbb P(V)^0_{\min}\ne\emptyset$ and $Z_{\min}\big(\mathbb P(V)\big)\ne\emptyset$. It is vacuous if $X$ is reduced. 

			Then we consider the linearisation $\hat U\curvearrowright(X,L)$. There exists $m\in\mathbb N_+$ such that $H^0(X,L^m)$ generates $\bigoplus_{r\in\mathbb N}H^0(X,L^{mr})$. We can define $\mathbb P\big(H^0(X,L^m)^*\big)^0_{\min}$ and $Z_{\min}\big(\mathbb P\big(H^0(X,L^m)^*\big)\big)$ as above. 

			{
				\begin{lemma}\label{lemma: non-emptiness of X^0_min does not depend on m when divisible}
					For the linearisation $\hat U\curvearrowright(X,L)$, if $m\in\mathbb N_+$ is such that $H^0(X,L^m)$ generates $\bigoplus_{r\in\mathbb N}H^0(X,L^{mr})$ and 
					\begin{equation}
						X\cap \mathbb P\big(H^0(X,L^m)^*\big)^0_{\min}\ne\emptyset
					\end{equation}
					where $X\hookrightarrow \mathbb P\big(H^0(X,L^m)^*\big)$ is a closed immersion, then for all $n\in\mathbb N_+$ we have 
					\begin{equation}
						X\cap \mathbb P\big(H^0(X,L^m)^*\big)^0_{\min}=X\cap \mathbb P\big(H^0(X,L^{mn})^*\big)^0_{\min}. 
					\end{equation}
				\end{lemma}
				\begin{proof}
					Let $m\in\mathbb N_+$ be as in the lemma. The linear system $H^0(X,L^{mn})$ induces a closed immersion $X\hookrightarrow \mathbb P\big(H^0(X,L^{mn})^*\big)$ for all $n\in\mathbb N_+$. 

					The open subscheme $X\cap \mathbb P\big(H^0(X,L^m)^*\big)^0_{\min}$ is the non-vanishing locus of $H^0(X,L^m)_{\lambda=\max}$, where $H^0(X,L^m)_{\lambda=\max}$ is the weight subspace of $H^0(X,L^m)$ with respect to $\lambda$ whose weight is the maximum. Similarly, the open subscheme $X\cap \mathbb P\big(H^0(X,L^{mn})^*\big)^0_{\min}$ is the non-vanishing locus of $H^0(X,L^{mn})_{\lambda=\max}$. 

					Let $w_{\max}$ denote the maximal weight of $\lambda$ on $H^0(X,L^m)$. The condition that $X\cap \mathbb P\big(H^0(X,L^m)^*\big)^0_{\min}\ne\emptyset$ implies that there exists $f\in H^0(X,L^m)_{\lambda=\max}$, which is not nilpotent in $\bigoplus_{r\in\mathbb N_+}H^0(X,L^r)$. Therefore $f^n\in H^0(X,L^{mn})\setminus\{0\}$. Since $H^0(X,L^m)$ generates $\bigoplus_{r\in\mathbb N_+}H^0(X,L^{mr})$, we have that the maximal weight on $H^0(X,L^{mn})$ is the weight of $f^n$, i.e. $nw_{\max}$. We have the following surjective map 
					\begin{equation}
						\mathrm{Sym}^nH^0(X,L^m)_{\lambda=\max}\to H^0(X,L^{mn})_{\lambda=\max}. 
					\end{equation}
					Therefore the non-vanishing locus of $H^0(X,L^{mn})_{\lambda=\max}$ and the non-vanishing locus of $H^0(X,L^m)_{\lambda=\max}$ coincide. This proves the lemma. 
				\end{proof}

				\begin{condition}\label{condition: non-emptiness of X^0_min}
					For the linearisation $\hat U\curvearrowright(X,L)$, assume that 
					\begin{equation}
						X\cap \mathbb P\big(H^0(X,L^m)^*\big)^0_{\min}\ne\emptyset
					\end{equation}
					for sufficiently divisible $m\in\mathbb N_+$. Denote 
					\begin{equation}
						X^0_{\min}:=X\cap \mathbb P\big(H^0(X,L^m)^*\big). 
					\end{equation}
				\end{condition}

				\begin{remark}
					Lemma \ref{lemma: non-emptiness of X^0_min does not depend on m when divisible} implies that $X^0_{\min}$ does not depend on $m\in\mathbb N_+$ when $m$ is sufficiently divisible and when Condition \ref{condition: non-emptiness of X^0_min} holds. 
				\end{remark}
			}

			The opposite of Condition \ref{condition: non-emptiness of X^0_min} is the condition that $H^0(X,L^m)_{\lambda=\max}$ is nilpotent for all sufficiently divisible $m\in\mathbb N_+$. When $X$ is non-empty and reduced, then Condition \ref{condition: non-emptiness of X^0_min} always holds. 

			With Condition \ref{condition: non-emptiness of X^0_min}, we also have that $X^0_{\min}$ is the non-vanishing locus of $\bigoplus_{r\in\mathbb N_+}H^0(X,L^r)_{\lambda=\max}$. Define the \emph{fixed point subscheme} $Z_{\min}\hookrightarrow X^0_{\min}$ 
			\begin{equation}
				Z_{\min}:=\big(X^0_{\min}\big)^\lambda. 
			\end{equation}
			Then $Z_{\min}=\mathrm{Proj}\Big(\bigoplus_{r\in\mathbb N_+}H^0(X,L^r)_{\lambda=\max}\Big)$. 
		}

		\subsubsection{Upstairs Unipotent Stabiliser conditions}
		{
			Let $\phi:\Omega_{X/\Bbbk}\to\mathfrak u^*\otimes_\Bbbk\mathcal O_X$ be the infinitesimal action of $\mathfrak u$ on $X$ (Definition \ref{definition: infinitesimal action of g on X}), and let $\mathrm{Fit}_k(\phi)\subseteq \mathcal O_X$ denote the $k$th Fitting ideal of $\mathrm{coker}(\phi)$. The following conditions generalise \cite{hoskins2021quotients} Assumption 4.40 and Assumption 4.42. 

			{
				\begin{condition}[Upstairs Unipotent Stabiliser]\label{condition: UU}
					For the linearisation $\hat U\curvearrowright(X,L)$ satisfying Condition \ref{condition: non-emptiness of X^0_min}, we say the \emph{upstairs unipotent stabiliser condition} (Condition UU) holds if $\mathrm{coker}(\phi)$ is locally free of constant rank on $X^0_{\min}$. 
				\end{condition}

				\begin{remark}
					By Lemma \ref{lemma: equivalence of dimStab>k and Fit_k in m_x}, Condition \hyperref[condition: UU]{UU} implies that $\dim\mathrm{Stab}_U(x)$ is constant for $\Bbbk$-points $x\in X^0_{\min}$. When $X^0_{\min}$ is reduced, the converse is true. 
				\end{remark}
			}

			{

				\begin{condition}[Weak Upstairs Unipotent Stabiliser]\label{condition: WUU}
					For the linearisation $\hat U\curvearrowright(X,L)$ satisfying Condition \ref{condition: non-emptiness of X^0_min}, 
					we say the \emph{weak upstairs unipotent stabiliser condition} (Condition WUU) holds if $\mathrm{coker}(\phi)$ is locally free of rank $k$ on an open subscheme which meets $Z_{\min}$, where $k:=\min\{i\in\mathbb N:\mathrm{Fit}_i(\phi)\ne0\}$. 
				\end{condition}
			}
		}
	}

	\subsection{Theorem with Condition \texorpdfstring{\hyperref[condition: UU]{UU}}{UU}}
	{
		When the linearisation $\hat U\curvearrowright(X,L)$ satisfies Condition \hyperref[condition: UU]{UU}, the open subscheme $X^0_{\min}$ admits a geometric quotient by $U$, and the open subscheme $X^0_{\min}\setminus UZ_{\min}$ admits a projective geometric quotient by $\hat U$. We will prove these statements. 

		Let $r:=\dim U$. Without loss of generality assume $r>0$. Let $\phi:\Omega_{X/\Bbbk}\to \mathfrak u^*\otimes_\Bbbk\mathcal O_X$ be the infinitesimal action of $\mathfrak u$. Condition \hyperref[condition: UU]{UU} implies that $X^0_{\min}\ne\emptyset$. 

		\subsubsection{Reduction to affine case}\label{subsubsection: reduction to affine case}
		{
			Condition \hyperref[condition: UU]{UU} assumes that $\mathrm{coker}(\phi)$ is locally free of rank $k$ on $X^0_{\min}$ for some integer $0\leq k\leq r$. The quotient $\mathfrak u^*\otimes_\Bbbk\mathcal O_{X^0_{\min}}\to \mathrm{coker}(\phi)|_{X^0_{\min}}$ corresponds to a morphism 
			\begin{equation}
				\Theta:X^0_{\min}\to \mathrm{Grass}(\mathfrak u^*,k)
			\end{equation}
			where $\mathrm{Grass}(\mathfrak u^*,k)$ is the Grassmannian of $k$-dimensional quotients of $\mathfrak u^*$. Let $\iota:Z_{\min}\hookrightarrow X^0_{\min}$ be the closed immersion. 

			There exists a finite family of $k$-dimensional vector subspaces $\big\{\mathfrak s_i^*\subseteq\mathfrak u^*\}_{i=1}^s$ such that 
			\begin{equation}
				\mathrm{Grass}(\mathfrak u^*,k)=\bigcup_{i=1}^s\mathrm{Grass}(\mathfrak u^*,k)_{\mathfrak s_i^*}
			\end{equation}
			where $\mathrm{Grass}(\mathfrak u^*,k)_{\mathfrak s_i^*}\subseteq\mathrm{Grass}(\mathfrak u^*,k)$ is the open subscheme whose closed points are quotients $\mathfrak u^*\to Q$ such that the composition $\mathfrak s_i^*\subseteq \mathfrak u^*\to Q$ is an isomorphism. 

			Open subsets $Z_{\min,f}\subseteq Z_{\min}$ for $m\in\mathbb N_+$ and $f\in H^0(X,L^m)_{\lambda=\max}$ form a basis of the Zariski topology of $Z_{\min}$. We can choose a sufficiently divisible $m\in\mathbb N_+$ and $f_1,\cdots,f_t\in H^0(X,L^m)_{\lambda=\max}$ such that $\{Z_{\min,f_i}\}_{i=1}^t$ is a covering of $Z_{\min}$, and the covering is subordinate to the covering of $\mathrm{Grass}(\mathfrak u^*,k)$, in the sense that for any $1\leq i\leq t$, there exists $1\leq j\leq s$ such that 
			\begin{equation}
				Z_{\min,f_i}\subseteq \Theta^{-1}\Big(\mathrm{Grass}(\mathfrak u^*,k)_{\mathfrak s_j^*}\Big). 
			\end{equation}
			The morphism $\Theta$ is $\mathbb G_m$-invariant and $\mathbb G_m$-equivariant, so 
			\begin{equation}
				\begin{split}
					&Z_{\min,f_i}\subseteq \Theta^{-1}\Big(\mathrm{Grass}(\mathfrak u^*,k)_{\mathfrak s_j^*}\Big)\\
					\textrm{if and only if }\;&X_{f_i}\subseteq \Theta^{-1}\Big(\mathrm{Grass}(\mathfrak u^*,k)_{\mathfrak s_j^*}\Big). 
				\end{split}
			\end{equation}

			We have achieved the following: 
			\begin{itemize}
				\item There exists $m\in\mathbb N_+$ and $f_1,\cdots,f_t\in H^0(X,L^m)_{\lambda=\max}$ such that $X^0_{\min}$ is covered by $\{X_{f_i}\}_{i=1}^t$; 
				\item For each $X_{f_i}$, there exists a linear subspace $\mathfrak s_i^*\subseteq \mathfrak u^*$ such that $\mathfrak s_i^*\otimes_\Bbbk\mathcal O_{X_{f_i}}\cong \mathrm{coker}(\phi)|_{X_{f_i}}$ as follows 
				\begin{equation}
					\begin{tikzcd}
						&\mathfrak s_i^*\otimes_\Bbbk\mathcal O_{X_{f_i}}\ar[d,sloped,"\subseteq"]\ar[rd,"\cong"]\\
						\Omega_{X_{f_i}/\Bbbk}\ar[r,"\phi|_{X_{f_i}}"]&\mathfrak u^*\otimes_\Bbbk\mathcal O_{X_{f_i}}\ar[r]&\mathrm{coker}(\phi)|_{X_{f_i}}\ar[r]&0. 
					\end{tikzcd}
				\end{equation}
			\end{itemize}
		}

		\subsubsection{Affine theory with Condition \texorpdfstring{\hyperref[condition: UU]{UU}}{UU}}
		{
			We work on each of $X_{f_i}$ and construct a universal geometric quotient by $U$. 

			Let $A$ be a finitely generated commutative $\Bbbk$-algebra. Let $\hat U$ act on $\mathrm{Spec}(A)$. Assume the following conditions for $\hat U\curvearrowright \mathrm{Spec}(A)$: 
			\begin{itemize}
				\item[(1)] If $A_{\lambda=w}\subseteq A$ denotes the weight subspace whose $\lambda$-weight is $w\in \mathbb Z$ for the $\mathbb G_m$-representation on $A$, then $A_{\lambda=w}=0$ for $w>0$; 
				\item[(2)] There exists a linear subspace $\mathfrak s^*\subseteq \mathfrak u^*$ such that $\mathfrak s^*\otimes_\Bbbk A\cong \mathrm{coker}(\phi_A)$ in the following diagram 
				\begin{equation}\label{diagram: affine UU}
					\begin{tikzcd}
						&\mathfrak s^*\otimes_\Bbbk A\ar[d,sloped,"\subseteq"]\ar[rd,"\cong"]\\
						\Omega_{A/\Bbbk}\ar[r,"\phi_A"]&\mathfrak u^*\otimes_\Bbbk A\ar[r]&\mathrm{coker}(\phi_A)\ar[r]&0. 
					\end{tikzcd}
				\end{equation}
			\end{itemize}
			The action $\hat U\curvearrowright X_{f_i}$ satisfies these conditions. 
			{
				\begin{remark}
					If $\rho:A\to \Bbbk[t,t^{-1}]\otimes_\Bbbk A$ denotes the co-action map of $\mathbb G_m\curvearrowright\mathrm{Spec}(A)$, then $A_{\lambda=w}=\{f\in A:\rho(f)=t^{-w}\otimes f\}$. Condition $(1)$ above is equivalent to that $\rho$ factors through $\Bbbk[t]\otimes_\Bbbk A$ 
					\begin{equation}
						\begin{tikzcd}
							A\ar[r,"\rho"]\ar[rd,dashed]&\Bbbk[t,t^{-1}]\otimes_\Bbbk A\\
							&\Bbbk[t]\otimes_\Bbbk A\ar[u,sloped,"\subseteq"]. 
						\end{tikzcd}
					\end{equation}
					Geometrically, Condition $(1)$ is equivalent to that the action morphism $\mathbb G_m\times\mathrm{Spec}(A)\to \mathrm{Spec}(A)$ extends to a morphism $\mathbb A^1\times\mathrm{Spec}(A)\to\mathrm{Spec}(A)$. 
				\end{remark}
			}

			Let $\mathfrak u\to\mathfrak s$ denote the surjective linear map dual to $\mathfrak s^*\subseteq \mathfrak u^*$. Let $\mathfrak u':=\ker(\mathfrak u\to\mathfrak s)$. Let $\sigma^*:A\to \mathcal O(U)\otimes_\Bbbk A$ denote the co-action map associated to $U\curvearrowright\mathrm{Spec}(A)$. Define sub-rings invariant under $\mathfrak u'$ and $\mathfrak u$
			\begin{equation}
				\begin{split}
					A^{\mathfrak u'}:=&\{f\in A:\xi.f=0\textrm{ for all }\xi\in\mathfrak u'\}\\
					A^{\mathfrak u}:=&\{f\in A:\xi.f=0\textrm{ for all }\xi\in\mathfrak u\}
				\end{split}
			\end{equation}
			where $\xi:A\to A$ is the $\Bbbk$-derivation defined as follows 
			\begin{equation}
				\begin{tikzcd}
					A\ar[rrrrr,bend left=10,"\xi",yshift=10]\ar[r,"\sigma^*"]&\mathcal O(U)\otimes_\Bbbk A\ar[rr,"\exp^*\otimes 1_A"]&&\mathfrak u^*\otimes_\Bbbk A\ar[rr,"-\xi\otimes 1_A"]&&A
				\end{tikzcd}
			\end{equation}
			for $\exp^*:\mathcal O(U)\to \mathfrak u^*$ sending $f$ to $\overline{f-f(e)}\in\mathfrak m_e/\mathfrak m_e^2\cong \mathfrak u^*$. It is easy to check $A^{\mathfrak u}\subseteq A^{\mathfrak u'}$ are sub-rings. 

			{
				\begin{lemma}\label{lemma: equality of invariant subrings with UU}
					We have $A^{\mathfrak u}=A^{\mathfrak u'}$. 
				\end{lemma}
				\begin{proof}
					Let $\{\xi_1,\cdots,\xi_r\}$ be a basis of $\mathfrak u$ such that $\xi_1,\cdots,\xi_{r-k}\in\mathfrak u'$. Let $\{u_1,\cdots,u_r\}$ be the basis of $\mathfrak u^*$ dual to $\{\xi_1,\cdots,\xi_r\}$. Then $\{\xi_{r-k+1},\cdots,\xi_r\}$ is a basis of $\mathfrak s^*$. The infinitesimal action $\phi_A:\Omega_{A/\Bbbk}\to \mathfrak u^*\otimes_\Bbbk A$ is determined by the equation 
					\begin{equation}
						\phi_A(\mathrm{d}f)=\sum_{i=1}^ru_i\otimes \xi_i.f,\quad f\in A. 
					\end{equation}

					For $\xi\in\mathfrak u$ and $f\in A$ we have 
					\begin{equation}
						\langle\phi_A(\mathrm{d}f),\;\xi\rangle=\sum_{i=1}^r\langle u_i,\xi\rangle\otimes\xi_i.f=\xi.f
					\end{equation}
					and therefore we have 
					\begin{equation}
						\begin{split}
							A^{\mathfrak u'}=&\{f\in A:\phi_A(\mathrm{d}f)\subseteq \mathfrak s^*\otimes_\Bbbk A\}\\
							A^{\mathfrak u}=&\{f\in A:\phi_A(\mathrm{d}f)=0\}. 
						\end{split}
					\end{equation}

					By diagram \eqref{diagram: affine UU}, we have $\mathrm{im}(\phi_A)\cap (\mathfrak s^*\otimes_\Bbbk A)=0$. The lemma is proved. 
				\end{proof}
			}

			We first construct the quotient of $\mathrm{Spec}(A)$ by the subgroup $U':=\mathrm{exp}(\mathfrak u')\subseteq U$. 
			{
				\begin{lemma}\label{lemma: existence of x_i, xi_j}
					There exists $\xi_1,\cdots,\xi_{r-k}\in\mathfrak u'$ and $x_1,\cdots,x_{r-k}\in A$ such that 
					\begin{itemize}
						\item $\xi_i.x_j=\delta_{ij}$ for $1\leq i,j\leq r-k$; 
						\item $\xi.x_j\in A_{\lambda=0}$ for $\xi\in\mathfrak u$ and $1\leq j\leq r-k$. 
					\end{itemize}
				\end{lemma}
				\begin{proof}
					Let $w\in\mathbb N_+$ be the unique $\lambda$-weight on $\mathfrak u$. Let $\{\xi_1,\cdots,\xi_r\}$ be a basis of $\mathfrak u$ such that $\xi_1,\cdots,\xi_{r-k}\in\mathfrak u'$. Let $\{u_1,\cdots,u_r\}$ be the basis of $\mathfrak u^*$ dual to $\{\xi_1,\cdots,\xi_r\}$. 

					By diagram \eqref{diagram: affine UU}, the following composition is surjective 
					\begin{equation}
						\begin{tikzcd}
							\Omega_{A/\Bbbk}\ar[r,"\phi_A"]&\mathfrak u^*\otimes_\Bbbk A\ar[r,->>]&(\mathfrak u')^*\otimes_\Bbbk A. 
						\end{tikzcd}
					\end{equation}
					Then for $1\leq i\leq r-k$, there exists $\omega_i\in\Omega_{A/\Bbbk}$ such that 
					\begin{equation}
						\phi_A(\omega_i)\in u_i\otimes 1+\mathfrak s^*\otimes_\Bbbk A. 
					\end{equation}

					Write $\omega_i=\sum_{n,\mu,\nu} a^{(i)}_{n,\mu}\;\mathrm{d}b^{(i)}_{n,\nu}$ for $a^{(i)}_{n,\mu}\in A_{\lambda=\mu},b^{(i)}_{n,\nu}\in A_{\lambda=\nu}$. Then we have 
					\begin{equation}
						\begin{split}
							\phi_A(\omega_i)=&\sum_{j=1}^r\sum_{n,\mu,\nu}u_j\otimes a^{(i)}_{n,\mu}\;\xi_j.b^{(i)}_{n,\nu}\\
							=&\sum_{j=1}^r\sum_nu_j\otimes a^{(i)}_{n,0}\;\xi_j.b^{(i)}_{n,-w}+\sum_{j=1}^r\sum_{\substack{n,\mu,\nu\\\mu<0\\\textrm{or }\nu<-w}}u_j\otimes a^{(i)}_{n,\mu}\;\xi_j.b^{(i)}_{n,\nu}. 
						\end{split}
					\end{equation}

					Let $x_i:=\sum_na^{(i)}_{n,0}b^{(i)}_{n,-w}\in A_{\lambda=-w}$ for $1\leq i\leq r-k$. Since $x_i\in A_{\lambda=-w}$ and any $\xi\in\mathfrak u\setminus\{0\}$ has weight $w\in\mathbb N_+$, we have $\xi.x_i\in A_{\lambda=0}$. This checks the second condition in the lemma. 

					Let $p:A\to A_{\lambda=0}$ be the projection onto the weight subspace of weight 0, and let $q:\mathfrak u^*\to(\mathfrak u')^*$ be the surjective linear map. Then in $(\mathfrak u')^*\otimes_\Bbbk A_{\lambda=0}$
					\begin{equation}
						\begin{split}
							u_i\otimes 1=&(q\otimes p)\circ\phi_A(\omega_i)\\
							=&(q\otimes p)\Big(\sum_{j=1}^r\sum_nu_j\otimes a^{(i)}_{n,0}\;\xi_j.b^{(i)}_{n,-w}+\sum_{j=1}^r\sum_{\substack{n,\mu,\nu\\\mu<0\\\textrm{or }\nu<-w}}u_j\otimes a^{(i)}_{n,\mu}\;\xi_j.b^{(i)}_{n,\nu}\Big)\\
							=&\sum_{j=1}^{r-k}u_j\otimes\xi_j.\Big(\sum_na^{(i)}_{n,0}b^{(i)}_{n,-w}\Big)\\
							=&\sum_{j=1}^{r-k}u_j\otimes\xi_j.x_i
						\end{split}
					\end{equation}
					that is $\xi_i.x_j=\delta_{ij}$ for $1\leq i,j\leq r-k$. This checks the first condition in the lemma. 
				\end{proof}
			}

			{
				\begin{proposition}\label{proposition: affine GIT with UU}
					Let $\mathrm{Spec}(A)$ be an algebraic affine scheme over $\Bbbk$. Let $\hat U$ act on $\mathrm{Spec}(A)$. Assume the action $\hat U\curvearrowright \mathrm{Spec}(A)$ satisfies: 
					\begin{itemize}
						\item[(1)] $A_{\lambda=w}=0$ if $w>0$; 
						\item[(2)] There exists a subspace $\mathfrak s^*\subseteq\mathfrak u^*$ such that $\mathfrak s^*\otimes_\Bbbk A\cong \mathrm{coker}(\phi_A)$ in the diagram 
						\begin{equation}
							\begin{tikzcd}
								&\mathfrak s^*\otimes_\Bbbk A\ar[d,sloped,"\subseteq"]\ar[rd,"\cong"]\\
								\Omega_{A/\Bbbk}\ar[r,"\phi_A"]&\mathfrak u^*\otimes_\Bbbk A\ar[r]&\mathrm{coker}(\phi_A)\ar[r]&0
							\end{tikzcd}
						\end{equation}
						where $\phi_A$ is the infinitesimal action of $\mathfrak u$ on $\mathrm{Spec}(A)$. 
					\end{itemize}
					Then: 
					\begin{itemize}
						\item $A$ is a polynomial ring over $A^{\mathfrak u}$; 
						\item The affine morphism $\mathrm{Spec}(A)\to\mathrm{Spec}(A^{\mathfrak u})$ is a universal geometric quotient by $U$. 
					\end{itemize}
				\end{proposition}
				\begin{proof}
					Let $\mathfrak u':=\ker(\mathfrak u\to \mathfrak s)$, where $\mathfrak u\to \mathfrak s$ is dual to $\mathfrak s^*\subseteq\mathfrak u^*$. We can apply \cite{Greuel1993GeometricQO} Theorem 3.10 for the action of $\mathfrak u'$ on $A$. Lemma \ref{lemma: existence of x_i, xi_j} shows that Condition $(3)$ in \cite{Greuel1993GeometricQO} Theorem 3.10 is satisfied and then we have Condition $(6)$ in that theorem that $\mathrm{Spec}(A)\to\mathrm{Spec}(A^{\mathfrak u'})$ is a trivial geometric quotient by $U'=\exp(\mathfrak u')$. Actually according to the proof of \cite{Greuel1993GeometricQO} Theorem 3.10, we have that $x_1,\cdots,x_{r-k}\in A$ are algebraically independent over $A^{\mathfrak u'}$ and $A=A^{\mathfrak u'}[x_1,\cdots,x_{r-k}]$. Since $A^{\mathfrak u}=A^{\mathfrak u'}$ by Lemma \ref{lemma: equality of invariant subrings with UU}, we have 
					\begin{equation}
						A=A^{\mathfrak u}[x_1,\cdots,x_{r-k}]. 
					\end{equation}
					This proves that $A$ is a polynomial ring over $A^{\mathfrak u}$. The statement that $\mathrm{Spec}(A)\to \mathrm{Spec}(A^{\mathfrak u})$ is a universal geometric quotient follows immediately. 
				\end{proof}

				\begin{corollary}\label{corollary: affine quotient represents quotient sheaf with UU}
					Let $\hat U\curvearrowright \mathrm{Spec}(A)$ and $\mathfrak s^*\subseteq\mathfrak u^*$ be as in Proposition \ref{proposition: affine GIT with UU}. Let $X:=\mathrm{Spec}(A)$ and $Y:=\mathrm{Spec}(A^{\mathfrak u})$. Then $Y$ represents the quotient presheaf $X/_pU$ 
					\begin{equation}
						X/_pU:(\mathrm{Sch}/\Bbbk)^{\mathrm{op}}\to\mathrm{Set},\quad T\mapsto \frac{\mathrm{Hom}_{\mathrm{Sch}/\Bbbk}(T,X)}{\mathrm{Hom}_{\mathrm{Sch}/\Bbbk}(T,U)}. 
					\end{equation}

					Moreover for any morphism $Y'\to Y$ in $\mathrm{Sch}/\Bbbk$, we have that $Y'$ represents the quotient sheaf $(X'/_pU)^\sharp$ on $(\mathrm{Sch}/\Bbbk)_{\acute etale}$ (See \href{https://stacks.math.columbia.edu/tag/02VG}{Definition 02VG} for quotient sheaves), where $X':=X\times_YY'$, and $(X'/_pU)^\sharp$ is the sheafification of $X'/_pU$. 
				\end{corollary}
				\begin{proof}
					For a scheme, we think it as a sheaf on $(\mathrm{Sch}/\Bbbk)_{\acute etale}$ via its functor of points. Let $\pi:X\to Y$ be the geometric quotient in Proposition \ref{proposition: affine GIT with UU}. Since $\pi:X\to Y$ is $U$-invariant, the map $\mathrm{Hom}_{\mathrm{Sch}/\Bbbk}(T,X)\to \mathrm{Hom}_{\mathrm{Sch}/\Bbbk}(T,Y)$ factors through $(X/_pU)(T)$ for any $T\in\mathrm{Sch}/\Bbbk$. This gives a morphism $\eta:X/_pU\to Y$ of presheaves. 

					By Lemma \ref{lemma: existence of x_i, xi_j} and Proposition \ref{proposition: affine GIT with UU}, we can choose a basis $\xi_1,\cdots,\xi_r$ of $\mathfrak u$ and $x_1,\cdots,x_{r-k}\in A$ such that: 
					\begin{itemize}
						\item $A=A^{\mathfrak u}[x_1,\cdots,x_{r-k}]$ and $x_1,\cdots,x_{r-k}$ are algebraically independent over $A^{\mathfrak u}$; 
						\item $\xi_i.x_j=\delta_{ij}$ if $1\leq i,j\leq r-k$. 
					\end{itemize}
					Let $u_1,\cdots,u_r\in \mathfrak u^*$ be the basis dual to $\xi_1,\cdots,\xi_r$. Lift $u_i\in\mathfrak u^*$ to $u_i\in\mathfrak m_e\subseteq\mathcal O(U)$. We identify $\mathcal O(U)\cong \Bbbk[u_1,\cdots,u_r]$. Then for $1\leq i\leq r-k$ 
					\begin{equation}
						\sigma^*:A\to\mathcal O(U)\otimes_\Bbbk A,\quad x_i\mapsto 1\otimes x_i+u_i\otimes 1+\mathrm{err}_i,\quad\mathrm{err}_i\in (u_{r-k+1}\otimes 1,\cdots,u_r\otimes 1). 
					\end{equation}

					We first show that $\eta:X/_pU\to Y$ is injective. Let $T\in\mathrm{Sch}/\Bbbk$. Let $g_i:T\to X$ for $i=1,2$ be two morphisms such that $\eta_T([g_1])=\eta_T([g_2])=h$, i.e. $h=\pi\circ g_1=\pi\circ g_2$. We identify $\mathcal O(U)\cong \Bbbk[u_1,\cdots,u_r]$ and define a ring map 
					\begin{equation}
						l^*:\Bbbk[u_1,\cdots,u_r]\to\Gamma(T,\mathcal O_T),\quad u_i\mapsto \begin{cases}
						g_2^*(x_i)-g_1^*(x_i),&\quad 1\leq i\leq r-k\\
						0,&\quad r-k<i\leq r. 
						\end{cases}
					\end{equation}
					Then the following diagram commutes 
					\begin{equation}
						\begin{tikzcd}
							\Gamma(T,\mathcal O_T)&A\ar[l,"g_2^*"]\\
							\Gamma(T,\mathcal O_T)\otimes_\Bbbk\Gamma(T,\mathcal O_T)\ar[u,"\Delta^*"]\\
							\mathcal O(U)\otimes_\Bbbk A\ar[u,"l^*\otimes g_1^*"]&A\ar[l,"\sigma^*"]\ar[uu,equal]
						\end{tikzcd}
					\end{equation}
					which corresponds to a commutative diagram of schemes 
					\begin{equation}
						\begin{tikzcd}
							T\ar[r,"g_2"]\ar[d,"\Delta"]&X\ar[dd,equal]\\
							T\times T\ar[d,"l\times g_1"]\\
							U\times X\ar[r,"\sigma"]&X. 
						\end{tikzcd}
					\end{equation}
					where $l:T\to U$ corresponds to the ring map $l^*:\mathcal O(U)\to\Gamma(T,\mathcal O_T)$. Then $g_2=\sigma(l,g_1)$ for this $l:T\to U$. This implies $[g_1]=[g_2]\in (X/_pU)(T)$, so $\eta$ is injective. 

					We then show $\eta$ is surjective as a morphism of presheaves. Let $T\in\mathrm{Sch}/\Bbbk$. Let $h:T\to Y$ be in $\mathrm{Hom}_{\mathrm{Sch}/\Bbbk}(T,Y)$. Let $g^*:A\to \Gamma(T,\mathcal O_T)$ be the composition $A= A^{\mathfrak u}[x_1,\cdots,x_{r-k}]\to A^{\mathfrak u}[x_1,\cdots,x_{r-k}]/(x_1,\cdots,x_{r-k})\cong A^{\mathfrak u}\xrightarrow{h^*}\Gamma(T,\mathcal O_T)$. Let $g:T\to X$ correspond to $g^*:A\to\Gamma(T,\mathcal O_T)$. Then $\eta_T([g])=h$. This shows that $\eta_T$ is surjective. Then $\eta:X/_pU\cong Y$. 

					We first prove that if $Y'\to Y$ is a morphism between affine schemes, then $X'/_pU\cong Y'$ as presheaves. Let $Y'\to Y$ be a morphism with $Y'$ affine with coordinate ring $B$. Then $X'=X\times_YY'$ is affine with coordinate ring $A\otimes_{A^{\mathfrak u}}B=B[x_1,\cdots,x_{r-k}]$. The same argument as above with $A^{\mathfrak u}$ replaced by $B$ shows that $X'/_pU\to Y'$ is injective and surjective as presheaves. 

					Next, we consider a general morphism $Y'\to Y$ in $\mathrm{Sch}/\Bbbk$. Choose a covering $\{Y'_j\to Y'\}_{j\in J}$ in $(\mathrm{Sch}/\Bbbk)_{\acute etale}$ such that each of $Y'_j$ is affine. Let $B_j$ be the coordinate ring of $Y'_j$. The composition $Y'_j\to Y'\to Y$ corresponds to a ring map $A^{\mathfrak u}\to B_j$. Let $X'_j:=X\times_YY'_j$. Then $X'_j$ is affine with coordinate ring $A\otimes_{A^{\mathfrak u}}B_j\cong B_j[x_1,\cdots,x_{r-k}]$. 

					We construct a morphism of presheaves $Y'\to (X'/_pU)^\sharp$. By the Yoneda lemma, this is equivalent to an element in $(X'/_pU)^\sharp(Y')$. The sheaf property of $(X'/_pU)^\sharp$ implies that the following is an equaliser diagram 
					\begin{equation}
						\begin{tikzcd}
							(X'/_pU)^\sharp(Y')\ar[r]&\displaystyle\prod_{j\in J}(X'/_pU)^\sharp(Y'_j)\ar[r,yshift=3]\ar[r,yshift=-3]&\displaystyle\prod_{k,l\in J}(X'/_pU)^\sharp\big(Y'_k\times_{Y'}Y'_l\big). 
						\end{tikzcd}
					\end{equation}

					We have that $X'_j/_pU\cong Y'_j$ as presheaves. Let $\zeta_j\in (X'/_pU)^\sharp(Y'_j)$ be elements corresponding to morphisms of presheaves $\zeta_j:Y'_j\to (X'/_pU)^\sharp$ as follows 
					\begin{equation}
						\begin{tikzcd}
							Y'_j\ar[rrr,bend right=15, yshift=-5,"\zeta_j"]\ar[r,"\cong"]&X'_j/_pU\ar[r]&X'/_pU\ar[r]&(X'/_pU)^\sharp
						\end{tikzcd}
					\end{equation}
					where $X'_j/_pU\to X'/_pU$ is induced by $X'_j\to X'$, which is the base change of $Y'_j\to Y'$. We will show the collection $\{\zeta_j\}_{j\in J}$ is a descent datum with respect to $\{Y'_j\to Y'\}_{j\in J}$. 

					Choose $k,l\in J$. Denote projections and base changes of covering morphisms 
					\begin{equation}
						\begin{split}
							\mathrm{pr}_k:&Y'_k\times_{Y'}Y'_k\to Y'_k,\quad \mathrm{pr}_l:Y'_k\times_{Y'}Y'_l\to Y'_l\\
							\iota_k:&X'_k\to X',\quad \iota_l:X'_l\to X'. 
						\end{split}
					\end{equation}
					Since $X'_k/_pU\cong Y'_k$, the quotient $X'_k\to Y'_k$ has sections. Let $\sigma:Y'_k\to X'_k$ be a section. The following diagram commutes 
					\begin{equation}
						\begin{tikzcd}
							Y'_k\times_{Y'}Y'_l\ar[r,"\sigma\circ\mathrm{pr}_k"]\ar[rd,"\mathrm{pr}_k"]\ar[dd,"\mathrm{pr}_l"]&X'_k\ar[r,"\iota_k"]\ar[d]&X'\ar[dd]\\
							&Y'_k\ar[rd]\\
							Y'_l\ar[rr]&&Y'. 
						\end{tikzcd}
					\end{equation}
					Then there is a morphism $\varrho:Y'_k\times_{Y'}Y'_k\to X'_l=X'\times_{Y'}Y'_l$ fitting into the diagram 
					\begin{equation}
						\begin{tikzcd}
							Y'_k\times_{Y'}Y'_l\ar[rd,"\varrho"]\ar[rrd,bend left=20, "{\iota_k\circ\sigma_k\circ\mathrm{pr}_k}"]\ar[rdd,bend right=20,"\mathrm{pr}_l"]\\
							&X'_l\ar[rd,phantom,very near start,"\lrcorner"]\ar[r,"\iota_l"]\ar[d]&X'\ar[d]\\
							&Y'_l\ar[r]&Y'.
						\end{tikzcd}
					\end{equation}
					Let $\sigma,\varrho$ be as above. Then the following diagram of presheaves commutes
					\begin{equation}
						\begin{tikzcd}
							Y'_k\times_{Y'}Y'_l\ar[rrr,"\mathrm{pr}_k"]\ar[rrd,"\sigma\circ\mathrm{pr}_k"]\ar[rdd,"\varrho"]\ar[ddd,"\mathrm{pr}_l"]&&&Y'_k\ar[d,"\cong"]\ar[ddddr,bend left=20,"\zeta_k"]\\
							&&X'_k\ar[r]\ar[d,"\iota_k"]&X'_k/_pU\ar[dd]\\
							&X'_l\ar[r,"\iota_l"]\ar[d]&X'\ar[rd]\\
							Y'_l\ar[r,"\cong"]\ar[rrrrd,bend right =10, "\zeta_l"]&X'_l/_pU\ar[rr]&&X'/_pU\ar[rd]\\
							&&&&(X'/_pU)^\sharp
						\end{tikzcd}
					\end{equation}
					Then for $\zeta_k\in (X'/_pU)^\sharp(Y'_k)$ and $\zeta_l\in(X'/_pU)^\sharp(Y'_l)$ we have 
					\begin{equation}
						\mathrm{pr}_k^*(\zeta_k)=\mathrm{pr}_l^*(\zeta_l)\in (X'/_pU)^\sharp\big(Y'_k\times_{Y'}Y'_l\big). 
					\end{equation}
					This holds for all $k,l\in J$, so $\{\zeta_j\}_{j\in J}$ is a descent datum. Since $(X'/_pU)^\sharp$ is a sheaf, there is a unique $\zeta\in (X'/_pU)^\sharp(Y')$ such that $\zeta|_{Y'_j}=\zeta_j$. This element $\zeta$ does not depend on the covering $\{Y'_j\to Y'\}_{j\in J}$, because $\zeta$ does not change when the covering is refined and any two coverings have a common refinement. The element $\zeta\in (X'/_pU)^\sharp(Y')$ is equivalent to a morphism $\zeta:Y'\to (X'/_pU)^\sharp$ of sheaves by the Yoneda lemma. We omit the verification that $\zeta$ is the inverse of $(X'/_pU)^\sharp\to Y'$. 
				\end{proof}
			}
		}

		\subsubsection{Gluing affine quotients}
		{
			We can glue the affine quotients into a quasi-projective quotient. 
			{
				\begin{theorem}\label{theorem: quotient for projective with UU}
					For the linearisation $\hat U\curvearrowright(X,L)$, assume Condition \hyperref[condition: UU]{UU}. Then $X^0_{\min}$ admits a universal geometric quotient by $U$ in $\mathrm{Sch}/\Bbbk$
					\begin{equation}
						\pi:X^0_{\min}\to Y. 
					\end{equation}
					and $\pi$ is locally isomorphic to $\mathbb A^{r-k}\times_\Bbbk Y$ over $Y$, where $r=\dim U$ and $k\in\mathbb N_+$ is the rank of $\mathrm{coker}(\phi)|_{X^0_{\min}}$. 

					Moreover there exists an ample line bundle $M$ on $Y$ such that $\pi^*M\cong L^m|_{X^0_{\min}}$ for some $m\in\mathbb N_+$. 
				\end{theorem}
				\begin{proof}
					Recall that the linearisation $\hat U\curvearrowright (X,L)$ satisfies Condition \hyperref[condition: UU]{UU}. At the end of Section \ref{subsubsection: reduction to affine case}, we have achieved the situation where: 
					\begin{itemize}
						\item There exists $m\in\mathbb N_+$ and $f_1,\cdots,f_t\in H^0(X,L^m)_{\lambda=\max}$ such that $X^0_{\min}$ is covered by $\{X_{f_i}\}_{i=1}^t$; 
						\item For each $X_{f_i}$, there exists a linear subspace $\mathfrak s_i^*\subseteq \mathfrak u^*$ such that $\mathfrak s_i^*\otimes_\Bbbk\mathcal O_{X_{f_i}}\cong \mathrm{coker}(\phi)|_{X_{f_i}}$ 
						\begin{equation}
							\begin{tikzcd}
								&\mathfrak s_i^*\otimes_\Bbbk\mathcal O_{X_{f_i}}\ar[d,sloped,"\subseteq"]\ar[rd,"\cong"]\\
								\Omega_{X_{f_i}/\Bbbk}\ar[r,"\phi|_{X_{f_i}}"]&\mathfrak u^*\otimes_\Bbbk\mathcal O_{X_{f_i}}\ar[r]&\mathrm{coker}(\phi)|_{X_{f_i}}\ar[r]&0. 
							\end{tikzcd}
						\end{equation}
					\end{itemize}

					We can apply Proposition \ref{proposition: affine GIT with UU} to each of $X_{f_i}$. For each $1\leq i\leq i\leq t$, there exists a universal geometric quotient $\pi_i:X_{f_i}\to Y_i:=\mathrm{Spec}\big(\mathcal O(X_{f_i})^U\big)$. 

					We first glue the affine schemes $Y_i$ for $1\leq i\leq t$. For $1\leq i,j\leq t$, consider $f_j/f_i\in\mathcal O(X_{f_i})$. We have $f_j/f_i\in \mathcal O(X_{f_i})^U$. Since $\pi_i:X_{f_i}\to Y_i$ is a geometric quotient by $U$, we have $\pi_i^*:\mathcal O(Y_i)\cong\mathcal O(X_{f_i})^U$. Then there exists a unique $g_{ji}\in \mathcal O(Y_i)$ such that $\pi_i^*(g_{ji})=f_j/f_i$. Let $Y_{ij}=Y_i\setminus\mathbb V(g_{ji})$ be a principal affine open subscheme. Then $\pi_i^{-1}(Y_{ij})=X_{f_if_j}=\pi_j^{-1}(Y_{ji})$. Since $\pi_i,\pi_j$ are universal geometric quotients, both $Y_{ij}$ and $Y_{ji}$ are geometric quotients of $X_{f_if_j}$. Therefore canonically $Y_{ij}\cong Y_{ji}$ and this isomorphism is the unique morphism such that the following diagram commutes
					\begin{equation}
						\begin{tikzcd}
							X_{f_if_j}\ar[r,equal]\ar[d,"\pi_i|_{X_{f_if_j}}"]&X_{f_jf_i}\ar[d,"\pi_j|_{X_{f_jf_i}}"]\\
							Y_{ij}\ar[r,swap,"\cong"]&Y_{ji}. 
						\end{tikzcd}
					\end{equation}
					It follows easily that $\{Y_i\}_{i=1}^t$ glue to a scheme $Y$ containing each $Y_i$ as an affine open subscheme. Moreover the morphisms $\pi_i:X_{f_i}\to Y_i$ glue to a morphism 
					\begin{equation}
						\pi:X^0_{\min}\to Y
					\end{equation}
					which is an affine universal geometric quotient by $U$. 

					We then construct an ample line bundle $M$ on $Y$ such that $\pi^*M\cong L^m|_{X^0_{\min}}$, where $m\in\mathbb N_+$ is the integer such that $f_i\in H^0(X,L^m)_{\lambda=\max}$. Recall functions $g_{ji}\in\mathcal O(Y_i)$ are invertible on $Y_{ij}$, i.e. $g_{ji}|_{Y_{ij}}\in\mathcal O(Y_{ij})^\times$. They satisfy the cocycle condition $g_{ij}g_{jk}g_{ki}=1$ in $\mathcal O(Y_i\cap Y_j\cap Y_k)$ for all $1\leq i,j,k\leq t$. Then there exists a line bundle $M$ whose transition functions are $\{g_{ij}\}_{i,j=1}^t$. It follows that $\pi^*M$ and $L^m|_{X^0_{\min}}$ both have transition functions $\{f_j/f_i\}_{i,j=1}^t$ with respect to the covering $\{X_{f_i}\}_{i=1}^t$. Therefore $\pi^*M\cong L^m|_{X^0_{\min}}$. 

					It remains to show that $M$ is ample. The scheme $Y=\bigcup_{i=1}^tY_i$ is quasi-compact since it is a finite union of affine schemes. We will construct sections $g^{(i)}\in\Gamma(Y,M)$ for $1\leq i\leq t$ such that each $Y_{g^{(i)}}=Y_i$ is affine. Fix $i$. Consider $g^{(i)}_j:=g_{ji}\in\mathcal O(Y_j)$ for $1\leq j\leq t$. Then $\{g^{(i)}_j\}_{j=1}^t$ define a global section $g^{(i)}\in\Gamma(Y,M)$ since 
					\begin{equation}
						g^{(i)}_{j_1}=g_{j_1j_2}g^{(i)}_{j_2},\quad \textrm{ for }1\leq j_1,j_2\leq t. 
					\end{equation}
					For each $j$, we have 
					\begin{equation}
						Y_j\cap Y_{g^{(i)}}=(Y_j)_{g^{(i)}_j}=(Y_j)_{g_{ji}}=Y_{ji}=Y_j\cap Y_i
					\end{equation}
					which proves that $Y_{g^{(i)}}=Y_i$. This proves that $M$ is ample (\href{https://stacks.math.columbia.edu/tag/01PS}{Definition 01PS}). 
				\end{proof}
			}

			There is a more concrete description of an immersion of the quotient $Y=X^0_{\min}/U$ into a projective space. 
			{
				\begin{lemma}\label{lemma: U-quotient embeds in P^N with UU}
					Assume Condition \hyperref[condition: UU]{UU}. Let $\pi:X^0_{\min}\to Y$ be the quotient in Theorem \ref{theorem: quotient for projective with UU}. Then for sufficiently divisible $m\in\mathbb N_+$, there exists a commutative diagram 
					\begin{equation}
						\begin{tikzcd}
							X\ar[r,hook]&\mathbb P\big(H^0(X,L^m)^*\big)\\
							X^0_{\min}\ar[r,hook]\ar[u,sloped,"\subseteq"]\ar[d,"\pi"]&\mathbb P\big(H^0(X,L^m)^*\big)^0_{\min}\ar[d]\ar[u,sloped,"\subseteq"]\\
							Y\ar[r,hook]&\mathbb P\Big(\big(H^0(X,L^m)^U\big)^*\Big)^0_{\min}
						\end{tikzcd}
					\end{equation}
				\end{lemma}
				\begin{proof}
					According to the proof of Theorem \ref{theorem: quotient for projective with UU}, there exists $m_1\in\mathbb N_+$ such that: 
					\begin{itemize}
						\item $H^0(X,L^{m_1})$ generates $\bigoplus_{d\geq0}H^0(X,L^{dm_1})$; 
						\item there exist $t\in\mathbb N_+$ and $f_1,\cdots,f_t\in H^0(X,L^{m_1})_{\lambda=\max}$; 
						\item $X^0_{\min}=\bigcup_{i=1}^tX_{f_i}$; 
						\item for each $1\leq i\leq t$, the invariant subring $\mathcal O(X_{f_i})^U$ is finitely generated. 
					\end{itemize}
					Then there exists $m_2\in\mathbb N_+$ such that for each $1\leq i\leq t$, we can choose $n_i\in\mathbb N_+$ generators of $\mathcal O(X_{f_i})^U$ of the form 
					\begin{equation}
						\frac{g^{(i)}_j}{f_i^{m_2}}\in\mathcal O(X_{f_i})^U\quad j=1,\cdots,n_i
					\end{equation}
					for $g^{(i)}_j\in H^0(X,L^{m_1m_2})$. Since these finitely many generators $\frac{g^{(i)}_j}{f_i^{m_2}}$ are invariant under $\mathfrak u$ and $\dim \mathfrak u<\infty$, there exists $m_3\in\mathbb N_+$ such that for all $i,j$ we have $f_i^{m_3}g^{(i)}_j\in H^0(X,L^{m_1(m_2+m_3)})^U$. 

					Let $m=m_1(m_2+m_3)\in\mathbb N_+$. The above shows that for each $1\leq i\leq t$, there is a surjective map 
					\begin{equation}
						\frac{H^0(X,L^m)^U}{f_i^{m_2+m_3}}\to \mathcal O(X_{f_i})^U
					\end{equation}
					and then the following ring map is surjective for $1\leq i\leq t$
					\begin{equation}
						\Big(\mathrm{Sym}^\bullet H^0(X,L^m)^U\Big)_{(f_i^{m_2+m_3})}\to \mathcal O(X_{f_i})^U
					\end{equation}
					which corresponds to a closed immersion 
					\begin{equation}
						X_{f_i}/U\hookrightarrow \mathbb P\Big(\big(H^0(X,L^m)^U\big)^*\Big)_{f_i^{m_2+m_3}}. 
					\end{equation}

					By the construction of $\pi:X^0_{\min}\to Y$, we have that $Y$ is glued from $X_{f_i}/U$. The above closed immersions for $1\leq i\leq t$ glue to a closed immersion 
					\begin{equation}
						Y\hookrightarrow \mathbb P\Big(\big(H^0(X,L^m)^U\big)^*\Big)^0_{\min}. 
					\end{equation}
					The closed immersions $X\hookrightarrow \mathbb P\big(H^0(X,L^m)^*\big)$ and $X^0_{\min}\hookrightarrow \mathbb P\big(H^0(X,L^m)^*\big)^0_{\min}$ are naturally defined. It is easy to see these three closed immersions fit into the commutative diagram as stated. 
				\end{proof}
			}

			{
				\begin{corollary}\label{corollary: geometric U-quotients represent quotient sheaves with UU}
					For the linearisation $\hat U\curvearrowright (X,L)$, assume Condition \hyperref[condition: UU]{UU}. Let $\pi:X^0_{\min}\to Y$ be the quotient in Theorem \ref{theorem: quotient for projective with UU}. Then $Y$ represents the quotient sheaf $(X^0_{\min}/_pU)^\sharp$ on $(\mathrm{Sch}/\Bbbk)_{\acute etale}$ (\href{https://stacks.math.columbia.edu/tag/03BD}{Definition 03BD}). 

					Moreover, for any morphism $Y'\to Y$ in $\mathrm{Sch}/\Bbbk$, the quotient sheaf $\big(X'/_pU\big)^\sharp$ is represented by $Y'$, where $X':=X^0_{\min}\times_YY'$. 
				\end{corollary}
				\begin{proof}
					According to the proof of Theorem \ref{theorem: quotient for projective with UU}, the quotient $\pi:X^0_{\min}\to Y$ is glued from $\pi|_{X_{f_i}}:X_{f_i}\to Y_i:=X_{f_i}/U$ for some $m\in\mathbb N_+$ and $f_1,\cdots,f_t\in H^0(X,L^m)_{\lambda=\max}$. Let $\{Y_j\to Y\}_{j\in J}$ denote this Zariski covering. which is also a covering in $(\mathrm{Sch}/\Bbbk)_{\acute etale}$. Let $\{X_j\to X^0_{\min}\}_{j\in J}$ be the base change covering of $\{Y_j\to Y\}$. By Corollary \ref{corollary: affine quotient represents quotient sheaf with UU}, we have $(X_j/_pU)^\sharp\cong Y_j$ for all $j\in J$, and the representability is stable under base change, in the sense that $\big(X_j\times_{Y_j}W\big/_pU\big)^\sharp\cong W$ for any morphism $W\to Y_j$ in $\mathrm{Sch}/\Bbbk$. 

					Let $Y'\to Y$ be a morphism in $\mathrm{Sch}/\Bbbk$. Denote 
					\begin{equation}
						X':=X\times_YY',\quad Y'_j:=Y_j\times_YY',\quad X'_j:=X'\times_{Y'}Y'_j. 
					\end{equation}
					Then $X'\to Y'$ is a geometric quotient as a base change of a universal geometric quotient. The quotient induces a morphism $(X'/_pU)^\sharp\to Y'$ of sheaves. We will construct an element $\zeta\in (X'/_pU)^\sharp(Y')$ such that the corresponding morphism $Y'\to (X'/_pU)^\sharp$ is inverse to $(X'/_pU)^\sharp\to Y'$. 

					By Corollary \ref{corollary: affine quotient represents quotient sheaf with UU}, we have $(X'_j/_pU)^\sharp\cong Y'_j$, which is the base change of $(X_j/_pU)^\sharp\cong Y_j$ along $Y'_j\to Y_j$. The morphism $X'_j\to (X'_j/_pU)^\sharp$ is surjective as a morphism of sheaves, so there exists a covering $\{\varphi_j:\widetilde{Y'_j}\to Y'_j\}$ in $(\mathrm{Sch}/\Bbbk)_{\acute etale}$ and a morphism $\sigma_j:\widetilde{Y'_j}\to X'_j$ fitting into the diagram 
					\begin{equation}
						\begin{tikzcd}
							\widetilde{Y'_j}\ar[r,"\sigma_j"]\ar[d,"\varphi_j"]&X'_j\ar[d]\\
							Y'_j\ar[r,"\cong"]&(X'_j/_pU)^\sharp. 
						\end{tikzcd}
					\end{equation}
					The composition of $\{\varphi_j:\widetilde{Y'_j}\to Y'_j\}$ and $\{Y'_j\to Y'\}$ is also a covering $\{\widetilde{Y'_j}\to Y'\}_{j\in J}$ in $(\mathrm{Sch}/\Bbbk)_{\acute etale}$. Since $(X'/_pU)^\sharp$ is a sheaf, the following is an equaliser diagram 
					\begin{equation}
						\begin{tikzcd}
							(X'/_pU)^\sharp(Y')\ar[r]&\displaystyle\prod_{j\in J}(X'/_pU)^\sharp\big(\widetilde{Y'_j}\big)\ar[r,yshift=3]\ar[r,yshift=-3]&\displaystyle\prod_{k,l\in J}(X'/_pU)^\sharp\big(\widetilde{Y'_k}\times_{Y'}\widetilde{Y'_l}\big). 
						\end{tikzcd}
					\end{equation}
					Let $\zeta_j\in (X'/_pU)^\sharp\big(\widetilde{Y'_j}\big)$ be the element corresponding to the morphism $\widetilde{Y'_j}\xrightarrow{\varphi_j} Y'_j\cong (X'_j/_pU)^\sharp\to (X'/_pU)^\sharp$. We will show $\{\zeta_j\}_{j\in J}$ is a descent datum with respect to $\{\widetilde{Y'_j}\to Y'\}_{j\in J}$. 

					Fix $k,l\in J$. The following diagram commutes
					\begin{equation}
						\begin{tikzcd}
							\widetilde{Y'_k}\times_{Y'}\widetilde{Y'_l}\ar[r,"\widetilde{\mathrm{pr}_k}"]\ar[d,"\widetilde{\mathrm{pr}_l}"]&\widetilde{Y'_k}\ar[r,"\sigma_k"]\ar[d,"\varphi_k"]&X'_k\ar[r,"\iota_k"]\ar[d]&X'\ar[d]\\
							\widetilde{Y'_l}\ar[d,"\varphi_l"]&Y'_k\ar[r,"\cong"]\ar[rrd]&(X'_k/_pU)^\sharp\ar[r]&(X'/_pU)^\sharp\ar[d]\\
							Y'_l\ar[rrr]&&&Y'. 
						\end{tikzcd}
					\end{equation}
					There exists a unique morphism $\varrho:\widetilde{Y'_k}\times_{Y'}\widetilde{Y'_l}\to X'_l=X'\times_{Y'}Y'_l$ fitting into the diagram 
					\begin{equation}
						\begin{tikzcd}
							\widetilde{Y'_k}\times_{Y'}\widetilde{Y'_l}\ar[rrd,bend left=20,"\iota_k\circ\sigma_k\circ\widetilde{\mathrm{pr}_k}"]\ar[rd,"\varrho"]\ar[rdd,bend right=20,swap,"\varphi_l\circ \widetilde{\mathrm{pr}_l}"]\\
							&X'_l\ar[r,"\iota_l"]\ar[rd,phantom,very near start,"\lrcorner"]\ar[d]&X'\ar[d]\\
							&Y'_l\ar[r]&Y'. 
						\end{tikzcd}
					\end{equation}
					The following diagram commutes
					\begin{equation}
						\begin{tikzcd}
							\widetilde{Y'_k}\times_{Y'}\widetilde{Y'_l}\ar[rrr,"\widetilde{\mathrm{pr}_k}"]\ar[ddd,"\widetilde{\mathrm{pr}_l}"]\ar[rrd,"\sigma_k\circ\widetilde{\mathrm{pr}_k}"]\ar[rdd,"\varrho"]&&&\widetilde{Y'_k}\ar[ld,"\sigma_k"]\ar[d,"\varphi_k"]\ar[ddd,bend left=20,"\zeta_k",xshift=15]\\
							&&X'_k\ar[d,"\iota_k"]\ar[rd]&Y'_k\ar[d,"\cong"]\\
							&X'_l\ar[r,"\iota_l"]\ar[rd]&X'\ar[rd]&(X'_k/_pU)^\sharp\ar[d]\\
							\widetilde{Y'_l}\ar[r,"\varphi_l"]\ar[rrr,bend right =10, yshift=-10,"\zeta_l"]&Y'_l\ar[r,"\cong"]&(X'_l/_pU)^\sharp\ar[r]&(X'/_pU)^\sharp. 
						\end{tikzcd}
					\end{equation}
					Therefore 
					\begin{equation}
						\widetilde{\mathrm{pr}_k}^*(\zeta_k)=\widetilde{\mathrm{pr}_l}^*(\zeta_l)\in (X'/_pU)^\sharp\big(\widetilde{Y'_k}\times_{Y'}\widetilde{Y'_l}\big)
					\end{equation}
					which holds for all $k,l\in J$. Then $\{\zeta_j\}_{j\in J}$ is a descent datum. There exists a unique $\zeta\in (X'/_pU)^\sharp(Y')$ such that $\zeta\big|_{\widetilde{Y'_j}}=\zeta_j$. We omit the verification that $\zeta$ does not depend on the covering $\{\widetilde{Y'_j}\to Y'\}_{j\in J}$ and that the corresponding morphism $\zeta:Y'\to (X'/_pU)^\sharp$ is inverse to $(X'/_pU)^\sharp \to Y'$. 
				\end{proof}
			}
		}

		\subsubsection{Quotients by \texorpdfstring{$\hat U$}{U-hat} and parabolic groups}
		{
			With Condition \hyperref[condition: UU]{UU}, there is a quasi-projective geometric quotient $\pi:X^0_{\min}\to X^0_{\min}/U$ by $U$. To obtain a quotient by $\hat U$, we consider the residual linear action of $\lambda$ on $X^0_{\min}/U$. 
			{
				\begin{theorem}[\cite{BércziGergely2018Gitf}]\label{theorem: projective geometric quotient by U-hat with UU}
					Let $X$ be a projective scheme over $\Bbbk$ with an ample line bundle $L$. Let $\hat U$ act on $X$ linearly with respect to $L$. Assume Condition \hyperref[condition: UU]{UU} for $\hat U\curvearrowright(X,L)$. Then the subset $X^0_{\min}\setminus UZ_{\min}$ is $\hat U$-invariant open, and it admits a projective universal geometric quotient by $\hat U$. 
				\end{theorem}
				\begin{remark}
					Condition \hyperref[condition: UU]{UU} implies $X^0_{\min}\ne\emptyset$. However it is possible that $X^0_{\min}\setminus UZ_{\min}=\emptyset$. 
				\end{remark}
				\begin{proof}
					By Theorem \ref{theorem: quotient for projective with UU}, we have a universal geometric quotient by $U$ 
					\begin{equation}
						\pi:X^0_{\min}\to Y. 
					\end{equation}
					By Lemma \ref{lemma: U-quotient embeds in P^N with UU}, there exists $m\in\mathbb N_+$ such that for $W:=H^0(X,L^m)^U$, we have an immersion $Y\to \mathbb P(W^*)$. Consider $Y$ as a locally closed subscheme of $\mathbb P(W^*)$. Let $\overline Y\subseteq \mathbb P(W^*)$ be the scheme theoretic image of $Y\subseteq \mathbb P(W^*)$. Then $\overline Y\subseteq \mathbb P(W^*)$ is $\lambda$-invariant and $Y=\overline Y\cap \mathbb P(W^*)^0_{\min}$. 

					There is a natural linearisation $\lambda\curvearrowright\big(\mathbb P(W^*),\mathcal O(1)\big)$, induced from the linear representation $\lambda\curvearrowright H^0(X,L^m)^U$. We will consider a twisted of the linearisation $\lambda\curvearrowright\big(\overline Y,\mathcal O(1)|_{\overline Y}\big)$. 

					Let $\max(\lambda, W),\mathrm{2nd}\max(\lambda,W)\in\mathbb Z$ denote the maximal and next maximal weights of $\lambda\curvearrowright W$. Let $\theta\in\mathbb Q$ be such that 
					\begin{equation}
						\mathrm{2nd}\max(\lambda,W)<-\theta<\max(\lambda, W). 
					\end{equation}
					Consider $\theta\in\mathbb Q$ as a rational character of $\mathbb G_m$. Denote $\mathcal O(1)_\theta$ the linearisation twisted by $\theta$. Then we have 
					\begin{equation}
						\mathbb P(W^*)^{\lambda,\mathcal O(1)_\theta,\mathrm{ss}}=\mathbb P(W^*)^{\lambda,\mathcal O(1)_\theta,\mathrm{s}}=\mathbb P(W^*)^0_{\min}\big\backslash \mathbb P\big((W_{\lambda=\max})^*\big). 
					\end{equation}
					Then by \cite{MumfordDavid1994Git} Theorem 1.19
					\begin{equation}
						\overline Y^{\lambda,\mathcal O(1)_\theta,\mathrm{(s)s}}=\overline Y\cap \mathbb P(W^*)^{\lambda,\mathcal O(1)_\theta,\mathrm{(s)s}}=Y\cap \Big(\mathbb P(W^*)^0_{\min}\big\backslash \mathbb P\big((W_{\lambda=\max})^*\big)\Big). 
					\end{equation}
					There is a projective universal geometric quotient by $\lambda$ 
					\begin{equation}
						\overline Y^{\lambda,\mathcal O(1)_\theta,\mathrm{(s)s}}\to \overline Y/\!/\lambda. 
					\end{equation}

					The subset $X^0_{\min}\setminus UZ_{\min}\subseteq X^0_{\min}$ is open since it is the preimage of $\overline Y^{\lambda,\mathcal O(1)_\theta,\mathrm{(s)s}}\subseteq Y$ along $\pi$
					\begin{equation}
						X^0_{\min}\setminus UZ_{\min}=\pi^{-1}\bigg(Y\cap \Big(\mathbb P(W^*)^0_{\min}\big\backslash \mathbb P\big((W_{\lambda=\max})^*\big)\Big)\bigg). 
					\end{equation}
					Then $X^0_{\min}\setminus UZ_{\min}\to \overline Y^{\lambda,\mathcal O(1)_\theta,\mathrm{(s)s}}$ is a universal geometric quotient by $U$ as a base change of $\pi:X^0_{\min}\to Y$. The composition 
					\begin{equation}
						X^0_{\min}\setminus UZ_{\min}\to \overline Y^{\lambda,\mathcal O(1)_\theta,\mathrm{(s)s}}\to \overline Y/\!/\lambda
					\end{equation}
					is then a projective universal geometric quotient by $\hat U$. 
				\end{proof}
			}

			We then consider a special type of parabolic group. Let $n_1,n_2\in\mathbb N_+$. Consider the parabolic subgroup $P\subseteq\mathrm{SL}(n_1+n_2)$ 
			\begin{equation}
				P=\bigg\{\begin{pmatrix}g_{11}&g_{12}\\&g_{22}\end{pmatrix}:\begin{matrix}g_{ii}\in\mathrm{GL}(n_i)\textrm{ for }i=1,2\\\det(g_{11})\det(g_{22})=1\end{matrix}\bigg\}. 
			\end{equation}
			Let $U\subseteq P$ be the unipotent radical, $T=\mathbb G_m$ and $K=\mathrm{SL}(n_1)\times\mathrm{SL}(n_2)$. Then there is a surjective homomorphism with finite kernel 
			\begin{equation}
				U\rtimes(T\times K)\to P,\quad (u,t,k)\mapsto u\begin{pmatrix}t^{\frac{n_2}{(n_1,n_2)}}I_{n_1}&\\&t^{-\frac{n_1}{(n_1,n_2)}}I_{n_2}\end{pmatrix}k
			\end{equation}
			where $(n_1,n_2):=\gcd(n_1,n_2)$. 

			Let $X$ be a projective scheme with an ample line bundle $L$. Let $P$ act on $X$ linearly with respect to $L$. We can consider the linearisation $U\rtimes(T\times K)\curvearrowright(X,L)$ via the above homomorphism onto $P$. In $U\rtimes(T\times K)$, the subgroup $U\rtimes T$ has the form $\hat U$. We can first apply our theorems for $\hat U$-linearisations and then consider the residual linearisation by $K$. 

			{
				\begin{theorem}\label{theorem: geometric quotients by parabolic with UU}
					Let $U\rtimes(T\times K)\to P$ be as above. Let $X$ be a projective scheme with an ample line bundle $L$. Let $P$ act on $X$ linearly with respect to $L$. Assume Condition \hyperref[condition: UU]{UU} for $\hat U\curvearrowright(X,L)$, where $\hat U:=U\rtimes T$. Let $Z_{\min}^{K,\mathrm{(s)s}}\subseteq Z_{\min}$ denote the $K$-(semi)stable locus of $K\curvearrowright \big(Z_{\min},L|_{Z_{\min}}\big)$. Let $p_\lambda:X^0_{\min}\to Z_{\min}$ be the retraction. Then the following $P$-invariant open subset 
					\begin{equation}
						\big(X^0_{\min}\setminus UZ_{\min}\big)\cap p_\lambda^{-1}\big(Z_{\min}^{K,\mathrm{s}}\big)
					\end{equation}
					has a quasi-projective universal geometric quotient by $P$. Moreover if $Z_{\min}^{K,\mathrm{s}}=Z_{\min}^{K,\mathrm{ss}}$, then the quotient is projective. 
				\end{theorem}
				\begin{proof}
					Assume $Z_{\min}^{K,\mathrm{s}}\ne\emptyset$. Let $\pi:X^0_{\min}\to Y$ be the quotient in Theorem \ref{theorem: quotient for projective with UU}. By Lemma \ref{lemma: U-quotient embeds in P^N with UU}, there exists $m\in\mathbb N_+$ such that for $W:=H^0(X,L^m)^U$, we can view $Y$ as a locally closed subscheme $Y\subseteq\mathbb P(W^*)$. Let $\overline Y\subseteq \mathbb P(W^*)$ be the scheme theoretic image of $Y\subseteq \mathbb P(W^*)$. Then $\overline Y\subseteq \mathbb P(W^*)$ is a $T\times K$-invariant closed subscheme. 

					Characters of $T\times K$ are characters of $T=\mathbb G_m$. Rational characters are represented by rational numbers. Let $\max(\lambda,W)\in\mathbb Z$ be the maximal weight of $\lambda\curvearrowright W$. Let $\epsilon\in\mathbb Q_{>0}$. Let $\theta\in\mathbb Q$ represents a rational character of $T\times K$ such that 
					\begin{equation}
						\max(\lambda,W)-\epsilon<-\theta<\max(\lambda,W). 
					\end{equation}

					Let $\mathbb P(W^*)^{T\times K,\widehat{\mathrm{s}}}\subseteq \mathbb P(W^*)^{T\times K,\widehat{\mathrm{ss}}}\subseteq\mathbb P(W^*)$ denote open subsets whose closed points are the following respectively 
					\begin{equation}
						\begin{split}
							\mathbb P(W^*)^{T\times K,\widehat{\mathrm{s}}}\overset{\mathrm{cl}}{=\!=}&\Big\{[v]\in\mathbb P(W^*):v\ne v_{\lambda=\max}\ne 0,\;[v_{\lambda=\max}]\in\mathbb P\big((W^*)_{\lambda=\max}\big)^{K,\mathrm{s}}\Big\}\\
							\mathbb P(W^*)^{T\times K,\widehat{\mathrm{ss}}}\overset{\mathrm{cl}}{=\!=}&\Big\{[v]\in\mathbb P(W^*):v\ne v_{\lambda=\max}\ne 0,\;[v_{\lambda=\max}]\in\mathbb P\big((W^*)_{\lambda=\max}\big)^{K,\mathrm{ss}}\Big\}
						\end{split}
					\end{equation}
					where $\overset{\mathrm{cl}}{=\!=}$ denotes that two sides have the same closed points. 

					When $\epsilon>0$ is sufficiently small, for the twisted linearisation $T\times K\curvearrowright\big(\mathbb P(W^*),\mathcal O(1)_\theta\big)$, we have by the Hilbert-Mumford criterion 
					\begin{equation}
						\mathbb P(W^*)^{T\times K,\widehat{\mathrm{s}}}\subseteq \mathbb P(W^*)^{T\times K,\mathcal O(1)_\theta,\mathrm{s}}\subseteq \mathbb P(W^*)^{T\times K,\mathcal O(1)_\theta,\mathrm{ss}}\subseteq\mathbb P(W^*)^{T\times K,\widehat{\mathrm{ss}}}. 
					\end{equation}

					Note that the four subsets above are all contained in $\mathbb P(W^*)^0_{\min}$. Recall that $Y=\overline Y\cap \mathbb P(W^*)^0_{\min}$. Then their intersections with $\overline Y$ form the following chain 
					\begin{equation}\label{equation: chain of hat-stable/stable/semistable/hat-semistable}
						Y\cap \mathbb P(W^*)^{T\times K,\widehat{\mathrm{s}}}\subseteq \overline Y^{T\times K,\mathcal O(1)_\theta,\mathrm{s}}\subseteq \overline Y^{T\times K,\mathcal O(1)_\theta,\mathrm{ss}}\subseteq Y\cap \mathbb P(W^*)^{T\times K,\widehat{\mathrm{ss}}}. 
					\end{equation}

					The stable locus $\overline Y^{T\times K,\mathcal O(1)_\theta,\mathrm{s}}$ has a quasi-projective universal geometric quotient by $T\times K$, and so does the $T\times K$-invariant open subset $Y\cap \mathbb P(W^*)^{T\times K,\widehat{\mathrm{s}}}$. The preimage of $Y\cap \mathbb P(W^*)^{T\times K,\widehat{\mathrm{s}}}$ along $\pi$ is 
					\begin{equation}
						\pi^{-1}\big(Y\cap \mathbb P(W^*)^{T\times K,\widehat{\mathrm{s}}}\big)=\big(X^0_{\min}\setminus UZ_{\min}\big)\cap p_\lambda^{-1}\big(Z_{\min}^{K,\mathrm{s}}\big). 
					\end{equation}
					Therefore the composition 
					\begin{equation}
						\big(X^0_{\min}\setminus UZ_{\min}\big)\cap p_\lambda^{-1}\big(Z_{\min}^{K,\mathrm{s}}\big)\to Y\cap \mathbb P(W^*)^{T\times K,\widehat{\mathrm{s}}}\to \big(Y\cap \mathbb P(W^*)^{T\times K,\widehat{\mathrm{s}}}\big)\big/T\times K
					\end{equation}
					is a quasi-projective universal geometric quotient by $U\rtimes(T\times K)$, the same as a quotient by $P$. 

					If $Z_{\min}^{K,\mathrm{s}}=Z_{\min}^{K,\mathrm{ss}}$, then we have
					\begin{equation}
						Y\cap \mathbb P(W^*)^{T\times K,\widehat{\mathrm{s}}}=Y\cap \mathbb P(W^*)^{T\times K,\widehat{\mathrm{ss}}}
					\end{equation}
					since 
					\begin{equation}
						Y\cap \mathbb P\big((W^*)_{\lambda=\max}\big)^{K,\mathrm{s}}=\pi\big(Z_{\min}^{K,\mathrm{s}}\big)=\pi\big(Z_{\min}^{K,\mathrm{s}}\big)=Y\cap \mathbb P\big((W^*)_{\lambda=\max}\big)^{K,\mathrm{ss}}. 
					\end{equation}
					Then four sets in the chain \eqref{equation: chain of hat-stable/stable/semistable/hat-semistable} are the same. In this situation the $P$-quotient is projective 
					\begin{equation}
						\Big(\big(X^0_{\min}\setminus UZ_{\min}\big)\cap p_\lambda^{-1}\big(Z_{\min}^{K,\mathrm{s}}\big)\Big)\Big/ P\cong \overline Y\big/\!\big/_\theta(T\times K). 
					\end{equation}
				\end{proof}
			}
		}
	}

	\subsection{Theorem with Condition \texorpdfstring{\hyperref[condition: WUU]{WUU}}{WUU}}
	{
		For the linearisation $\hat U\curvearrowright (X,L)$, Condition \hyperref[condition: WUU]{WUU} is weaker than Condition \hyperref[condition: UU]{UU}. When Condition \hyperref[condition: WUU]{WUU} is satisfied, we can show that there exists a $\hat U$-equivariant blow-up $\widetilde X\to X$ and an ample linearisation $\hat U\curvearrowright(\widetilde X,\widetilde L)$, such that Condition \hyperref[condition: UU]{UU} is satisfied. 

		\subsubsection{The blow-up centre}
		{
			Let $\iota:Z_{\min}\hookrightarrow X$ denote the closed immersion and let $\iota^\sharp:\mathcal O_X\to \iota_*\mathcal O_{Z_{\min}}$ denote the associated morphism between sheaves of rings. Then $\ker(\iota^\sharp)\subseteq\mathcal O_X$ is the sheaf of ideals associated to $Z_{\min}\hookrightarrow X$. Recall that $k\in\mathbb N$ is the integer such that $\mathrm{Fit}_{k-1}(\phi)=0$ and $\dim\mathrm{Stab}_U(z)=k$ for some closed point $z\in Z_{\min}$. Let $\mathcal I\subseteq\mathcal O_X$ be the sheaf of ideals 
			\begin{equation}
				\mathcal I:=\ker(\iota^\sharp)+\mathrm{Fit}_k(\phi)\subseteq\mathcal O_X. 
			\end{equation}
			Let $\mathbb V(\mathcal I)\hookrightarrow X$ denote the closed subscheme associated to $\mathcal I$. Then scheme theoretically $\mathbb V(\mathcal I)=Z_{\min}\cap \mathbb V(\mathrm{Fit}_k(\phi))$. 

			Let $\sigma:U\times_\Bbbk X\to X$ denote the action morphism. Since $U$ is affine, we have $\sigma_*(\mathcal O_{U\times_\Bbbk X})\cong \mathcal O(U)\otimes_\Bbbk\mathcal O_X$. Let $\mathcal J\subseteq \mathcal O_X$ be the kernel of the following composition 
			\begin{equation}
				\begin{tikzcd}
					\mathcal O_X\ar[r,"\sigma^\sharp"]&\sigma_*(\mathcal O_{U\times_\Bbbk X})\ar[r,"\cong"]&\mathcal O(U)\otimes_\Bbbk\mathcal O_X\ar[r,->>]&\mathcal O(U)\otimes_\Bbbk\mathcal O_X/\mathcal I. 
				\end{tikzcd}
			\end{equation}
			The closed subscheme $\mathbb V(\mathcal J)\hookrightarrow X$ is the scheme theoretic image of 
			\begin{equation}
				\begin{tikzcd}
					U\times_\Bbbk\big(Z_{\min}\cap \mathbb V(\mathrm{Fit}_k(\phi))\big)\ar[r,equal]&U\times_\Bbbk \mathbb V(\mathcal I)\ar[r,hook]&U\times_\Bbbk X\ar[r,"\sigma"]&X
				\end{tikzcd}
			\end{equation}
			which is the scheme theoretic version of the closure of the $U$-sweep of $\{z\in Z_{\min}:\dim\mathrm{Stab}_U(z)>k\}$. 

			We will choose $\mathbb V(\mathcal J)\hookrightarrow X$ to be the centre of the blow-up. Condition \hyperref[condition: WUU]{WUU} implies that $X\setminus \mathbb V(\mathcal J)\ne\emptyset$, so $\widetilde X\ne\emptyset$. When Condition \hyperref[condition: UU]{UU} holds, we have $\mathcal I=\mathcal J=\mathcal O_X$ and the centre is empty. Blow-ups are constructed locally. Since we will only consider $X^0_{\min}$, we only need the centre to be a closed subscheme $C\hookrightarrow X^0_{\min}$. 

			Let $m\in\mathbb N_+$ be such that $H^0(X,L^m)$ generates $\bigoplus_{n\in\mathbb N}H^0(X,L^{mn})$. By the existence of $z\in Z_{\min}$ and Lemma \ref{lemma: equivalence of dimStab>k and Fit_k in m_x}, there exists $f\in H^0(X,L^m)_{\lambda=\max}$ such that $\mathrm{Fit}_k(\phi)(X_f)_{\lambda=0}\not\subseteq\mathrm{rad}(\mathcal O(X_f))$. 

			Let $f\in H^0(X,L^m)_{\lambda=\max}$ such that $X_f\ne\emptyset$. Since $X_f$ is $U$-invariant and scheme theoretic images are constructed locally, we have that $J_f:=\mathcal J(X_f)\subseteq\mathcal O(X_f)$ is the kernel of the following composition 
			\begin{equation}
				\begin{tikzcd}
					\mathcal O(X_f)\ar[r,"(\sigma|_{X_f})^*"]&\mathcal O(U)\otimes_\Bbbk\mathcal O(X_f)\ar[r,->>]&\mathcal O(U)\otimes_\Bbbk\mathcal O(X_f)/I_f
				\end{tikzcd}
			\end{equation}
			where $I_f:=\mathcal I(X_f)=\mathcal O(X_f)_{\lambda<0}+\mathrm{Fit}_k(\phi)(X_f)$. 

			Let $\pi:\widetilde X\to X$ be the blow-up. Let $E\hookrightarrow \widetilde X$ be the exceptional divisor. Then $\mathcal O(-E)$ is relatively ample with respect to $\pi$. For $x\ll y\in\mathbb N_+$, let 
			\begin{equation}
				\widetilde L_{x,y}:=\mathcal O(-E)^{\otimes x}\otimes_{\mathcal O_{\widetilde X}}\pi^* L^{\otimes y}
			\end{equation}
			be a line bundle on $\widetilde X$. Since $x\ll y$ we have that $\widetilde L_{x,y}$ is ample by \href{https://stacks.math.columbia.edu/tag/0892}{Lemma 0892}. All constructions are equivariant under $\hat U$, and then there is a natural linearisation 
			\begin{equation}
				\hat U\curvearrowright\big(\widetilde X,\widetilde L_{x,y}\big). 
			\end{equation}

			The \emph{blow-up algebra} associated to the ideal and ring $J_f\subseteq \mathcal O(X_f)$ is the following graded $\mathcal O(X_f)$-algebra 
			\begin{equation}
				\mathrm{Bl}_{J_f}\big((\mathcal O(X_f)\big):=\bigoplus_{n\in\mathbb N}(J_f)^n=\mathcal O(X_f)\oplus J_f\oplus (J_f)^2\oplus\cdots
			\end{equation}
			where the summand $(J_f)^n$ has degree $n$. For $a\in (J_f)^1$ of degree one, the \emph{affine blow-up algebra} is the homogeneous localisation 
			\begin{equation}
				\mathcal O(X_f)\Big[\frac{J_f}{a}\Big]:=\Big(\mathrm{Bl}_{J_f}\big(\mathcal O(X_f)\big)\Big)_{(a)}. 
			\end{equation}
			Denote 
			\begin{equation}
				\widetilde X_{f,a}:=\mathrm{Spec}\Big(\mathcal O(X_f)\Big[\frac{J_f}{a}\Big]\Big)\subseteq \widetilde X. 
			\end{equation}
			Then $\{\widetilde X_{f,a}\}_{f,a}$ form an affine open covering of $\widetilde X$ for $f\in H^0(X,L^m)$ and $a\in J_f$. 
			We have
			\begin{equation}
				\pi^{-1}(X_f)=\bigcup_{a\in J_f}\widetilde X_{f,a}. 
			\end{equation}

			It is easy to see 
			\begin{equation}
				\widetilde X^0_{\min}=\bigcup_{\substack{f\in H^0(X,L^m)_{\lambda=\max}\\a\in (J_f)_{\lambda=0}}}\widetilde X_{f,a}
			\end{equation}
			and we can choose the sub-covering $\big\{\widetilde X_{f,a}\big\}$ for $(f,a)$ such that $a\notin\mathrm{rad}(\mathcal O(X_f))$ and $a\in (J_f)_{\lambda=0}=\mathrm{Fit}_k(\phi)(X_f)_{\lambda=0}$ is of the form 
			\begin{equation}
				a=\det\begin{pmatrix}
					\xi_{i_1}.g_1&\cdots&\xi_{i_1}.g_{r-k}\\
					\vdots&\ddots&\vdots\\
					\xi_{i_{r-k}}.g_1&\cdots&\xi_{i_{r-k}}.g_{r-k}
				\end{pmatrix}
			\end{equation}
			for $1\leq i_1<\cdots<i_{r-k}\leq r$ and $g_1,\cdots,g_{r-k}\in \mathcal O(X_f)_{\lambda=-w}$. The condition that $a\notin \mathrm{rad}(\mathcal O(X_f))$ implies that $\mathrm{Fit}_k(\phi)(X_f)_{\lambda=0}\not\subseteq \mathrm{rad}(\mathcal O(X_f))$. 
		}

		\subsubsection{Main theorem}
		{
			Recall that $U$ is an additive group and the induced $\mathbb G_m$-representation on $\mathfrak u$ has one weight, denoted $w>0$. 

			{
				\begin{theorem}\label{theorem: blowup from WUU to UU}
					For the linearisation $\hat U\curvearrowright (X,L)$, assume Condition \hyperref[condition: WUU]{WUU}. Then the blow-up $\pi:\widetilde X\to X$ and the linearisation $\hat U\curvearrowright\big(\widetilde X,\widetilde L_{x,y}\big)$ for $x\ll y\in\mathbb N_+$ satisfy Condition \hyperref[condition: UU]{UU}. 
				\end{theorem}
				\begin{proof}
					Condition \hyperref[condition: WUU]{WUU} implies that there exists $k\in\mathbb N$ such that $\mathrm{Fit}_{k-1}(\phi)=0$ and $\dim\mathrm{Stab}_U(z)=k$ for some closed point $z\in Z_{\min}$. We will prove that 
					\begin{equation}
						\mathrm{Fit}_{k-1}(\widetilde\phi)|_{\widetilde X^0_{\min}}=0,\quad \mathrm{Fit}_k(\widetilde\phi)|_{\widetilde X^0_{\min}}=\mathcal O_{\widetilde X^0_{\min}}
					\end{equation}
					where $\widetilde\phi:\Omega_{\widetilde X/\Bbbk}\to \mathfrak u^*\otimes_\Bbbk\mathcal O_{\widetilde X}$ is the infinitesimal action of $\mathfrak u$ on $\widetilde X$. 

					We can check the equalities of Fitting ideals locally. Recall that $\widetilde X^0_{\min}$ is covered by $\tilde X_{f,a}$ for $m\in\mathbb N_+$ sufficiently divisible, $f\in H^0(X,L^m)_{\lambda=\max}$ and non-nilpotent $a\in (J_f)_{\lambda=0}$ of the form 
					\begin{equation}
						a=\det\begin{pmatrix}
							\xi_{i_1}.g_1&\cdots&\xi_{i_1}.g_{r-k}\\
							\vdots&\ddots&\vdots\\
							\xi_{i_{r-k}}.g_1&\cdots&\xi_{i_{r-k}}.g_{r-k}
						\end{pmatrix}
					\end{equation}
					for $1\leq i_1<\cdots<i_{r-k}\leq r$ and $g_1,\cdots,g_{r-k}\in \mathcal O(X_f)_{\lambda=-w}$. Without loss of generality, we can assume $(i_1,\cdots,i_{r-k})=(1,\cdots,r-k)$. We need to show 
					\begin{equation}
						\mathrm{Fit}_{k-1}(\widetilde\phi)(\widetilde X_{f,a})=0,\quad \mathrm{Fit}_k(\widetilde\phi)(\widetilde X_{f,a})=\mathcal O(X_f)\Big[\frac{J_f}{a}\Big]. 
					\end{equation}

					There is a surjective map $\mathcal O(\widetilde X_{f,a})\otimes_\Bbbk\mathrm{d}\big(\frac{J_f}{a}\big)\twoheadrightarrow \Omega_{\mathcal O(\widetilde X_{f,a})/\Bbbk}$. Therefore Fitting ideals of $\mathrm{coker}(\widetilde\phi)(\widetilde X_{f,a})$ can be calculated from the following presentation
					\begin{equation}
						\begin{tikzcd}
							\mathcal O(\widetilde X_{f,a})\otimes_\Bbbk\mathrm{d}\big(\frac{J_f}{a}\big)\ar[r]&\mathfrak u^*\otimes_\Bbbk\mathcal O(\widetilde X_{f,a})\ar[r]&\mathrm{coker}(\widetilde\phi)(\widetilde X_{f,a})\ar[r]&0. 
						\end{tikzcd}
					\end{equation}

					For the $(k-1)$th Fitting ideal, we need to show 
					\begin{equation}
						\det\begin{pmatrix}
							\zeta_1.\frac{h_1}{a}&\cdots&\zeta_1.\frac{h_{r-k+1}}{a}\\
							\vdots&\ddots&\vdots\\
							\zeta_{r-k+1}.\frac{h_1}{a}&\cdots&\zeta_{r-k+1}.\frac{h_{r-k+1}}{a}
						\end{pmatrix}=0\quad \textrm{in }\mathcal O(\widetilde X_{f,a})
					\end{equation}
					for any $h_1,\cdots,h_{r-k+1}\in J_f$ and any $\zeta_1,\cdots,\zeta_{r-k+1}\in \mathfrak u$. The above determinant times $a^{r-k+1}$ is an element in $\mathrm{Fit}_{k-1}(\phi)(X_f)=0$. Since $a\in\mathcal O(\widetilde X_{f,a})$ is a non-zerodivisor, the determinant above is zero. This proves $\mathrm{Fit}_{k-1}(\widetilde\phi)(\widetilde X_{f,a})=0$. 

					For the $k$th Fitting ideal, we first observe that for 
					\begin{equation}
						b_i:=\det\begin{pmatrix}
							\xi_1.g_1&\cdots&\xi_1.g_{r-k}\\
							\vdots&\cdots&\vdots\\
							g_1&\cdots&g_{r-k}\\
							\vdots&\cdots&\vdots\\
							\xi_{r-k}.g_1&\cdots&\xi_{r-k}.g_{r-k}
						\end{pmatrix}\in\mathcal O(X_f),\quad 1\leq i\leq r-k
					\end{equation}
					where $(g_1,\cdots,g_{r-k})$ occupies the $i$th row, we have the following equations 
					\begin{equation}
						\xi_i.b_j=\delta_{ij}a,\quad 1\leq i,j\leq r-k. 
					\end{equation}
					We also have $b_i\in J_f$ since they map to zero along $\mathcal O(X_f)\to \mathcal O(U)\otimes_\Bbbk \mathcal O(X_f)/I_f$. Then $\mathrm{d}\big(\frac{b_i}{a}\big)\in \mathrm{d}\big(\frac{J_f}{a}\big)$ and then the following determinant is in $\mathrm{Fit}_k(\widetilde\phi)(\widetilde X_{f,a})$ 
					\begin{equation}
						\det\begin{pmatrix}
							\xi_1.\frac{b_1}{a}&\cdots&\xi_1.\frac{b_{r-k}}{a}\\
							\vdots&\ddots&\vdots\\
							\xi_{r-k}.\frac{b_1}{a}&\cdots&\xi_{r-k}.\frac{b_{r-k}}{a}
						\end{pmatrix}=\det\begin{pmatrix}
							1&\cdots&0\\
							\vdots&\ddots&\vdots\\
							0&\cdots&1
						\end{pmatrix}=1. 
					\end{equation}
					This proves $\mathrm{Fit}_k(\widetilde X_{f,a})=\mathcal O(\widetilde X_{f,a})$. 
				\end{proof}

				\begin{corollary}\label{corollary: geometric U and U-hat quotients with WUU}
					For the linearisation $\hat U\curvearrowright(X,L)$, assume Condition \hyperref[condition: WUU]{WUU}. Let $k\in\mathbb N$ be the minimal dimension of $\mathrm{Stab}_U(x)$ for closed points $x\in X^0_{\min}$. Let $\phi:\Omega_{X^0_{\min}/\Bbbk}\to\mathfrak u^*\otimes_\Bbbk\mathcal O_{X^0_{\min}}$ be the infinitesimal action of $\mathfrak u$ on $X^0_{\min}$. Let $\mathrm{Fit}_k(\phi)\subseteq \mathcal O_{X^0_{\min}}$ denote the $k$th Fitting ideal of $\mathrm{coker}(\phi)$. Then: 
					\begin{itemize}
						\item the subset $X^0_{\min}\setminus \mathbb V(\mathrm{Fit}_k(\phi))$ has a quasi-projective universal geometric quotient by $U$; 
						\item the subset $X^0_{\min}\setminus \mathbb V(\mathrm{Fit}_k(\phi))\setminus UZ_{\min}$ is $\hat U$-invariant open in $X^0_{\min}$ and it has a quasi-projective universal geometric quotient by $\hat U$. 
					\end{itemize}
				\end{corollary}
				\begin{proof}
					By Theorem \ref{theorem: blowup from WUU to UU}, there is $\hat U$-equivariant blow-up $\pi:\widetilde X\to X$, and a linearisation $\hat U\curvearrowright\big(\widetilde X,\widetilde L\big)$ which satisfies Condition \hyperref[condition: UU]{UU}. Moreover, the centre of the blow-up is disjoint from $X^0_{\min}\setminus \mathbb V(\mathrm{Fit}_k(\phi))$. 

					By Theorem \ref{theorem: quotient for projective with UU}, we have: 
					\begin{itemize}
						\item $\widetilde X^0_{\min}$ has a quasi-projective universal geometric quotient by $U$; 
						\item $\widetilde X^0_{\min}\setminus U\widetilde Z_{\min}$ is open in $\widetilde X^0_{\min}$ and it has a projective universal geometric quotient by $\hat U$. 
					\end{itemize}

					We have 
					\begin{equation}
						\widetilde X^0_{\min}\cap \pi^{-1}\big(X^0_{\min}\setminus\mathbb V(\mathrm{Fit}_k(\phi))\big)=\pi^{-1}\big(X^0_{\min}\setminus\mathbb V(\mathrm{Fit}_k(\phi))\big)\cong X^0_{\min}\setminus\mathbb V(\mathrm{Fit}_k(\phi)). 
					\end{equation}
					Then $X^0_{\min}\setminus\mathbb V(\mathrm{Fit}_k(\phi))$ has a quasi-projective universal geometric quotient by $U$, since it is isomorphic to a $U$-invariant open subset of $\widetilde X^0_{\min}$. 

					Similarly, we have
					\begin{equation}
						\big(\widetilde X^0_{\min}\setminus U\widetilde Z_{\min}\big)\cap \pi^{-1}\big(X^0_{\min}\setminus \mathbb V(\mathrm{Fit}_k(\phi))\big)=\pi^{-1}\big(X^0_{\min}\setminus \mathbb V(\mathrm{Fit}_k(\phi))\setminus UZ_{\min}\big)\cong X^0_{\min}\setminus\mathbb V(\mathrm{Fit}_k(\phi))\setminus UZ_{\min}. 
					\end{equation}
					Then $X^0_{\min}\setminus\mathbb V(\mathrm{Fit}_k(\phi))\setminus UZ_{\min}$ is open and it has a quasi-projective universal geometric quotient by $\hat U$. 
				\end{proof}
			}
		}
	}

	\subsection{Stratification by unipotent stabilisers}
	{
		Let $U$ be an additive group of dimension $r$ over $\Bbbk$. Let $T\in\mathrm{Sch}/\Bbbk$. Let $U$ act on $T$. Let $\phi:\Omega_{T/\Bbbk}\to \mathfrak u^*\otimes_\Bbbk \mathcal O_T$ be the infinitesimal action of $\mathfrak u$ on $T$. 

		Consider the closed subschemes 
		\begin{equation}
			T=Z_{-1}\supset Z_0\supset\cdots\supset Z_r=\emptyset
		\end{equation}
		defined by the Fitting ideals of $\mathrm{coker}(\phi)$. Then by \href{https://stacks.math.columbia.edu/tag/05P8}{Lemma 05P8}, for each $\delta\in\mathbb N$ the locally closed subscheme $Z_{\delta-1}\setminus Z_\delta\subseteq T$ represents the functor $(\mathrm{Sch}/T)^{\mathrm{op}}\to \mathrm{Set}$
		\begin{equation}
			(T'\to T)\mapsto \begin{cases}
				\{\mathrm{pt}\}&,\quad \textrm{if }\mathcal (T'\to T)^*\mathrm{coker}(\phi)\textrm{ is locally free of rank }\delta\\
				\emptyset&,\quad \textrm{otherwise. }
			\end{cases}
		\end{equation}

		For $\delta\in\mathbb N$, let $S_\delta(T)\subseteq T$ denote the locally closed subscheme that is universal with respect to the property that the pullback of $\mathrm{coker}(\phi)$ is locally free of rank $\delta$. We call the following morphism the \emph{stratification by unipotent stabilisers}
		\begin{equation}
			\coprod_{\delta\in\mathbb N}S_\delta(T)\to T. 
		\end{equation}

		{
			\begin{lemma}\label{lemma: stratification by unipotent stabilisers commutes with equivariant immersions}
				Let $f:T_1\to T_2$ be an immersion. Let $U$ act on $T_1,T_2$ such that $f$ is $U$-equivariant. Let $\coprod_{\delta\in\mathbb N}S_\delta(T_i)\to T_i$ be the stratification by unipotent stabilisers for $i=1,2$. Then for all $\delta\in\mathbb N$, the following diagram is Cartesian 
				\begin{equation}
					\begin{tikzcd}
						S_\delta(T_1)\ar[r]\ar[d]\ar[rd,phantom,very near start,"\lrcorner"]&S_\delta(T_2)\ar[d]\\
						T_1\ar[r,"f"]&T_2. 
					\end{tikzcd}
				\end{equation}
			\end{lemma}
			\begin{proof}
				Let $\phi_i:\Omega_{T_i/\Bbbk}\to\mathfrak u^*\otimes_\Bbbk\mathcal O_{T_i}$ denote the infinitesimal action of $\mathfrak u:=\mathrm{Lie}(U)$ on $T_i$ for $i=1,2$. Recall that $S_\delta(T_i)\to T_i$ is universal with respect to the property that the pullback of $\mathrm{coker}(\phi_i)$ is locally free of rank $\delta$. By Lemma \ref{lemma: coker of infinitesimal stabiliser restricts to invariant subschemes}, we have $f^*\mathrm{coker}(\phi_2)\cong\mathrm{coker}(\phi_1)$. It is then clear that the diagram is Cartesian. 
			\end{proof}
		}

		Let $\hat U\curvearrowright(X,L)$. Assume Condition \ref{condition: non-emptiness of X^0_min}. We will apply the above stratification to the $U$-action on $X^0_{\min}$, and construct quotients of $S_\delta(X^0_{\min})$ for strata such that $S_\delta(X^0_{\min})\cap Z_{\min}\ne\emptyset$. 

		{
			\begin{theorem}\label{theorem: geometric U and U-hat quotients of strata}
				Let $\hat U\curvearrowright(X,L)$. Assume Condition \ref{condition: non-emptiness of X^0_min}. Let $\coprod_{\delta\in\mathbb N}S_\delta(X^0_{\min})\to X^0_{\min}$ be the stratification by unipotent stabilisers of $U\curvearrowright X^0_{\min}$. Then for each $\delta\in\mathbb N$, we have that:
				\begin{itemize}
					\item $S_\delta(X^0_{\min})\setminus UZ_{\min}\subseteq S_\delta(X^0_{\min})$ is an open subset in the stratum; 
					\item If $S_\delta(X^0_{\min})\cap Z_{\min}\ne\emptyset$, then $S_\delta(X^0_{\min})$ has a universal quasi-projective geometric quotient by $U$; 
					\item If $S_\delta(X^0_{\min})\cap Z_{\min}\ne\emptyset$, then $S_\delta(X^0_{\min})\setminus UZ_{\min}$ is open in $S_\delta(X^0_{\min})$ and it has a universal quasi-projective geometric quotient by $\hat U$. 
				\end{itemize}
			\end{theorem}
			\begin{proof}
				Fix $\delta\in\mathbb N$. Assume $S_\delta(X^0_{\min})\ne\emptyset$. Let $\phi:\Omega_{X^0_{\min}/\Bbbk}\to\mathfrak u^*\otimes_\Bbbk\mathcal O_{X^0_{\min}}$ be the infinitesimal action of $\mathfrak u$ on $X^0_{\min}$. Let $Z_k(\mathrm{coker}(\phi))\subseteq X^0_{\min}$ denote the closed subscheme associated to $\mathrm{Fit}_k(\phi)\subseteq \mathcal O_{X^0_{\min}}$ for every $k\in\mathbb N$. Then $S_\delta(X^0_{\min})=Z_{\delta-1}(\mathrm{coker}(\phi))\setminus Z_\delta(\mathrm{coker}(\phi))$. 

				When $S_\delta(X^0_{\min})\cap Z_{\min}=\emptyset$, we have that $S_\delta(X^0_{\min})\cap UZ_{\min}=\emptyset$, since $S_\delta(X^0_{\min})$ is $U$-invariant. Therefore $S_\delta(X^0_{\min})\setminus UZ_{\min}=S_\delta(X^0_{\min})$ is open in $S_\delta(X^0_{\min})$. 

				Assume $S_\delta(X^0_{\min})\cap Z_{\min}\ne\emptyset$. For simplicity write $Z_{\delta-1}:=Z_{\delta-1}(\mathrm{coker}(\phi))$. Let $\overline{Z_{\delta-1}}\hookrightarrow X$ be the scheme theoretic image of $Z_{\delta-1}\hookrightarrow X^0_{\min}\subseteq X$. Consider $\hat U\curvearrowright\big(\overline{Z_{\delta-1}},L\big|_{\overline{Z_{\delta-1}}}\big)$. It is easy to see the following equalities of subschemes: 
				\begin{equation}
					(\overline{Z_{\delta-1}})^0_{\min}=Z_{\delta-1},\quad \big((\overline{Z_{\delta-1}})^0_{\min}\big)^\lambda=Z_{\delta-1}\cap Z_{\min}. 
				\end{equation}
				Then for the open subset $S_\delta(X^0_{\min})$ of $(\overline{Z_{\delta-1}})^0_{\min}=Z_{\delta-1}$ we have: 
				\begin{itemize}
					\item $S_\delta(X^0_{\min})\cap \big((\overline{Z_{\delta-1}})^0_{\min}\big)^\lambda=S_\delta(X^0_{\min})\cap Z_{\min}\ne\emptyset$; 
					\item the cokernel of the infinitesimal action of $\mathfrak u$ on $S_\delta(X^0_{\min})$ is locally free, since this cokernel is isomorphic to $\mathrm{coker}(\phi)\big|_{S_\delta(X^0_{\min})}$ by Lemma \ref{lemma: coker of infinitesimal stabiliser restricts to invariant subschemes}. 
				\end{itemize}
				These two conditions above imply Condition \hyperref[condition: WUU]{WUU} for $\hat U\curvearrowright\big(\overline{Z_{\delta-1}},L\big|_{\overline{Z_{\delta-1}}}\big)$. We can apply Corollary \ref{corollary: geometric U and U-hat quotients with WUU}. The subsets in Corollary \ref{corollary: geometric U and U-hat quotients with WUU} are 
				\begin{equation}
					Z_\delta\Big\backslash \mathbb V\Big(\mathrm{Fit}_\delta\big(\mathrm{coker}(\phi)\big|_{Z_{\delta-1}}\big)\Big)=S_\delta(X^0_{\min})
				\end{equation}
				and 
				\begin{equation}
					Z_{\delta-1}\Big\backslash \mathbb V\Big(\mathrm{Fit}_\delta\big(\mathrm{coker}(\phi)|_{Z_{\delta-1}}\big)\Big)\Big\backslash U(Z_{\delta-1}\cap Z_{\min})=S_\delta(X^0_{\min})\setminus UZ_{\min}. 
				\end{equation}
				The statements then follow from Corollary \ref{corollary: geometric U and U-hat quotients with WUU}. 
			\end{proof}
		}

		We can prove a result similar to Theorem \ref{theorem: geometric quotients by parabolic with UU}. 
		{
			\begin{theorem}\label{theorem: geometric quotients of strata by parabolic}
				Let $U\rtimes(T\times K)\curvearrowright(X,L)$. Assume Condition \ref{condition: non-emptiness of X^0_min} for $\hat U\curvearrowright(X,L)$, where $\hat U:=U\rtimes T$. Let $Z_{\min}^{K,\mathrm{s}}\subseteq Z_{\min}$ denote the $K$-stable locus of $K\curvearrowright \big(Z_{\min},L\big|_{Z_{\min}}\big)$. Let $p_\lambda:X^0_{\min}\to Z_{\min}$ be the retraction. Let $\coprod_{\delta\in\mathbb N}S_\delta(X^0_{\min})\to X^0_{\min}$ be the stratification by unipotent stabilisers of $U\curvearrowright X^0_{\min}$. Let $\delta\in\mathbb N$. If $S_\delta(X^0_{\min})\cap Z_{\min}\ne\emptyset$, then the invariant open subset of $S_\delta(X^0_{\min})$ 
				\begin{equation}
					\big(S_\delta(X^0_{\min})\setminus UZ_{\min}\big)\cap p_\lambda^{-1}\big(Z_{\min}^{K,\mathrm{s}}\big)
				\end{equation}
				has a quasi-projective universal geometric quotient by $U\rtimes(T\times K)$. 
			\end{theorem}
			\begin{proof}
				The openness is from Theorem \ref{theorem: geometric U and U-hat quotients of strata}. In the proof of Theorem \ref{theorem: geometric U and U-hat quotients of strata}, we have seen that $\hat U\curvearrowright\big(\overline{Z_{\delta-1}},L\big|_{\overline{Z_{\delta-1}}}\big)$ satisfies Condition \hyperref[condition: WUU]{WUU}. By Theorem \ref{theorem: blowup from WUU to UU}, there is an equivariant blow-up of $\overline{Z_{\delta-1}}$ and a linearisation such that Condition \hyperref[condition: UU]{UU} is satisfied. 

				The subset $\big(S_\delta(X^0_{\min})\setminus UZ_{\min}\big)\cap p_\lambda^{-1}\big(Z_{\min}^{K,\mathrm{s}}\big)$ is open in $\overline{Z_{\delta-1}}$ and it is disjoint from the blow-up centre. So we can think it as an open subset of the blow-up. By Theorem \ref{theorem: geometric quotients by parabolic with UU}, there is an open subset of the blow-up which has a quasi-projective universal geometric quotient by $U\rtimes(T\times K)$. It is easy to see this open subset contains $\big(S_\delta(X^0_{\min})\setminus UZ_{\min}\big)\cap p_\lambda^{-1}\big(Z_{\min}^{K,\mathrm{s}}\big)$. So $\big(S_\delta(X^0_{\min})\setminus UZ_{\min}\big)\cap p_\lambda^{-1}\big(Z_{\min}^{K,\mathrm{s}}\big)$ has a quasi-projective universal geometric quotient by $U\rtimes (T\times K)$. 
			\end{proof}
		}
	}
}

\section{Sheaves of length 2}
{
	Let $\tau$ be a Harder-Narasimhan type of length 2. Let $P(t),P_i(t)\in\mathbb Q[t]$ for $1\leq i\leq 2$ be the Hilbert polynomials associated to $\tau$, in the sense that if $\mathcal E$ is a sheaf of type $\tau$, then its Harder-Narasimhan filtration
	\begin{equation}
		0=\mathcal E_{\leq 0}\subsetneq\mathcal E_{\leq 1}\subsetneq \mathcal E_{\leq 2}=\mathcal E
	\end{equation}
	has length 2, $P(t)$ is the Hilbert polynomial of $\mathcal E$ and $P_i(t)$ is the Hilbert polynomial of $\mathcal E_i:=\mathcal E_{\leq i}/\mathcal E_{\leq i-1}$ for $i=1,2$. Recall that $B$ is a connected projective scheme over $\Bbbk$ and $\mathcal L$ is an ample line bundle on $B$. We will assume that there exists a $(\tau,\delta)$-stable sheaf for some $\delta\in\mathbb N$. 
	
	Recall $\mathrm P_\tau\curvearrowright\big(\overline{Y_\tau^{\mathrm{s}}},\mathcal M\big|_{\overline{Y_\tau^{\mathrm{s}}}}\big)$ described at the end of Section \ref{subsubsection of moduli of unstable sheaves}. Let $\mathrm U_\tau\subseteq \mathrm P_\tau$ be the unipotent radical and let $\mathrm K_\tau\subseteq \mathrm P_\tau$ be its semisimple part $\mathrm K_\tau=\prod_{i=1}^2\mathrm{SL}(V_i)$. Let $\mathrm T_\tau=\mathbb G_m$ and let $\mathrm T_\tau\to \mathrm P_\tau$ be the following homomorphism 
	\begin{equation}
		t\mapsto \begin{pmatrix}
			t^{\frac{P_2}{(P_1,P_2)}}1_{V_1}\\
			&t^{\frac{-P_1}{(P_1,P_2)}}1_{V_2}
		\end{pmatrix}
	\end{equation}
	where $P_i:=P_i(m)=\dim V_i$ and $(P_1,P_2):=\mathrm{gcd}(P_1,P_2)$. Then $\mathrm U_\tau\rtimes(\mathrm T_\tau\times \mathrm K_\tau)\to\mathrm P_\tau$ is surjective with finite kernel. 
	
	Consider the linearisation
	\begin{equation}
		\mathrm U_\tau\rtimes(\mathrm T_\tau\times \mathrm K_\tau)\curvearrowright\big(\overline{Y_\tau^{\mathrm{s}}},\mathcal M\big|_{\overline{Y_\tau^{\mathrm{s}}}}\big). 
	\end{equation}
	By Lemma \ref{lemma: stratification by unipotent stabilisers commutes with equivariant immersions} we have 
	\begin{itemize}
		\item[(1)] $S_\delta\big((\overline{Y_\tau^{\mathrm{s}}})^0_{\min}\big)\cap Y_\tau^{\mathrm{s}}=S_\delta(Y_\tau^{\mathrm{s}})$ for all $\delta\in\mathbb N$, where $\coprod_{\delta\in\mathbb N}S_\delta(Y_\tau^{\mathrm{s}})\to Y_\tau^{\mathrm{s}}$ is the stratification by unipotent stabilisers. 
	\end{itemize}
	By Proposition \ref{proposition: Y_tau^s closure and its Z_min^K-s} we have 
	\begin{itemize}
		\item[(2)] $p_\tau^{-1}\big(Z_{\min}^{\mathrm K_\tau,\mathrm{s}}(\overline{Y_\tau^{\mathrm{s}}})\big)=Y_\tau^{\mathrm{s}}$. 
	\end{itemize}

	If $S_\delta\big((\overline{Y_\tau^{\mathrm{s}}})^0_{\min}\big)\cap Z_{\min}(\overline{Y_\tau^{\mathrm{s}}})\ne\emptyset$, then by Theorem \ref{theorem: geometric U and U-hat quotients of strata} the following subscheme is invariant open in $S_\delta\big((\overline{Y_\tau^{\mathrm{s}}})^0_{\min}\big)$ and it has a quasi-projective geometric quotient by $\mathrm U_\tau\rtimes(\mathrm T_\tau\times \mathrm K_\tau)$ 
	\begin{equation}
		\begin{split}
			&\big(S_\delta\big((\overline{Y_\tau^{\mathrm{s}}})^0_{\min}\big)\big\backslash \mathrm U_\tau Z_{\min}(\overline{Y_\tau^{\mathrm{s}}})\big)\cap p_\tau^{-1}\big(Z_{\min}^{\mathrm K_\tau,\mathrm{s}}(\overline{Y_\tau^{\mathrm{s}}})\big)\\
			=&\Big(S_\delta\big((\overline{Y_\tau^{\mathrm{s}}})^0_{\min}\big)\cap p_\tau^{-1}\big(Z_{\min}^{\mathrm K_\tau,\mathrm{s}}(\overline{Y_\tau^{\mathrm{s}}})\big)\Big)\Big\backslash \mathrm U_\tau Z_{\min}(\overline{Y_\tau^{\mathrm{s}}})\\
			\overset{(2)}{=\!=}&\Big(S_\delta\big((\overline{Y_\tau^{\mathrm{s}}})^0_{\min}\big)\cap Y_\tau^{\mathrm{s}}\Big)\Big\backslash \mathrm U_\tau Z_{\min}(\overline{Y_\tau^{\mathrm{s}}})\\
			\overset{(1)}{=\!=}&S_\delta(Y_\tau^{\mathrm{s}})\big\backslash \mathrm U_\tau Z_{\min}(\overline{Y_\tau^{\mathrm{s}}}). 
		\end{split}
	\end{equation}

	Denote 
	\begin{equation}\label{equation: definition of Y_(tau,delta)^s}
		Y_{\tau,\delta}^{\mathrm{s}}:=S_\delta(Y_\tau^{\mathrm{s}})\setminus \mathrm U_\tau Z_{\min}(\overline{Y_\tau^{\mathrm{s}}}). 
	\end{equation}
	Then we have a quasi-projective geometric quotient by $\mathrm P_\tau$ 
	\begin{equation}
		Y_{\tau,\delta}^{\mathrm{s}}\to Y_{\tau,\delta}^{\mathrm{s}}/\mathrm P_\tau. 
	\end{equation}
	The condition $S_\delta\big((\overline{Y_\tau^{\mathrm{s}}})^0_{\min}\big)\cap Z_{\min}(\overline{Y_\tau^{\mathrm{s}}})\ne\emptyset$ is satisfied if $Y_{\tau,\delta}^{\mathrm{s}}\ne\emptyset$ by Lemma \ref{lemma: dim Stab_U constant along p_tau for length 2 sheaves}. The condition $Y_{\tau,\delta}^{\mathrm{s}}\ne\emptyset$ is equivalent to the existence of a $(\tau,\delta)$-stable sheaf on $B$. We will always assume $(\tau,\delta)$ is such that there exists a $(\tau,\delta)$-stable sheaf. 

	Similar to $Y_\tau$ and $Y_\tau^{\mathrm{s}}$, the schemes $S_\delta(Y_\tau^{\mathrm{s}})$ and $Y_{\tau,\delta}^{\mathrm{s}}$ also have moduli-theoretic descriptions. We will describe their functors of points (Lemma \ref{lemma: moduli descriptions of S_alpha(Y_tau^s) and Y_(tau,alpha)^s}) in Section \ref{subsection of unipotent endomorphisms}. The relation between $\mathrm P_\tau\curvearrowright Y_{\tau,\delta}^{\mathrm{s}}$ and the moduli functor $\mathbf M'_{\tau,\delta}$ (Lemma \ref{lemma: theta from Y_(tau,delta)^s/_pP_tau to M'_(tau,delta) sheafifies to an isomorphism}) will be discussed in Section \ref{subsection of the moduli functor M'_(tau,alpha) in length 2}.

	\subsection{Unipotent endomorphisms}\label{subsection of unipotent endomorphisms}
	{
		Let $\mathrm U_\tau\subseteq \mathrm P_\tau$ be the unipotent radical and let $\mathfrak u$ denote its Lie algebra. Then $\mathfrak u:=\mathrm{Hom}_\Bbbk(V_2,V_1)$. 
		
		Let $T\in\mathrm{Sch}/\Bbbk$. Recall that a morphism $T\to Y^{\mathrm{s}}_\tau$ is equivalent to an isomorphism class $[(\mathcal E,q)]$ such that 
		\begin{equation}
			q:V(-m)_T\to \mathcal E
		\end{equation}
		is surjective, each $\mathcal E_i$ if flat over $T$ with Hilbert polynomial $P_i(t)$, and fibres of each $\mathcal E_i$ at closed points of $T$ are stable, where $\mathcal E_{\leq i}:=q\big(V_{\leq i}(-m)_T\big)$ and $\mathcal E_i:=\mathcal E_{\leq i}/\mathcal E_{\leq i-1}$. 

		{
			Recall Definition \ref{definition: T-valued infinitesimal transformation in G} and Definition \ref{definition: T-infinitesimal transformations stabilising T-points} for scheme-valued infinitesimal transformations and when they stabilise scheme-valued points. 
			\begin{lemma}\label{lemma: natural isomorphisms of infinitesimal stabilisers with unipotent endomorphisms}
				When $m\gg0$, we have the following commutative diagram of $\Gamma(T,\mathcal O_T)$-modules for any $T\in\mathrm{Sch}/\Bbbk$ and any $T$-valued point $[(\mathcal E,q)]$ in $Y_\tau^{\mathrm{s}}$ 
				\begin{equation}\label{equation: diagram in lemma that infinitesimal stabilisers isomorphic to suitable unipotent endomorphisms}
					\begin{tikzcd}
						\Big\{\begin{matrix}T\textrm{-valued infinitesimal}\\\textrm{transformations in }\mathrm U_\tau\end{matrix}\Big\}\ar[r,"\cong"{sloped}]&\mathrm{Hom}_{\mathcal O_T}\big(V_2\otimes_\Bbbk\mathcal O_T,V_1\otimes_\Bbbk\mathcal O_T\big)\\
						\left\{\begin{matrix}\rho:\rho\textrm{ is a }T\textrm{-valued}\\\textrm{infinitesimal}\\\textrm{transformation in }\mathrm U_\tau\\\textrm{and }\rho\textrm{ stabilises }[(\mathcal E,q)]\end{matrix}\right\}\ar[r,"\cong"{sloped}]\ar[u,"\subseteq"{sloped}]&\mathrm{Hom}_{\mathcal O_{B_T}}(\mathcal E_2,\mathcal E_1)\ar[u,"\subseteq"{sloped}]
					\end{tikzcd}
				\end{equation}
				where the horizontal maps are isomorphisms and the vertical maps are inclusions. Moreover, the top map is natural in $T\in\mathrm{Sch}/\Bbbk$ and the bottom map is natural in $T\in\mathrm{Sch}/Y_\tau^{\mathrm{s}}$. 
			\end{lemma}
			\begin{proof}
				Let $f:B\to \mathrm{Spec}(\Bbbk)$ denote the structure morphism. For $T\in\mathrm{Sch}/\Bbbk$, let $f_T:B_T\to T$ denote the base change of $f$. For a coherent sheaf $\mathcal F$ on $B_T$, let $\mathcal F(m)$ denote $\mathcal F\otimes_{\mathcal O_{B_T}}\mathcal L_T^m$. For the right vertical arrow, let $\varphi:\mathcal E\to \mathcal E$ be $\mathcal O_{B_T}$-linear such that $\varphi(\mathcal E_{\leq i})\subseteq\mathcal E_{\leq i-1}$ for $i=1,2$. 

				There is an exact sequence 
				\begin{equation}
					\begin{tikzcd}
						0\ar[r]&\mathrm{Hom}_{\mathcal O_{B_T}}(\mathcal E_2,\mathcal E_1)\ar[r]&\mathrm{Hom}_{\mathcal O_{B_T}}(\mathcal E_{\leq 2},\mathcal E_1)\ar[r]&\mathrm{Hom}_{\mathcal O_{B_T}}(\mathcal E_{\leq 1},\mathcal E_1). 
					\end{tikzcd}
				\end{equation}
				We can think $\varphi\in\mathrm{Hom}_{\mathcal O_{B_T}}(\mathcal E_2,\mathcal E_1)$ as a morphism $\varphi:\mathcal E_{\leq 2}\to \mathcal E_1$ such that $\varphi(\mathcal E_{\leq 1})=0$. Consider the diagram on $B_T$ 
				\begin{equation}
					\begin{tikzcd}
						V_{\leq 2}(-m)_T\ar[d,"q_{\leq 2}"]&V_1(-m)_T\ar[d,"q_1"]\\
						\mathcal E_{\leq 2}\ar[r,"\varphi"]&\mathcal E_1. 
					\end{tikzcd}
				\end{equation}

				Apply the functor $\mathcal F\mapsto \mathcal F(m)\mapsto f_{T,*}(\mathcal F(m))$ to the diagram above to obtain the following diagram 
				\begin{equation}\label{equation: diagram of defining eta from varphi}
					\begin{tikzcd}
						V_{\leq 2}\otimes_\Bbbk\mathcal O_T\ar[d,"f_{T,*}(q_{\leq 2}(m))"]\ar[rr,dashed,"\eta"]&&V_1\otimes_\Bbbk\mathcal O_T\ar[d,"f_{T,*}(q_1(m))"]\\
						f_{T,*}(\mathcal E_{\leq 2}(m))\ar[rr,"f_{T,*}(\varphi(m))"]&&f_{T,*}(\mathcal E_1(m)). 
					\end{tikzcd}
				\end{equation}
				When $m\gg0$ is large, we have 
				\begin{itemize}
					\item $f_{T,*}(\mathcal E_j(m))$ is locally free of rank $P_j(m)$ for $j=1,2$; 
					\item $f_{T,*}(\mathcal E_{\leq j}(m))$ is locally free of rank $P_{\leq j}(m)$ for $j=1,2$
				\end{itemize}
				for all $T\in\mathrm{Sch}/\Bbbk$ and $[(\mathcal E,q)]$. So $f_{T,*}(q_1(m))$ is a morphism between locally free sheaves of the same rank $P_1(m)$. For any closed point $\{t\}\hookrightarrow T$, the morphism $f_{T,*}(q_1(m))$ pulls back to an isomorphism $V_1\cong H^0(B,\mathcal E_1(m)|_{B_t})$. By \href{https://stacks.math.columbia.edu/tag/00O0}{Lemma 00O0}, we have that $f_{T,*}(q_1(m))$ is an isomorphism since its pullbacks at closed points are isomorphisms. Similarly $f_{T,*}(q_{\leq 2}(m))$ is an isomorphism. Therefore the morphism $\eta:V_{\leq 2}\otimes_\Bbbk\mathcal O_T\to V_1\otimes_\Bbbk\mathcal O_T$ exists in diagram \eqref{equation: diagram of defining eta from varphi}. Since $\varphi(\mathcal E_{\leq 1})=0$, we have $\eta(V_{\leq 1}\otimes_\Bbbk\mathcal O_T)=0$, i.e. $\eta$ defines a morphism $\bar\eta:V_2\otimes_\Bbbk\mathcal O_T\to V_1\otimes_\Bbbk\mathcal O_T$. The map $\varphi\mapsto \bar\eta$ defines the right vertical map in \eqref{equation: diagram in lemma that infinitesimal stabilisers isomorphic to suitable unipotent endomorphisms}. 
				
				The right vertical map in diagram \eqref{equation: diagram in lemma that infinitesimal stabilisers isomorphic to suitable unipotent endomorphisms} is injective, since if $\varphi\mapsto \bar\eta$, then $\varphi,\bar\eta$ satisfy the following diagram 
				\begin{equation}
					\begin{tikzcd}
						V(-m)_T\ar[r,"\eta(-m)"]\ar[d,->>,"q"]&V(-m)_T\ar[d,->>,"q"]\\\mathcal E\ar[r,"\varphi"]&\mathcal E. 
					\end{tikzcd}
				\end{equation}
				where $\eta$ is the composition $V_{\leq 2}\otimes_\Bbbk\mathcal O_T\to V_2\otimes_\Bbbk\mathcal O_T\xrightarrow{\bar\eta}V_1\otimes_\Bbbk\mathcal O_T$ and $\eta(-m):=\big((f_T)^*\eta\big)\otimes_{\mathcal O_{B_T}}1_{\mathcal L_T^{-m}}$. 
				
				For the top map in diagram \eqref{equation: diagram in lemma that infinitesimal stabilisers isomorphic to suitable unipotent endomorphisms}, let $\rho:T[\epsilon]\to \mathrm U_\tau$ be a $T$-valued infinitesimal transformation in $\mathrm U_\tau$. A morphism $\rho:T[\epsilon]\to\mathrm U_\tau$ is equivalent to a morphism of coherent sheaves on $T[\epsilon]$ 
				\begin{equation}
					u:V\otimes_\Bbbk\mathcal O_{T[\epsilon]}\to V\otimes_\Bbbk\mathcal O_{T[\epsilon]}
				\end{equation}
				such that: 
				\begin{itemize}
					\item[(1)] $u(V_{\leq i}\otimes_\Bbbk\mathcal O_{T[\epsilon]})\subseteq V_{\leq i}\otimes_\Bbbk\mathcal O_{T[\epsilon]}$, and we denote $u_{\leq i}:V_{\leq i}\otimes_\Bbbk\mathcal O_{T[\epsilon]}\to V_{\leq i}\otimes_\Bbbk\mathcal O_{T[\epsilon]}$; 
					\item[(2)] the induced morphism $u_i:V_i\otimes_\Bbbk\mathcal O_{T[\epsilon]}\to V_i\otimes_\Bbbk\mathcal O_{T[\epsilon]}$ is the identity 
					\begin{equation}
						\begin{tikzcd}
							0\ar[r]&V_{\leq i-1}\otimes_\Bbbk\mathcal O_{T[\epsilon]}\ar[r,"\subseteq"]\ar[d,"u_{\leq i-1}"]&V_{\leq i}\otimes_\Bbbk\mathcal O_{T[\epsilon]}\ar[r]\ar[d,"u_{\leq i}"]&V_i\otimes_\Bbbk\mathcal O_{T[\epsilon]}\ar[r]\ar[d,"u_i=1"]&0\\
							0\ar[r]&V_{\leq i-1}\otimes_\Bbbk\mathcal O_{T[\epsilon]}\ar[r,"\subseteq"]&V_{\leq i}\otimes_\Bbbk\mathcal O_{T[\epsilon]}\ar[r]&V_i\otimes_\Bbbk\mathcal O_{T[\epsilon]}\ar[r]&0. 
						\end{tikzcd}
					\end{equation}
				\end{itemize}
				Note that the underlying topological spaces of $T$ and $T[\epsilon]$ are the same, and $\mathcal O_{T[\epsilon]}\cong\mathcal O_T\oplus\epsilon\mathcal O_T$ as $\mathcal O_T$-modules. The morphism $u:V\otimes_\Bbbk\mathcal O_{T[\epsilon]}\to V\otimes_\Bbbk\mathcal O_{T[\epsilon]}$ is equivalent to two $\mathcal O_T$-linear morphisms 
				\begin{equation}
					\mu:V\otimes_\Bbbk\mathcal O_T\to V\otimes_\Bbbk\mathcal O_T,\quad \eta:V\otimes_\Bbbk\mathcal O_T\to V\otimes_\Bbbk\mathcal O_T
				\end{equation}
				in the sense that according to the decomposition $V\otimes_\Bbbk\mathcal O_{T[\epsilon]}\cong (V\otimes_\Bbbk\mathcal O_T)\oplus\epsilon (V\otimes_\Bbbk\mathcal O_T)$, we have 
				\begin{equation}
					u=\begin{pmatrix}\mu&\\\epsilon\eta&\mu_\epsilon\end{pmatrix}
				\end{equation}
				where $\mu_\epsilon:\epsilon V\otimes_\Bbbk\mathcal O_T\to \epsilon V\otimes_\Bbbk\mathcal O_T$ is the morphism equivalent to $\mu$. 
				
				For a morphism $\rho:T[\epsilon]\to\mathrm U_\tau$ to be an infinitesimal transformation, the diagram $\begin{tikzcd}T\ar[r,hook]\ar[d]&T[\epsilon]\ar[d,"\rho"]\\\mathrm{Spec}(\Bbbk)\ar[r,"e",hook]&\mathrm U_\tau\end{tikzcd}$ has to commute, which is equivalent to the following condition for $u$: 
				\begin{itemize}
					\item[(3)] $\mu=1_{V\otimes_\Bbbk\mathcal O_T}$, i.e. $u=\begin{pmatrix}1&\\\epsilon\eta&1_\epsilon\end{pmatrix}$. 
				\end{itemize}

				Condition $(3)$ for $u$ says the map $\rho\mapsto \eta$ is injective when $\rho$ is a $T$-valued infinitesimal transformation in $\mathrm U_\tau$. With $(3)$, condition $(1)$ and $(2)$ together for $u$ are equivalent to the following 
				\begin{equation}
					\eta(V_{\leq i}\otimes_\Bbbk\mathcal O_T)\subseteq V_{\leq i-1}\otimes \mathcal O_T,\quad i=1,2. 
				\end{equation}
				Such $\eta$ is equivalent to $\bar\eta:V_2\otimes_\Bbbk\mathcal O_T\to V_1\otimes_\Bbbk\mathcal O_T$. Therefore the association $\rho\mapsto \bar\eta$ defines the top map in diagram \eqref{equation: diagram in lemma that infinitesimal stabilisers isomorphic to suitable unipotent endomorphisms}, which is bijective from its construction. 
				
				For the bottom map in diagram \eqref{equation: diagram in lemma that infinitesimal stabilisers isomorphic to suitable unipotent endomorphisms}, let $\rho:T[\epsilon]\to\mathrm U_\tau$ be a $T$-valued infinitesimal transformation in $\mathrm U_\tau$ and assume $\rho$ is mapped to $\bar\eta:V_2\otimes_\Bbbk\mathcal O_T\to V_1\otimes_\Bbbk\mathcal O_T$ along the top map in diagram \eqref{equation: diagram in lemma that infinitesimal stabilisers isomorphic to suitable unipotent endomorphisms}. Let $\eta:V\otimes_\Bbbk\mathcal O_T\to V\otimes_\Bbbk\mathcal O_T$ be the composition $V\otimes_\Bbbk\mathcal O_T\twoheadrightarrow V_2\otimes_\Bbbk\mathcal O_T\xrightarrow{\bar\eta}V_1\otimes_\Bbbk\mathcal O_T\subseteq V\otimes_\Bbbk\mathcal O_T$. We have $\eta(V_{\leq i}\otimes_\Bbbk\mathcal O_T)\subseteq V_{\leq i-1}\otimes_\Bbbk\mathcal O_T$ for $i=1,2$. 
				
				Let $h:T\to Y_\tau^{\mathrm{s}}$ be the morphism corresponding to $[(\mathcal E,q)]$. Then $\rho:T[\epsilon]\to \mathrm U_\tau$ stabilises $[(\mathcal E,q)]$ if and only if the following diagram commutes 
				\begin{equation}
					\begin{tikzcd}
						T[\epsilon]\ar[r,"{\mathrm{pr}_T}"]\ar[d,"{(\rho,h\circ \mathrm{pr}_T)}"]&T\ar[d,"h"]\\
						\mathrm U_\tau\times_\Bbbk Y_\tau^{\mathrm{s}}\ar[r,"\sigma"]&Y_\tau^{\mathrm{s}}
					\end{tikzcd}
				\end{equation}
				where $\sigma:\mathrm U_\tau\times_\Bbbk Y_\tau^{\mathrm{s}}\to Y_\tau^{\mathrm{s}}$ is the action morphism of $\mathrm U_\tau\curvearrowright Y_\tau^{\mathrm{s}}$. There are two families parametrised by $T[\epsilon]$ of quotients corresponding to two morphisms $h\circ \mathrm{pr}_T$ and $\sigma\circ(\rho,h\circ\mathrm{pr}_T)$. The diagram above commutes if and only if $h\circ \mathrm{pr}_T=\sigma\circ(\rho,h\circ\mathrm{pr}_T)$, if and only if the two families over $T[\epsilon]$ are equivalent. 
				
				The family over $T[\epsilon]$ corresponding to $h\circ \mathrm{pr}_T$ is the following 
				\begin{equation}
					\begin{tikzcd}[ampersand replacement =\&]
						V(-m)_T\oplus \epsilon V(-m)_T\ar[r,"{\begin{pmatrix}q&\\ &q_\epsilon\end{pmatrix}}"]\&\mathcal E\oplus\epsilon\mathcal E\ar[r]\&0
					\end{tikzcd}
				\end{equation}
				where $q_\epsilon:\epsilon V(-m)_T\to \epsilon\mathcal E$ is equivalent to $q$. The family corresponding to $\sigma\circ (\rho,h\circ\mathrm{pr}_T)$ is the composition 
				\begin{equation}
					\begin{tikzcd}[ampersand replacement =\&]
						V(-m)_T\oplus \epsilon V(-m)_T\ar[rr,"{\begin{pmatrix}1&\\-\epsilon\eta(-m)&1_\epsilon\end{pmatrix}}"]\&\&V(-m)_T\oplus\epsilon V(-m)_T\ar[r,"{\begin{pmatrix}q&\\&q_\epsilon\end{pmatrix}}"]\&\mathcal E\oplus\epsilon\mathcal E
					\end{tikzcd}
				\end{equation}
				where $\eta(-m):V(-m)_T\to V(-m)_T$ denotes $\eta(-m):=\big((f_T)^*\eta\big)\otimes_{\mathcal O_{B_T}}1_{\mathcal L_T^{-m}}$. It is easy to see the two families are equivalent if and only if $\eta(-m)(\ker q)\subseteq \ker q$, which is equivalent to the existence of $\varphi:\mathcal E\to\mathcal E$ such that the following diagram commutes 
				\begin{equation}
					\begin{tikzcd}
						V(-m)_T\ar[r,"\eta(-m)"]\ar[d,"q",->>]&V(-m)_T\ar[d,"q",->>]\\
						\mathcal E\ar[r,"\varphi"]&\mathcal E
					\end{tikzcd}
				\end{equation}
				which is equivalent to that $\bar\eta$ is in the image of the right vertical map in diagram \eqref{equation: diagram in lemma that infinitesimal stabilisers isomorphic to suitable unipotent endomorphisms}, i.e. the association $\rho\mapsto \varphi$ defines the bottom map in diagram \eqref{equation: diagram in lemma that infinitesimal stabilisers isomorphic to suitable unipotent endomorphisms}. 
				
				We have now defined all the maps in diagram \eqref{equation: diagram in lemma that infinitesimal stabilisers isomorphic to suitable unipotent endomorphisms}, with the horizontal maps being bijective and the vertical maps being injective. It is easy to check that they are all $\Gamma(T,\mathcal O_T)$-linear and natural in $T\in\mathrm{Sch}/\Bbbk$ and $T\in\mathrm{Sch}/Y_\tau^{\mathrm{s}}$ respectively. 
			\end{proof}
		}
		
		{
			Recall that $Y_{\tau,\delta}^{\mathrm{s}}$ is defined in \eqref{equation: definition of Y_(tau,delta)^s}, and $f:B\to\mathrm{Spec}(\Bbbk)$ is the structure morphism. We describe the functors of points of $Y_{\tau,\delta}^{\mathrm{s}}\subseteq S_\delta(Y_\tau^{\mathrm{s}})$. 

			\begin{lemma}\label{lemma: moduli descriptions of S_alpha(Y_tau^s) and Y_(tau,alpha)^s}
				The scheme $S_\delta(Y_\tau^{\mathrm{s}})$ represents the functor which sends $T\in\mathrm{Sch}/\Bbbk$ to the set of isomorphism classes $[(\mathcal E,q)]$ such that: 
				\begin{itemize}
					\item[(1)] $q:V(-m)_T\to \mathcal E$ is a surjective morphism of coherent sheaves on $B_T$, and $\mathcal E_i$ is flat over $T$ with Hilbert polynomial $P_i(t)$ for $i=1,2$; 
					\item[(2)] fibres of $\mathcal E_i$ at closed points are stable for $i=1,2$; 
					\item[(3)] for any closed subscheme $T'\hookrightarrow T$, the coherent sheaf 
					\begin{equation}
						f_{T',*}\mathcal Hom_{\mathcal O_{B_{T'}}}\big(\mathcal E_2|_{B_{T'}},\mathcal E_1|_{B_{T'}}\big)
					\end{equation}
					is locally free of rank $\delta$ on $T'$. 
				\end{itemize}
				
				The open subscheme $Y_{\tau,\delta}^{\mathrm{s}}\subseteq S_\delta(Y_\tau^{\mathrm{s}})$ represents the sub-functor which sends $T\in\mathrm{Sch}/\Bbbk$ to the set of isomorphism classes $[(\mathcal E,q)]$ such that $(1),(2),(3)$ above hold and: 
				\begin{itemize}
					\item[(4)] fibres of $\mathcal E$ at closed points $t\in T$ are non-split, in the sense that $\mathcal E|_{B_t}\not\cong \mathcal E_1|_{B_t}\oplus\mathcal E_2|_{B_t}$. 
				\end{itemize}
			\end{lemma}
			\begin{proof}
				Recall in Section \ref{subsubsection of moduli of unstable sheaves}, the schemes $Y_\tau^{\mathrm{s}}\subseteq Y_\tau$ are defined by functors of points. Let $T\in\mathrm{Sch}/\Bbbk$. Then $(1)$ is equivalent to that $[(\mathcal E,q)]$ corresponds to a $T$-valued point in $Y_\tau$, and $(2)$ is equivalent to that $[(\mathcal E,q)]$ corresponds to a $T$-valued point in $Y_\tau^{\mathrm{s}}$. 
				
				By definition of the stratification by unipotent stabilisers, the morphism $T\to Y_\tau^{\mathrm{s}}$ factors through $S_\delta(Y_\tau^{\mathrm{s}})$ if and only if $\mathrm{coker}(\phi)|_T$ is locally free of rank $\delta$, where $\phi:\Omega_{Y_\tau^{\mathrm{s}}/\Bbbk}\to \mathfrak u^*\otimes_\Bbbk\mathcal O_{Y_\tau^{\mathrm{s}}}$ is the infinitesimal action. By Lemma \ref{lemma: Stab_g(t') criterion for local freedom of coker(phi)|_T} and Lemma \ref{lemma: bijections between infinitesimal transformations stabilising a point with Stab_g(t)} this is further equivalent to 
				\begin{itemize}
					\item[(3')] for any closed subscheme $T'\hookrightarrow T$, the following defines a locally free sheaf of rank $\delta$ on $T'$ 
					\begin{equation}
						\begin{split}
							W'&\mapsto\Big\{\begin{matrix}\rho:\rho\textrm{ is a }W'\textrm{-valued infinitesimal transformation in }\mathrm U_\tau\\\textrm{and }\rho\textrm{ stabilises }W'\subseteq T'\hookrightarrow T\to Y_\tau^{\mathrm{s}}\end{matrix}\Big\}\\
							&\textrm{for open subsets }W'\subseteq T'. 
						\end{split}
					\end{equation}
				\end{itemize}
				
				By Lemma \ref{lemma: natural isomorphisms of infinitesimal stabilisers with unipotent endomorphisms}, the sheaf in $(3')$ is isomorphic to the following sheaf on $T'$ 
				\begin{equation}
					\begin{split}
						W'\mapsto& \mathrm{Hom}_{\mathcal O_{B_{W'}}}\big(\mathcal E_2|_{B_{W'}},\mathcal E_1|_{B_{W'}}\big)\\
						\cong&\mathcal Hom_{\mathcal O_{B_{T'}}}\big(\mathcal E_2|_{B_{T'}},\mathcal E_1|_{B_{T'}}\big)\big((f_{T'})^{-1}(W')\big)\\
						=&\Big(f_{T',*}\mathcal Hom_{\mathcal O_{B_{T'}}}\big(\mathcal E_2|_{B_{T'}},\mathcal E_1|_{B_{T'}}\big)\Big)(W'). 
					\end{split}
				\end{equation}
				So $(3)$ and $(3')$ are equivalent. This proves the moduli description of $S_\delta(Y_\tau^{\mathrm{s}})$. 
				
				We have that $Y_{\tau,\delta}^{\mathrm{s}}\subseteq S_\delta(Y_\tau^{\mathrm{s}})$ is an open subscheme, since $Y_{\tau,\delta}^{\mathrm{s}}$ is open in a larger scheme $S_\delta\big((\overline{Y_\tau^{\mathrm{s}}})^0_{\min}\big)$ by Theorem \ref{theorem: geometric U and U-hat quotients of strata}. The morphism $T\to S_\delta(Y_\tau^{\mathrm{s}})$ factors through $Y_{\tau,\delta}^{\mathrm{s}}$ if and only if any closed point $t\in T$ is mapped into $Y_{\tau,\delta}^{\mathrm{s}}$, which is equivalent to: 
				\begin{itemize}
					\item[(4')] fibres of $\mathcal E$ at closed points $t\in T$, denoted by $\mathcal E|_{B_t}$, satisfy the condition that the following two closed points in $Y_\tau^{\mathrm{s}}$ are not in one $\mathrm U_\tau$-orbit
					\begin{equation}
						\big[\big(\mathcal E|_{B_t},q|_{B_t}\big)\big],\quad \Big[\Big(\bigoplus_{i=1}^2\mathcal E_i|_{B_t},\; \bigoplus_{i=1}^2q_i|_{B_t}\Big)\Big]\in Y_\tau^{\mathrm{s}}. 
					\end{equation}
				\end{itemize}
				If $(4)\iff(4')$, then the moduli description of $Y_{\tau,\delta}^{\mathrm{s}}$ is proved, hence the lemma. 
				
				Obviously we have $(4)\implies(4')$. Conversely, it suffices to prove that if $[(\mathcal E,q)]\in Y_\tau^{\mathrm{s}}$ is a closed point such that $\mathcal E\cong \bigoplus_{i=1}^2\mathcal E_i$, then $\Big[\Big(\bigoplus_{i=1}^2\mathcal E_i,\bigoplus_{i=1}^2q_i\Big)\Big]$ is in the $\mathrm U_\tau$-orbit of $[(\mathcal E,q)]$. Let $[(\mathcal E,q)]\in Y_\tau^{\mathrm{s}}$ be such a closed point. Then the following diagram commutes 
				\begin{equation}
					\begin{tikzcd}
						0\ar[r]&V_1(-m)\ar[d,->>,"q_1"]\ar[r]&V(-m)\ar[d,->>,"q"]\ar[r]&V_2(-m)\ar[d,->>,"q_2"]\ar[r]&0\\
						0\ar[r]&\mathcal E_1\ar[r]\ar[d,equal]&\mathcal E\ar[r]\ar[d,"\cong"]&\mathcal E_2\ar[r]\ar[d,equal]&0\\
						0\ar[r]&\mathcal E_1\ar[r]&\mathcal E_1\oplus\mathcal E_2\ar[r]&\mathcal E_2\ar[r]&0\\
						0\ar[r]&V_1(-m)\ar[u,->>,"q_1"]\ar[r]&V(-m)\ar[u,->>,"q_1\oplus q_2"]\ar[r]&V_2(-m)\ar[u,->>,"q_2"]\ar[r]&0. 
					\end{tikzcd}
				\end{equation}

				Apply the functor $\mathcal F\mapsto \mathcal F\otimes_{\mathcal O_B}\mathcal L^m\mapsto f_*(\mathcal F\otimes_{\mathcal O_B}\mathcal L^m)$ to the diagram above, and we have the following 
				\begin{equation}
					\begin{tikzcd}
						0\ar[r]&V_1\ar[d,"f_*(q_1(m))"]\ar[r]&V\ar[d,"f_*(q(m))"]\ar[r]&V_2\ar[d,"f_*(q_2(m))"]\ar[r]&0\\
						0\ar[r]&f_*(\mathcal E_1(m))\ar[r]\ar[d,equal]&f_*(\mathcal E(m))\ar[r]\ar[d,"\mu"]&f_*(\mathcal E_2(m))\ar[r]\ar[d,equal]&0\\
						0\ar[r]&f_*(\mathcal E_1(m))\ar[r]&f_*(\mathcal E_1(m))\oplus f_*(\mathcal E_2(m))\ar[r]&f_*(\mathcal E_2(m))\ar[r]&0\\
						0\ar[r]&V_1\ar[u,"f_*(q_1(m))"]\ar[r]&V\ar[u,"f_*(q_1(m))\oplus f_*(q_2(m))"]\ar[r]&V_2\ar[u,"f_*(q_2(m))"]\ar[r]&0
					\end{tikzcd}
				\end{equation}
				where $\mathcal E(m):=\mathcal E\otimes_{\mathcal O_B}\mathcal L^m$ and $q(m):=q\otimes 1_{\mathcal L^m}$ and similarly for $\mathcal E_i(m)$, $q_i(m)$. We denote by $\mu:f_*(\mathcal E(m))\to f_*(\mathcal E_1(m))\oplus f_*(\mathcal E_2(m))$ the map induced from $\mathcal E\to \mathcal E_1\oplus\mathcal E_2$. We have that $f_*(q(m))$ and $f_*(q_i(m))$ are isomorphisms. 
				
				Let $g:V\to V$ be the composition 
				\begin{equation}
					g:=\big(f_*(q_1(m))^{-1}\oplus f_*(q_2(m))^{-1}\big)\circ \mu\circ f_*(q(m)). 
				\end{equation}
				Then $g\in\mathrm{GL}(V)$ and the following diagram commutes 
				\begin{equation}
					\begin{tikzcd}
						V(-m)\ar[r,"g\otimes 1_{\mathcal L^{-m}}"]\ar[d,->>,"q"]&V(-m)\ar[d,->>,"q_1\oplus q_2"]\\
						\mathcal E\ar[r,"\cong"{sloped}]&\mathcal E_1\oplus\mathcal E_2
					\end{tikzcd}
				\end{equation}
				i.e. $g.[(\mathcal E,q)]=\Big[\Big(\bigoplus_{i=1}^2\mathcal E_i,\bigoplus_{i=1}^2q_i\Big)\Big]$. Finally we have $g\in \mathrm U_\tau$, since the following diagram commutes
				\begin{equation}
					\begin{tikzcd}
						0\ar[r]&V_1\ar[r]\ar[d,equal]&V\ar[r]\ar[d,"g"]&V_2\ar[r]\ar[d,equal]&0\\
						0\ar[r]&V_1\ar[r]&V\ar[r]&V_2\ar[r]&0. 
					\end{tikzcd}
				\end{equation}
				This proves $(4)\iff(4')$. 
			\end{proof}
		}
	}

	\subsection{The moduli functor \texorpdfstring{$\mathbf M'_{\tau,\delta}$}{M'\_{tau,alpha}}}\label{subsection of the moduli functor M'_(tau,alpha) in length 2}
	{
		Recall that $(\tau,\delta)$-stability is defined in Definition \ref{definition: Jackson's (tau,alpha)-stability} and families of $(\tau,\delta)$-stable sheaves are defined in Definition \ref{definition: family of (tau,delta)-stable sheaves}. Recall that $f:B\to\mathrm{Spec}(\Bbbk)$ is the structure morphism. 

		Recall that the moduli functor $\mathbf M'_{\tau,\delta}$ is defined in section \ref{subsection: the moduli problem} which sends $T\in\mathrm{Sch}/\Bbbk$ to the set 
		\begin{equation}
			\mathbf M'_{\tau,\delta}(T):=\Big\{\begin{matrix}\textrm{isomorphism classes of flat families}\\\textrm{over }T\textrm{ of }(\tau,\delta)\textrm{-stable sheaves on }B\end{matrix}\Big\}. 
		\end{equation}
		
		Let $m\gg0$. Consider $Y_{\tau,\delta}^{\mathrm{s}}$ defined in \eqref{equation: definition of Y_(tau,delta)^s}. The Yoneda functor is a fully faithful embedding $\mathrm{Sch}/\Bbbk\to \mathrm{PSh}(\mathrm{Sch}/\Bbbk)$, and we will not distinguish between a scheme and its functor of points. We will define a $\mathrm P_\tau$-invariant natural transformation
		\begin{equation}
			\rho:Y_{\tau,\delta}^{\mathrm{s}}\to \mathbf M'_{\tau,\delta}. 
		\end{equation}
		Let $T\in\mathrm{Sch}/\Bbbk$, and let $T\to Y_{\tau,\delta}^{\mathrm{s}}$ be a morphism. By Lemma \ref{lemma: moduli descriptions of S_alpha(Y_tau^s) and Y_(tau,alpha)^s}, the morphism $T\to Y_{\tau,\delta}^{\mathrm{s}}$ corresponds to an isomorphism class $[(\mathcal E,q)]$ satisfying $(1)-(4)$ in Lemma \ref{lemma: moduli descriptions of S_alpha(Y_tau^s) and Y_(tau,alpha)^s}. Then the isomorphism class of $\mathcal E$ is an element in $\mathbf M'_{\tau,\delta}(T)$. This defines a map 
		\begin{equation}
			\rho_T:\mathrm{Hom}_{\mathrm{Sch}/\Bbbk}(T,Y_{\tau,\delta}^{\mathrm{s}})\to \mathbf M'_{\tau,\delta}(T)
		\end{equation}
		It is easy to see $\rho_T$ is $\mathrm{Hom}_{\mathrm{Sch}/\Bbbk}(T,\mathrm P_\tau)$-invariant and it is natural in $T\in\mathrm{Sch}/\Bbbk$. This defines a $\mathrm P_\tau$-invariant natural transformation $\rho:Y_{\tau,\delta}^{\mathrm{s}}\to \mathbf M'_{\tau,\delta}$.

		Consider the quotient presheaf 
		\begin{equation}
			Y_{\tau,\delta}^{\mathrm{s}}/_p\mathrm P_\tau:(\mathrm{Sch}/\Bbbk)^{\mathrm{op}}\to\mathrm{Set},\quad T\mapsto \frac{\mathrm{Hom}_{\mathrm{Sch}/\Bbbk}(T,Y_{\tau,\delta}^{\mathrm{s}})}{\mathrm{Hom}_{\mathrm{Sch}/\Bbbk}(T,\mathrm P_\tau)}
		\end{equation}
		where $/_p$ indicates it is a quotient presheaf. Since $\rho:Y_{\tau,\delta}^{\mathrm{s}}\to \mathbf M'_{\tau,\delta}$ is $\mathrm P_\tau$-invariant, it factors through $Y_{\tau,\delta}^{\mathrm{s}}/_p\mathrm P_\tau$, i.e. there exists a unique natural transformation $\theta$ such that the following diagram in $\mathrm{PSh}(\mathrm{Sch}/\Bbbk)$ commutes 
		\begin{equation}
			\begin{tikzcd}
				Y_{\tau,\delta}^{\mathrm{s}}\ar[r,"\rho"]\ar[d,"\varpi"]&\mathbf M'_{\tau,\delta}\\
				Y_{\tau,\delta}^{\mathrm{s}}/_p\mathrm P_\tau\ar[ru,"\theta"]
			\end{tikzcd}
		\end{equation}
		
		{
			\begin{lemma}\label{lemma: theta from Y_(tau,delta)^s/_pP_tau to M'_(tau,delta) sheafifies to an isomorphism}
				When $m\gg0$ is large, the natural transformation $\theta:Y_{\tau,\delta}^{\mathrm{s}}/_p\mathrm P_\tau\to\mathbf M'_{\tau,\delta}$ satisfies: 
				\begin{itemize}
					\item[(1)] $\theta:Y_{\tau,\delta}^{\mathrm{s}}/_p\mathrm P_\tau\to\mathbf M'_{\tau,\delta}$ is injective. 

					\item[(2)] For any $T\in\mathrm{Sch}/\Bbbk$ and any $x\in\mathbf M'_{\tau,\delta}(T)$, there exists an open covering $T=\bigcup_{j\in J}T_j$ such that $x|_{T_j}\in \mathbf M'_{\tau,\delta}(T_j)$ is in the image of $\theta_{T_j}:\big(Y_{\tau,\delta}^{\mathrm{s}}/_p\mathrm P_\tau\big)(T_j)\to\mathbf M'_{\tau,\delta}(T_j)$ for each $j\in J$, where $x|_{T_j}$ is the pullback of $x$ along $T_j\to T$. 
				\end{itemize}
				Moreover, the morphism $\theta$ induces an isomorphism of sheaves on $(\mathrm{Sch}/\Bbbk)_{\acute etale}$ 
				\begin{equation}
					\theta^\sharp:\big(Y_{\tau,\delta}^{\mathrm{s}}/_p\mathrm P_\tau\big)^\sharp\cong (\mathbf M'_{\tau,\delta})^\sharp
				\end{equation}
				where ${}^\sharp$ refers to the sheafification with respect to the étale topology. 
			\end{lemma}
			\begin{proof}
				For $i=1,2$ the collection of stable sheaves of Hilbert polynomial $P_i(t)$ on $B$ is bounded and then regularities of such stable sheaves are bounded. Let $m\gg0$ be larger than the regularity of any stable sheaf of Hilbert polynomial in $\{P_1(t),P_2(t)\}$. We will prove the lemma for this $m$. For $T\in\mathrm{Sch}/\Bbbk$, recall that $f_T:B_T\to T$ is the projection. Denote $\mathcal E(m):=\mathcal E\otimes_{\mathcal O_{B_T}}\mathcal L_T^m$. For any morphism $T'\to T$ in $\mathrm{Sch}/\Bbbk$, denote by $\mathcal E|_{B_{T'}}$ the pullback of $\mathcal E$ along $B_{T'}\to B_T$. 
				
				For $(1)$, let $T\in\mathrm{Sch}/\Bbbk$. It suffices to prove that if $\rho_T(h)=\rho_T(h')$ for two elements $h,h'\in\mathrm{Hom}_{\mathrm{Sch}/\Bbbk}(T,Y_{\tau,\delta}^{\mathrm{s}})$, then $h,h'$ are in one $\mathrm{Hom}_{\mathrm{Sch}/\Bbbk}(T,\mathrm P_\tau)$-orbit. By Lemma \ref{lemma: moduli descriptions of S_alpha(Y_tau^s) and Y_(tau,alpha)^s}, elements in $\mathrm{Hom}_{\mathrm{Sch}/\Bbbk}(T,Y_{\tau,\delta}^{\mathrm{s}})$ correspond to isomorphism classes $[(\mathcal E,q)]$ satisfying $(1)-(4)$ in Lemma \ref{lemma: moduli descriptions of S_alpha(Y_tau^s) and Y_(tau,alpha)^s}. Let $[(\mathcal E,q)]$ and $[(\mathcal E',q')]$ be two isomorphism classes corresponding to $h$ and $h'$ respectively 
				\begin{equation}
					q:V(-m)_T\to \mathcal E,\quad q':V(-m)_T\to\mathcal E'. 
				\end{equation}
				The equality $\rho_T(h)=\rho_T(h')$ implies $\mathcal E\cong\mathcal E'$ on $B_T$ and the isomorphism preserves the filtrations $\{\mathcal E_{\leq i}\}_{i=0}^2$ and $\{\mathcal E'_{\leq i}\}_{i=0}^2$. 
				
				Consider the following commutative diagram with exact rows
				\begin{equation}\label{equation: diagram of V_T to f_(T,*)(E(m))}
					\begin{tikzcd}
						0\ar[r]&V_{\leq i-1}\otimes_\Bbbk\mathcal O_T\ar[r]\ar[d,"f_{T,*}(q_{\leq i-1}(m))"]&V_{\leq i}\otimes_\Bbbk\mathcal O_T\ar[r]\ar[d,"f_{T,*}(q_{\leq i}(m))"]&V_i\otimes_\Bbbk\mathcal O_T\ar[d,"f_{T,*}(q_i(m))"]\ar[r]&0\\
						0\ar[r]&f_{T,*}(\mathcal E_{\leq i-1}(m))\ar[r]&f_{T,*}(\mathcal E_{\leq i}(m))\ar[r]&f_{T,*}(\mathcal E_i(m))\ar[r]&0. 
					\end{tikzcd}
				\end{equation}
				Recall that $m$ is larger than the regularity of $\mathcal E_i$ for $i=1,2$. We have that $f_{T,*}(\mathcal E_{\leq i}(m))$ and $f_{T,*}(\mathcal E_i(m))$ are locally free of ranks $P_{\leq i}(m)$ and $P_i(m)$ on $T$ for $i=1,2$ respectively. 
				Therefore $f_{T,*}(q_{\leq i}(m))$ and $f_{T,*}(q_i(m))$ in diagram \eqref{equation: diagram of V_T to f_(T,*)(E(m))} are morphisms between locally free sheaves. For any closed point $t\in T$, we have that $[(\mathcal E|_{B_t},q|_{B_t})]$ corresponds to a closed point in $Y_{\tau,\delta}^{\mathrm{s}}$, so 
				\begin{equation}
					\begin{split}
						f_{t,*}(q_i|_{B_t}(m)):&V_i\to f_{t,*}(\mathcal E_i|_{B_t}(m))\\
						f_{t,*}(q_{\leq i}|_{B_t}(m)):&V_{\leq i}\to f_{t,*}(\mathcal E_{\leq i}|_{B_t}(m))
					\end{split}
				\end{equation}
				are isomorphisms, i.e. fibres at $t\in T$ of $f_{T,*}(q_{\leq i}(m))$ and $f_{T,*}(q_i(m))$ in diagram \eqref{equation: diagram of V_T to f_(T,*)(E(m))} are isomorphisms. This proves that $f_{T,*}(q_{\leq i}(m))$ and $f_{T,*}(q_i(m))$ are isomorphisms by \href{https://stacks.math.columbia.edu/tag/00O0}{Lemma 00O0}
				\begin{equation}
					\begin{split}
						f_{T,*}(q_{\leq i}(m)):&V_{\leq i}\otimes_\Bbbk\mathcal O_T\cong f_{T,*}(\mathcal E_{\leq i}(m))\\
						f_{T,*}(q_i(m)):&V_{\leq i}\otimes_\Bbbk\mathcal O_T\cong f_{T,*}(\mathcal E_i(m)). 
					\end{split}
				\end{equation}
				Similarly, we have that $f_{T,*}(q'_{\leq i}(m))$ and $f_{T,*}(q'_i(m))$ are isomorphisms for $i=1,2$. 
				
				Consider the morphism $g:V\otimes_\Bbbk\mathcal O_T\to V\otimes_\Bbbk\mathcal O_T$ defined as the following composition on $T$ 
				\begin{equation}
					\begin{tikzcd}
						V\otimes_\Bbbk\mathcal O_T\ar[rr,"f_{T,*}(q(m))"]\ar[rrrr,bend right=15,"g",yshift=-5]&&f_{T,*}(\mathcal E(m))\ar[rr,"f_{T,*}(q'(m))^{-1}"]&&V\otimes_\Bbbk\mathcal O_T. 
					\end{tikzcd}
				\end{equation}
				It is easy to check: 
				\begin{itemize}
					\item the morphism $g:V\otimes_\Bbbk\mathcal O_T\to V\otimes_\Bbbk\mathcal O_T$ preserves the filtration $\{V_{\leq i}\}_{i=0}^2$, thus it defines a morphism $g\in\mathrm{Hom}_{\mathrm{Sch}/\Bbbk}(T,\mathrm P_{\mathrm{GL},\tau})$, where $\mathrm P_{\mathrm{GL},\tau}\subseteq \mathrm{GL}(V)$ is the parabolic subgroup stabilising the flag $\{V_{\leq i}\}_{i=0}^2$; 

					\item the following diagram commutes 
					\begin{equation}
						\begin{tikzcd}
							V(-m)_T\ar[r,"g(-m)"]\ar[d,"q"]&V(-m)_T\ar[d,"q'"]\\
							\mathcal E\ar[r,equal]&\mathcal E
						\end{tikzcd}
					\end{equation}
					that is $g.[(\mathcal E,q)]=[(\mathcal E',q')]$. 
				\end{itemize}
				Thus $h,h'\in\mathrm{Hom}_{\mathrm{Sch}/\Bbbk}(T,Y_{\tau,\delta}^{\mathrm{s}})$ are in one $\mathrm{Hom}_{\mathrm{Sch}/\Bbbk}(T,\mathrm P_{\mathrm{GL},\tau})$-orbit, which is a $\mathrm{Hom}_{\mathrm{Sch}/\Bbbk}(T,\mathrm P_{\tau})$-orbit, since the centre $\mathbb G_m\subseteq\mathrm{GL}(V)$ acts trivially. This proves $(1)$.

				For $(2)$, let $T\in\mathrm{Sch}/\Bbbk$ and let $x\in\mathbf M'_{\tau,\delta}(T)$. Let $\mathcal E$ represent the isomorphism class $x\in\mathbf M'_{\tau,\delta}(T)$. Note that $\mathcal E$ and $\mathcal E_{\leq i}$ are coherent sheaves on $B_T$ and $\mathcal E_i:=\mathcal E_{\leq i}/\mathcal E_{\leq i-1}$ is flat over $T$ with Hilbert polynomial $P_i(t)$. 
				
				Since $m$ is larger than any regularity, we have that: 
				\begin{itemize}
					\item $f_{T,*}(\mathcal E_{\leq i}(m))$ and $f_{T,*}(\mathcal E_i(m))$ are locally free of ranks $P_{\leq i}(m)$ and $P_i(m)$ on $T$ for $i=1,2$ respectively; 
					\item $R^jf_{T,*}(\mathcal E_{\leq i}(m))=0$ and $R^jf_{T,*}(\mathcal E_i(m))=0$ for $i=1,2$ and $j>0$; 
					\item $(f_T)^*f_{T,*}(\mathcal E_{\leq i}(m))\to\mathcal E_{\leq i}(m)$ and $(f_T)^*f_{T,*}(\mathcal E_i(m))\to \mathcal E_i(m)$ are surjective for $i=1,2$. 
				\end{itemize}

				There is an open covering $T=\bigcup_{j\in J}T_j$ such that $f_{T,*}(\mathcal E_i(m))|_{T_j}$ is free of rank $P_i(m)$ for $i=1,2$. Then $f_{T,*}(\mathcal E_{\leq i}(m))|_{T_j}$ is free of rank $P_{\leq i}(m)$ for $i=1,2$. On each $T_j$, we can choose generating global sections compatible with the filtration $\big\{f_{T,*}(\mathcal E_{\leq i}(m))\big\}_{i=0}^2$ and obtain isomorphisms in the following diagram 
				\begin{equation}
					\begin{tikzcd}
						0\ar[r,"\subseteq"]&V_{\leq 1}\otimes_\Bbbk \mathcal O_{T_j}\ar[r,"\subseteq"]\ar[d,"\cong"]&V_{\leq 2}\otimes_\Bbbk\mathcal O_{T_j}\ar[r,equal]\ar[d,"\cong"]&V\otimes_\Bbbk\mathcal O_{T_j}\ar[d,"\cong"]\\
						0\ar[r,"\subseteq"]&f_{T,*}(\mathcal E_{\leq 1}(m))|_{T_j}\ar[r,"\subseteq"]&f_{T,*}(\mathcal E_{\leq 2}(m))|_{T_j}\ar[r,equal]&f_{T,*}(\mathcal E(m))|_{T_j}\;.
					\end{tikzcd}
				\end{equation}
				
				On each $T_j$ and $i=1,2$ we have a surjective morphism
				\begin{equation}
					\begin{split}
						V_{\leq i}(-m)_{T_j}:=&(f_{T_j})^*\big(V_{\leq j}\otimes_\Bbbk\mathcal O_{T_j}\big)\otimes_{\mathcal O_{B_{T_j}}}\mathcal L_{T_j}^{-m}\\
						\cong&(f_{T_j})^*\big(f_{T,*}(\mathcal E_{\leq i}(m))|_{T_j}\big)\otimes_{\mathcal O_{B_{T_j}}}\mathcal L_{T_j}^{-m}\\
						\cong&(f_T)^*f_{T,*}(\mathcal E_{\leq i}(m))\big|_{B_{T_j}}\otimes_{\mathcal O_{B_{T_j}}}\mathcal L_{T_j}^{-m}\\
						\twoheadrightarrow&\mathcal E_{\leq i}(m)\big|_{B_{T_j}}\otimes_{\mathcal O_{B_{T_j}}}\mathcal L_{T_j}^{-m}\\
						\cong&\mathcal E_{\leq i}|_{B_{T_j}}. 
					\end{split}
				\end{equation}
				When $i=2$, let $q_j:V(-m)_{T_j}\to \mathcal E|_{B_{T_j}}$ denote the surjective morphism. Then we have $q_j\big(V_{\leq i}(-m)_{B_{T_j}}\big)=\mathcal E_{\leq i}|_{B_{T_j}}$ for $i=1,2$. Therefore $[(\mathcal E|_{B_{T_j}},q_j)]$ defines a morphism $T_j\to Y_\tau$. It is easy to check that it defines a morphism $h_j:T_j\to Y_{\tau,\delta}^{\mathrm{s}}$, for example by Lemma \ref{lemma: moduli descriptions of S_alpha(Y_tau^s) and Y_(tau,alpha)^s}. We have that $\rho_{T_j}$ maps $h_j$ to $x|_{T_j}=[\mathcal E|_{B_{T_j}}]\in\mathbf M'_{\tau,\delta}(T_j)$. Since the image of $\rho_{T_j}$ and the image of $\theta_{T_j}$ coincide, we have that $x|_{T_j}$ is in the image of $\theta_{T_j}$. This proves $(2)$. 

				Let $\theta^\sharp:\big(Y_{\tau,\delta}^{\mathrm{s}}/_p\mathrm P_\tau\big)^\sharp\to (\mathbf M'_{\tau,\delta})^\sharp$ be sheafification of $\theta$ on $(\mathrm{Sch}/\Bbbk)_{\acute etale}$. We have that $\theta^\sharp$ is injective since $\theta$ is by $(1)$. Zariski coverings are étale coverings, so $\theta^\sharp$ is surjective as a morphism of sheaves by $(2)$. Then $\theta^\sharp$ is an isomorphism. 
			\end{proof}
		}
	}

	\subsection{Local properties of the quotient}
	{
		Recall in Corollary \ref{corollary: geometric U-quotients represent quotient sheaves with UU}, with Condition \hyperref[condition: UU]{UU}, the quotient sheaf is represented by the geometric quotient, and the representability commutes with base change. We will prove the quotient $\pi:Y_{\tau,\delta}^{\mathrm{s}}\to Y_{\tau,\delta}^{\mathrm{s}}/\mathrm P_\tau$ represents the associated quotient presheaf. Therefore $(\mathbf M'_{\tau,\delta})^\sharp\cong Y_{\tau,\delta}^{\mathrm{s}}/\mathrm P_\tau$ by Lemma \ref{lemma: theta from Y_(tau,delta)^s/_pP_tau to M'_(tau,delta) sheafifies to an isomorphism}. 

		{
			\begin{lemma}\label{lemma: reductive stabiliser in the unipotent quotient is central}
				Let $m\gg 0$. Let $[(\mathcal E,q)]$ be a closed point of $Y_\tau^{\mathrm{s}}$. Assume $\mathcal E\not\cong \mathcal E_1\oplus\mathcal E_2$. If $(u,t,k)\in\mathrm U_\tau\rtimes(\mathrm T_\tau\times\mathrm K_\tau)$ stabilises $[(\mathcal E,q)]$, then the image of $(t,k)$ along $\mathrm T_\tau\times \mathrm K_\tau\to\mathrm P_\tau$ is a scalar multiplication of $1_V\in \mathrm P_\tau$. 
			\end{lemma}
			\begin{proof}
				Let $q:V(-m)\to \mathcal E$ be a quotient in the class $[(\mathcal E,q)]$. Let $p_\tau^{\mathrm{s}}:Y_\tau^{\mathrm{s}}\to \prod_{i=1}^2R_i^{\mathrm{s}}$ be the retraction. See Theorem \ref{theorem: simpson's construction of moduli of semistable sheaves} for $R_i^{\mathrm{s}}$. Let $x:=[(\mathcal E,q)]$ denote the point of $Y_\tau^{\mathrm{s}}$. 

				We have that $p_\tau^{\mathrm{s}}$ is $\mathrm U_\tau$-invariant and $\mathrm T_\tau\times \mathrm K_\tau$-equivariant. Since $(u,t,k)$ stabilises $x\in Y_\tau^{\mathrm{s}}$, we have
				\begin{equation}
					p_\tau^{\mathrm{s}}(x)=p_\tau^{\mathrm{s}}\big((u,t,k).x\big)=p_\tau^{\mathrm{s}}\big((t,k).x\big)=(t,k).p_\tau^{\mathrm{s}}(x). 
				\end{equation}
				Thus $(t,k)$ stabilises $p_\tau^{\mathrm{s}}(x)$. 

				Let $\begin{pmatrix}g_1&\\&g_2\end{pmatrix}\in\mathrm P_\tau$ be the image of $(t,k)\in\mathrm T_\tau\times\mathrm K_\tau$. Then $g_i\in\mathrm{GL}(V_i)$ stabilises $[(\mathcal E_i,q_i)]\in R_i^{\mathrm{s}}$ for $i=1,2$. By \cite{HuybrechtsDaniel2010Tgom} Lemma 4.3.2, we have that $g_i$ is in the image of $\mathrm{Aut}(\mathcal E_i)\to \mathrm{GL}(V_i)$. Since $\mathcal E_i$ is stable, we have $\mathrm{Aut}(\mathcal E_i)\cong \Bbbk^\times$. Then $g_i=c_i1_{V_i}$ for some $c_i\in\Bbbk^\times$. 

				We need to prove $c_1=c_2$ when $\mathcal E\not\cong\mathcal E_1\oplus\mathcal E_2$. Suppose otherwise $c_1\ne c_2$. Let $\begin{pmatrix}1_{V_1}&A\\&1_{V_2}\end{pmatrix}\in \mathrm P_\tau$ be the image of $u\in\mathrm U_\tau$, where $A:V_2\to V_1$ is a linear map. The image of $(u,t,k)$ in $\mathrm P_\tau$ is $\begin{pmatrix}c_11_{V_1}&c_2A\\&c_21_{V_2}\end{pmatrix}$. Let $g:V\to V$ denote this linear map. Consider another linear map 
				\begin{equation}
					\sigma:V_2\to V,\quad v_2\mapsto \frac{c_2Av_2}{c_2-c_1}+v_2. 
				\end{equation}
				If $\mathrm{pr}_2:V\to V_2$ denotes the projection, then 
				\begin{equation}
					\sigma\circ\mathrm{pr}_2=\frac{1}{c_2-c_1}g-\frac{c_1}{c_2-c_1}1_V,\quad \textrm{in }\mathrm{End}_\Bbbk(V). 
				\end{equation}

				The linear maps $g,\sigma,\mathrm{pr}_2$ induce morphisms 
				\begin{equation}
					\begin{split}
						g(-m):&V(-m)\to V(-m)\\
						\sigma(-m):&V_2(-m)\to V(-m)\\
						\mathrm{pr}_2(-m)&:V(-m)\to V_2(-m)
					\end{split}
				\end{equation}
				where $g(-m):=g\otimes_\Bbbk 1_{\mathcal L^m}$ and $\sigma(-m),\mathrm{pr}_2(-m)$ are defined similarly. They satisfy 
				\begin{equation}
					\sigma(-m)\circ \mathrm{pr}_2(-m)=\frac{1}{c_2-c_1}g(-m)-\frac{c_1}{c_2-c_1}1_{V(-m)}. 
				\end{equation}

				Since $g.x=x$, we have $g(-m)(\ker q)\subseteq\ker(q)$ by \cite{HuybrechtsDaniel2010Tgom} Lemma 4.3.2. By the above equation we have
				\begin{equation}
					\sigma(-m)\big(\mathrm{pr}_2(-m)(\ker q)\big)\subseteq \ker q. 
				\end{equation}
				Since $\mathrm{pr}_2(-m)(\ker q)=\ker q_2$, we have $\sigma(-m)(\ker q_2)\subseteq\ker q$. Then $\sigma(-m)$ induces a morphism $\bar\sigma:\mathcal E_2\to\mathcal E$ fitting in the following diagram 
				\begin{equation}
					\begin{tikzcd}
						V_2(-m)\ar[r,"\sigma(-m)"]\ar[d,->>,"q_2"]&V(-m)\ar[d,->>,"q"]\\
						\mathcal E_2\ar[r,dashed,"\bar\sigma"]&\mathcal E. 
					\end{tikzcd}
				\end{equation}

				Since $\mathrm{pr}_2\circ\sigma=1_{V_2}$, we have $\overline{\mathrm{pr}_2}\circ\bar\sigma=1_{\mathcal E_2}$, where $\overline{\mathrm{pr}_2}:\mathcal E\to\mathcal E_2$ is the quotient. Then $\mathcal E\to\mathcal E_2$ has a section, i.e. $\mathcal E\cong \mathcal E_1\oplus\mathcal E_2$, which is a contradiction. 
			\end{proof}
		}

		{
			\begin{lemma}\label{lemma: quotient of Y_(tau,delta)^s by P_tau represents quotient sheaf}
				The quotient $\pi:Y_{\tau,\delta}^{\mathrm{s}}\to Y_{\tau,\delta}^{\mathrm{s}}/\mathrm P_\tau$ represents the quotient sheaf $\big(Y_{\tau,\delta}^{\mathrm{s}}/_p\mathrm P_\tau\big)^\sharp$ on $(\mathrm{Sch}/\Bbbk)_{\acute etale}$. 
			\end{lemma}
			\begin{proof}
				Let $Z(\mathrm P_\tau)\subseteq \mathrm P_\tau$ denote the centre. Then $Z(\mathrm P_\tau)=\{t1_V\in \mathrm P_\tau:t\in\mathbb G_m\}$. The action of $\mathrm P_\tau$ on $Y_{\tau,\delta}^{\mathrm{s}}$ factors through $\mathrm P_\tau/Z(\mathrm P_\tau)$. Therefore $Y_{\tau,\delta}^{\mathrm{s}}$ has a quasi-projective universal geometric quotient by $\mathrm P_\tau/Z(\mathrm P_\tau)$. 

				We can think the quotient of $Y_{\tau,\delta}^{\mathrm{s}}$ by $\mathrm P_\tau/Z(\mathrm P_\tau)$ as a composition of universal geometric quotients
				\begin{equation}
					Y_{\tau,\delta}^{\mathrm{s}}\to Y_{\tau,\delta}^{\mathrm{s}}/\mathrm U_\tau\to \big(Y_{\tau,\delta}^{\mathrm{s}}/\mathrm U_\tau\big)\big/ \big(\mathrm P_\tau/\mathrm U_\tau Z(\mathrm P_\tau)\big). 
				\end{equation}

				These two quotients both represent their quotient sheaves: 
				\begin{itemize}
					\item[(1)] The first quotient $Y_{\tau,\delta}^{\mathrm{s}}\to Y_{\tau,\delta}^{\mathrm{s}}/\mathrm U_\tau$ represents the quotient sheaf $\big(Y_{\tau,\delta}^{\mathrm{s}}/_p\mathrm U_\tau\big)^\sharp$ by Corollary \ref{corollary: geometric U-quotients represent quotient sheaves with UU}; 
					\item[(2)] The second quotient $Y_{\tau,\delta}^{\mathrm{s}}/\mathrm U_\tau\to \big(Y_{\tau,\delta}^{\mathrm{s}}/\mathrm U_\tau\big)\big/ \big(\mathrm P_\tau/\mathrm U_\tau Z(\mathrm P_\tau)\big)$ is a geometric quotient of a reductive group, and closed points have trivial stabilisers by Lemma \ref{lemma: reductive stabiliser in the unipotent quotient is central}. We can apply Luna's étale slice theorem to conclude that it is a principal $\mathrm P_\tau/\mathrm U_\tau Z(\mathrm P_\tau)$-bundle by \cite{HuybrechtsDaniel2010Tgom} Corollary 4.2.13. In particular, it also represents the associated quotient sheaf. 
				\end{itemize}

				As a composition of $(1)$ and $(2)$, we have that $\pi:Y_{\tau,\delta}^{\mathrm{s}}\to Y_{\tau,\delta}^{\mathrm{s}}/\mathrm P_\tau$ induces an isomorphism $Y_{\tau,\delta}^{\mathrm{s}}/\mathrm P_\tau\cong \big(Y_{\tau,\delta}^{\mathrm{s}}/_p\mathrm P_\tau\big)^\sharp$. 
			\end{proof}
		}

		{
			\begin{theorem}\label{theorem: sheafification of M'_(tau,delta) is representable}
				The sheafification of $\mathbf M'_{\tau,\delta}$ on $(\mathrm{Sch}/\Bbbk)_{\acute etale}$ is representable by $Y_{\tau,\delta}^{\mathrm{s}}/\mathrm P_\tau$ 
				\begin{equation}
					(\mathbf M'_{\tau,\delta})^\sharp\cong Y_{\tau,\delta}^{\mathrm{s}}/\mathrm P_\tau. 
				\end{equation}
				In particular, the moduli functor $\mathbf M'_{\tau,\delta}$ has a coarse moduli space $Y_{\tau,\delta}^{\mathrm{s}}/\mathrm P_\tau$. 
			\end{theorem}
			\begin{proof}
				By Lemma \ref{lemma: theta from Y_(tau,delta)^s/_pP_tau to M'_(tau,delta) sheafifies to an isomorphism}, the sheafification of $\theta$ is an isomorphism $\theta^\sharp: \big(Y_{\tau,\delta}^{\mathrm{s}}/_p\mathrm P_\tau\big)^\sharp \cong (\mathbf M'_{\tau,\delta})^\sharp$. By Lemma \ref{lemma: quotient of Y_(tau,delta)^s by P_tau represents quotient sheaf}, we have $\big(Y_{\tau,\delta}^{\mathrm{s}}/_p\mathrm P_\tau\big)^\sharp \cong Y_{\tau,\delta}^{\mathrm{s}}/P_\tau$. Then $(\mathbf M'_{\tau,\delta})^\sharp\cong Y_{\tau,\delta}^{\mathrm{s}}/\mathrm P_\tau$. 
			\end{proof}
		}

	}
}

\appendix
{
	\section{On the boundaries of \texorpdfstring{$Y_\tau$}{Y\_tau} and \texorpdfstring{$Y_\tau^{\mathrm{s}}$}{Y\_tau\^{}s}}\label{appendix: on the boundaries of Y_tau and Y_tau^s}
	{
		Let $\tau$ be a Harder-Narasimhan type of length 2. Let $P_i(t),P(t)$ be the Hilbert polynomials associated to $\tau$. For $1\leq i\leq 2$, let $P_{\leq i}(t):=P_1(t)+\cdots+P_i(t)$. Recall that $\mathrm P_\tau\subseteq\mathrm{SL}(V)$ is the parabolic subgroup stabilising the flag $\{V_{\leq i}\}_{i=0}^2$ (See the first paragraph of Section \ref{subsubsection of moduli of unstable sheaves}). Recall that $Y_\tau\subseteq\mathrm{Quot}_{V(-m)/B/\Bbbk}^{P(t),\mathcal L}$ is the locally closed subscheme that represents the sub-functor $\mathbf Y_\tau$ (See Lemma \ref{lemma: representability of Y_tau subset Quot}). There is a retraction $p_\tau:Y_\tau\to \prod_{i=1}^2\mathrm{Quot}_{V_i(-m)/B/\Bbbk}^{P_i(t),\mathcal L}$ (See diagram \eqref{equation: diagram of retraction p_tau on Y_tau}), and $Y_\tau^{\mathrm{s}}\subseteq Y_\tau$ is the preimage of $\prod_{i=1}^2R_i^{\mathrm{s}}\subseteq \prod_{i=1}^2\mathrm{Quot}_{V_i(-m)/B/\Bbbk}^{P_i(t),\mathcal L}$ under $p_\tau$ (See diagram \eqref{equation: defining diagram of Y_tau^(s)s}). 

		We want to consider the $\mathrm P_\tau$-quotient of $Y_\tau^{\mathrm{s}}$. However, $Y_\tau^{\mathrm{s}}$ is not projective and our non-reductive GIT requires a linear action on a projective scheme. Therefore we need to consider the $\mathrm P_\tau$-action on a suitable projective completion $\overline{Y_\tau^{\mathrm{s}}}$. We will take $\overline{Y_\tau^{\mathrm{s}}}$ to be the scheme theoretic image of $Y_\tau^{\mathrm{s}}\to \mathrm{Quot}_{V(-m)/B/\Bbbk}^{P(t),\mathcal L}$. The main result of this section is Proposition \ref{proposition: Y_tau^s closure and its Z_min^K-s}, which says that we do not need to worry about the boundary of $Y_\tau^{\mathrm{s}}\subseteq \overline{Y_\tau^{\mathrm{s}}}$. 

		For $T\in\mathrm{Sch}/\Bbbk$ and a coherent sheaf $\mathcal F$ on $B$, recall the notation
		\begin{equation}
			\begin{split}
				&B_T:=B\times_\Bbbk T\\
				&f_T:B_T\to T\quad\textrm{the projection}\\
				&\mathcal F_T\quad \textrm{the pullback of }\mathcal F\textrm{ along }B_T\to B. 
			\end{split}
		\end{equation}
		For a family $(\mathcal E,q)$ of quotients parametrised by $T$, recall the notation
		\begin{equation}
			\begin{split}
				&\mathcal E_{\leq i}:=q\big(V_{\leq i}(-m)_T\big)\\
				&\mathcal E_i:=\mathcal E_{\leq i}/\mathcal E_{\leq i-1}\\
				&q_i:V_i(-m)_T\to \mathcal E_i\quad\textrm{the induced surjective morphism}. 
			\end{split}
		\end{equation}

		Recall that there exist closed immersions
		\begin{equation}
			\begin{split}
				\mathrm{Quot}_{V(-m)/B/\Bbbk}^{P(t),\mathcal L}&\hookrightarrow \mathrm{Grass}\big(V\otimes_\Bbbk H,P(M)\big)\\
				&\hookrightarrow\mathbb P\bigg(\Big(\bigwedge^{P(M)}(V\otimes_\Bbbk H)\Big)^*\bigg)
			\end{split}
		\end{equation}
		where $H:=H^0(B,\mathcal L^{M-m})=f_*(\mathcal L^{M-m})$. Note that points in the Grassmannian are quotients while points in the projective space are one-dimensional subspaces. Recall that $\mathcal M$ denotes a very ample line bundle on $\mathrm{Quot}_{V(-m)/B/\Bbbk}^{P(t),\mathcal L}$, which is the pullback of $\mathcal O(1)$ along the immersions above. 

		{
			Let $\overline{Y_\tau}\hookrightarrow \mathrm{Quot}_{V(-m)/B/\Bbbk}^{P(t),\mathcal L}$ be the scheme theoretic image of $Y_\tau\to \mathrm{Quot}_{V(-m)/B/\Bbbk}^{P(t),\mathcal L}$. 
			\begin{proposition}\label{lemma: Y_tau closure and Y_tau}
				For $M\gg m\gg0$, there exists a closed subscheme 
				\begin{equation}
					\mathbb P^{N'}\hookrightarrow \mathbb P\bigg(\Big(\bigwedge^{P(M)}(V\otimes_\Bbbk H)\Big)^*\bigg)
				\end{equation}
				such that: 
				\begin{itemize}
					\item[(1)] $\mathbb P^{N'}$ is a $\mathrm T_\tau$-invariant linear subspace; 
					\item[(2)] the immersion $Y_\tau\to \mathbb P\Big(\big(\bigwedge^{P(M)}(V\otimes_\Bbbk H)\big)^*\Big)$ factors through $\mathbb P^{N'}$, i.e. there exists a closed immersion $\overline{Y_\tau}\hookrightarrow \mathbb P^{N'}$; 
					\item[(3)] $\big(\mathbb P^{N'}\big)^0_{\min}\cap \overline{Y_\tau}=Y_\tau$, where $\lambda:\mathrm T_\tau\to\mathrm P_\tau$. 
				\end{itemize}
			\end{proposition}
			\begin{proof}
				Only in this proof, for simplicity we denote 
				\begin{equation}
					P:=P(M),\quad P_i:=P_i(M),\quad P_{\leq i}:=P_{\leq i}(M)
				\end{equation}
				and we simplify $(-)\otimes_\Bbbk(-)$ to $(-)\otimes(-)$
				
				We first construct the linear subspace $\mathbb P^{N'}$, which is determined by a quotient space of $\bigwedge^P(V\otimes H)$. We have 
				\begin{equation}
					\begin{split}
						\bigwedge^P(V\otimes H)=&\bigwedge^P\Big(\bigoplus_{i=1}^2(V_i\otimes H)\Big)\\
						=&\bigoplus_{(m_1,m_2)\in\Lambda}\bigwedge^{m_1}(V_1\otimes H)\otimes\;\bigwedge^{m_2}(V_2\otimes H)
					\end{split}
				\end{equation}
				where the direct sum is over the index set $\Lambda$ 
				\begin{equation}
					\Lambda:=\big\{(m_1,m_2)\in\mathbb N^2:m_1+m_2=P(M)\big\}. 
				\end{equation}
				Let $\Lambda_0\subseteq\Lambda$ be the subset 
				\begin{equation}
					\Lambda_0:=\big\{(m_1,m_2)\in\Lambda:m_1>P_1(M)\big\}. 
				\end{equation}
				and let $W\subseteq \bigwedge^P(V\otimes H)$ be the following subspace 
				\begin{equation}
					W:=\bigoplus_{(m_1,m_2)\in\Lambda_0}\bigwedge^{m_1}(V_1\otimes H)\otimes\;\bigwedge^{m_2}(V_2\otimes H). 
				\end{equation}
				Then the following quotient space of $\bigwedge^P(V\otimes H)$
				\begin{equation}
					\begin{tikzcd}
						\bigwedge^P(V\otimes H)\ar[r]&\bigwedge^P(V\otimes H)\Big/W\ar[r]&0
					\end{tikzcd}
				\end{equation}
				defines a linear subspace 
				\begin{equation}
					\mathbb P^{N'}\hookrightarrow \mathbb P\Big(\big(\bigwedge^P(V\otimes H)\big)^*\Big). 
				\end{equation}
				Since $W\subseteq \bigwedge^P(V\otimes H)$ is invariant under the torus $\mathrm T_\tau$, we have that $\mathbb P^{N'}$ is invariant under $\mathrm T_\tau$. This proves $(1)$. 
				
				Next, we prove $(2)$. Assume first we have proved that $Y_\tau\to \mathbb P\Big(\big(\bigwedge^P(V\otimes H)\big)^*\Big)$ factors through $\mathbb P^{N'}$, then we have the following commutative diagram 
				\begin{equation}
					\begin{tikzcd}
						Y_\tau\ar[r]\ar[d]&\mathrm{Quot}_{V(-m)/B/\Bbbk}^{P(t),\mathcal L}\ar[d,hook]\\
						\mathbb P^{N'}\ar[r,hook]&\mathbb P\Big(\big(\bigwedge^P(V\otimes H)\big)^*\Big). 
					\end{tikzcd}
				\end{equation}
				It is easy to see that if we have morphisms $T\to Z_1\hookrightarrow Z_2$ where $Z_1\hookrightarrow Z_2$ is a closed immersion, then the scheme theoretic image of $T\to Z_1$ is canonically isomorphic to the scheme theoretic image of $T\to Z_2$. Apply this observation to the diagram above. Then $\overline{Y_\tau}$ is the scheme theoretic image of all three morphisms from $Y_\tau$ in the diagram. In particular, we have $\overline{Y_\tau}\hookrightarrow \mathbb P^{N'}$, which is the second part of $(2)$. 
				
				The first part of $(2)$ is equivalent to the following 
				\begin{itemize}
					\item[(2')] for any $T\in\mathrm{Sch}/\Bbbk$ and any $T\to Y_\tau$, the composition $T\to Y_\tau\to \mathbb P\Big(\big(\bigwedge^P(V\otimes H)\big)^*\Big)$ factors through $\mathbb P^{N'}$. 
				\end{itemize}
				
				Let $T\in\mathrm{Sch}/\Bbbk$ and let $T\to Y_\tau$ be a morphism. Let $(\mathcal E,q)$ be a family over $T$ of quotients whose classifying morphism is $T\to Y_\tau$. The composition $T\to Y_\tau\to \mathbb P\Big(\big(\bigwedge^P(V\otimes H)\big)^*\Big)$ is the classifying morphism of the following surjective morphism 
				\begin{equation}
					\begin{tikzcd}
						\bigwedge^P(V\otimes H)\otimes \mathcal O_T\ar[r,"\bigwedge^P\rho"]&\bigwedge^Pf_{T,*}\big(\mathcal E(M)\big)\ar[r]&0
					\end{tikzcd}
				\end{equation}
				where 
				\begin{itemize}
					\item $\mathcal E(M):=\mathcal E\otimes_{\mathcal O_{B_T}}\mathcal L_T^{M}$; 
					\item $\rho:(V\otimes H)\otimes \mathcal O_T\to f_{T,*}\big(\mathcal E(M)\big)$ is the morphism 
					\begin{equation}
						\rho:=f_{T,*}\Big(q\otimes_{\mathcal O_{B_T}}1_{\mathcal L_T^{M}}\Big); 
					\end{equation}
					\item $f_{T,*}\big(\mathcal E(M)\big)$ is locally free of rank $P=P(M)$ on $T$ when $M$ is sufficiently large. 
				\end{itemize}
				
				The morphism $T\to \mathbb P\Big(\big(\bigwedge^P(V\otimes H)\big)^*\Big)$ factors through $\mathbb P^{N'}$ if and only if the following composition is zero 
				\begin{equation}
					\begin{tikzcd}
						W\otimes\mathcal O_T\ar[r,"\subseteq"]&\bigwedge^P(V\otimes H)\otimes \mathcal O_T\ar[r,"\bigwedge^P\rho"]&\bigwedge^Pf_{T,*}\big(\mathcal E(M)\big). 
					\end{tikzcd}
				\end{equation}
				
				Since the isomorphism class of $(\mathcal E,q)$ corresponds to a morphism to $Y_\tau$, we have that $\mathcal E_{\leq i}$ is flat over $T$ with Hilbert polynomial $P_{\leq i}(t)$. We can enlarge $M$ and assume that $f_{T,*}\big(\mathcal E_{\leq i}(M)\big)$ is locally free of rank $P_{\leq i}$ for all $1\leq i\leq 2$ and all $T\in\mathrm{Sch}/\Bbbk$ and all families over $T$ of quotients $(\mathcal E,q)$. 
				
				Let $(m_1,m_2)\in\Lambda$. Choose 
				\begin{equation}
					\begin{split}
						s^{(1)}_1,&\cdots,s^{(1)}_{m_1}\in V_1\otimes H\\
						S^{(2)}&\in\bigwedge^{m_2}(V_2\otimes H). 
					\end{split}
				\end{equation}
				Then vectors of the following form span the vector space $\bigwedge^{m_1}(V_1\otimes H)\otimes\;\bigwedge^{m_2}(V_2\otimes H)$
				\begin{equation}
					v:=\big(s^{(1)}_1\wedge\cdots\wedge s^{(1)}_{m_1}\big)\otimes S^{(2)},\quad (m_1,m_2)\in\Lambda. 
				\end{equation}
				If we consider all $(m_1,m_2)\in\Lambda_0$, then vectors of the above form span $W$. Therefore to show $\Big(\bigwedge^P\rho\Big)\big(W\otimes \mathcal O_T\big)=0$, it suffices to show that $v\otimes 1\in W\otimes\mathcal O_T$ maps to zero under $\bigwedge^P\rho$. 
				
				The image of $v\otimes 1$ under $\bigwedge^P\rho$ is 
				\begin{equation}
					\Big(\rho\big(s^{(1)}_1\otimes 1\big)\wedge\cdots\wedge \rho\big(s^{(1)}_{m_1}\otimes 1\big)\Big)\wedge\; \rho\big(S^{(2)}\otimes1\big). 
				\end{equation}
				We have 
				\begin{equation}
					\rho\big(s^{(1)}_1\otimes 1\big)\wedge\cdots\wedge \rho\big(s^{(1)}_{m_1}\otimes 1\big)=0
				\end{equation}
				since it is a section of the sheaf 
				\begin{equation}
					\bigwedge^{m_1}f_{T,*}\big(\mathcal E_1(M)\big)=0
				\end{equation}
				which is the zero sheaf since $f_{T,*}\big(\mathcal E_1(M)\big)$ is locally free of rank $P_1(M)<m_1$. This proves $(2)'$, thus $(2)$. 
				
				For $(3)$, let $w_1>w_2$ be the integers such that $w_1\dim V_1+w_2\dim V_2=0$ and $\lambda(t)=t^{w_1}1_{V_1}\oplus t^{w_2}1_{V_2}$. By \href{https://stacks.math.columbia.edu/tag/01QV}{Lemma 01QV} and \href{https://stacks.math.columbia.edu/tag/01R8}{Lemma 01R8}, we have that $Y_\tau\subseteq\overline{Y_\tau}$ is an open subscheme. Then both sides of $(3)$ are open subschemes of $\overline{Y_\tau}$, so it suffices to show that they contain the same closed points. 
				
				Recall that 
				\begin{equation}
					\begin{split}
						H^0\big(\mathbb P^{N'},\mathcal O(1)\big)\cong &\bigwedge^P(V\otimes H)\Big/W\\
						\cong& \bigoplus_{(m_1,m_2)\in\Lambda\setminus\Lambda_0}\bigwedge^{m_1}(V_1\otimes H)\otimes\; \bigwedge^{m_2}(V_2\otimes H)
					\end{split}
				\end{equation}
				The weights for $\lambda\curvearrowright H^0\big(\mathbb P^{N'},\mathcal O(1)\big)$ are 
				\begin{equation}
					\sum_{i=1}^2w_im_i\dim V_i\dim H,\quad (m_1,m_2)\in\Lambda\setminus\Lambda_0. 
				\end{equation}
				When $(m_1,m_2)\in\Lambda\setminus\Lambda_0$, we have $m_1\leq P_1(M)$ and $m_1+m_2=P(M)$ and 
				\begin{equation}
					\begin{split}
						&\sum_{i=1}^2w_im_i\dim V_i\dim H\\
						=&w_2P\dim V\dim H+(w_1-w_2)m_1\dim V_1\dim H\\
						\leq &w_2P\dim V\dim H+(w_1-w_2)P_1\dim V_1\dim H\\
						=&\sum_{i=1}^2w_iP_i\dim V_i\dim H. 
					\end{split}
				\end{equation}
				So the maximal $\lambda$-weight on $H^0\big(\mathbb P^{N'},\mathcal O(1)\big)$ is $\sum_{i=1}^2w_iP_i\dim V_i\dim H$ and 
				\begin{equation}
					H^0\big(\mathbb P^{N'},\mathcal O(1)\big)_{\lambda=\max}\cong \bigwedge^{P_1}(V_1\otimes H)\otimes\;\bigwedge^{P_2}(V_2\otimes H). 
				\end{equation}
				
				For one direction of inclusion $Y_\tau\subseteq \big(\mathbb P^{N'}\big)^0_{\min}\cap \overline{Y_\tau}$, let $[(\mathcal E,q)]\in Y_\tau$ be a closed point. By \cite{HuybrechtsDaniel2010Tgom} Lemma 4.4.3, we have when considered in $\mathrm{Quot}_{V(-m)/B/\Bbbk}^{P(t),\mathcal L}$ 
				\begin{equation}
					\lim_{t\to 0}\lambda(t).[(\mathcal E,q)]=\Big[\Big(\bigoplus_{i=1}^2\mathcal E_i,\bigoplus_{i=1}^2q_i\Big)\Big]\in\prod_{i=1}^2\mathrm{Quot}_{V_i(-m)/B/\Bbbk}^{P_i(t),\mathcal L}. 
				\end{equation}
				The limit point is still in $Y_\tau$ and its image in $\mathbb P^{N'}$ represents the one-dimensional quotient 
				\begin{equation}
					\begin{split}
						\bigwedge^P(V\otimes H)\Big/W\twoheadrightarrow& \bigwedge^P\Big(\bigoplus_{i=1}^2H^0\big(B,\mathcal E_i(M)\big)\Big)\\
						\cong&\bigwedge^{P_1}H^0\big(B,\mathcal E_1(M)\big)\otimes\;\bigwedge^{P_2}H^0\big(B,\mathcal E_2(M)\big). 
					\end{split}
				\end{equation}
				The point $\lim\lambda(t).[(\mathcal E,q)]$ is in $\big(\mathbb P^{N'}\big)^0_{\min}$ if and only if the following composition is not zero 
				\begin{equation}
					\begin{split}
						&\bigwedge^{P_1}(V_1\otimes H)\otimes\;\bigwedge^{P_2}(V_2\otimes H)\\
						\subseteq&\bigwedge^P(V\otimes H)\Big/W\\
						\twoheadrightarrow&\bigwedge^{P_1}H^0\big(B,\mathcal E_1(M)\big)\otimes\;\bigwedge^{P_2}H^0\big(B,\mathcal E_2(M)\big). 
					\end{split}
				\end{equation}
				The composition above is surjective since it is the tensor product of exterior products of surjective morphisms
				\begin{equation}
					V_i\otimes H:=V_i\otimes H^0\big(B,\mathcal L^{mM'}\big)\twoheadrightarrow H^0\big(B,\mathcal E_i(M)\big),\quad 1\leq i\leq 2. 
				\end{equation}
				This proves $\lim\lambda(t).[(\mathcal E,q)]\in\big(\mathbb P^{N'}\big)^0_{\min}$ and this shows $[(\mathcal E,q)]\in \big(\mathbb P^{N'}\big)^0_{\min}$. This proves $Y_\tau\subseteq \big(\mathbb P^{N'}\big)^0_{\min}\cap \overline{Y_\tau}$. 
				
				For the other direction $\big(\mathbb P^{N'}\big)^0_{\min}\cap \overline{Y_\tau}\subseteq Y_\tau$, let $[(\mathcal E,q)]\in\mathrm{Quot}_{V(-m)/B/\Bbbk}^{P(t),\mathcal L}$ be a closed point such that $[(\mathcal E,q)]\in \big(\mathbb P^{N'}\big)^0_{\min}\cap \overline{Y_\tau}$. The point $[(\mathcal E,q)]\in \mathbb P\Big(\big(\bigwedge^P(V\otimes H)\big)^*\Big)$ corresponds to the one-dimensional quotient 
				\begin{equation}
					\pi:\bigwedge^P(V\otimes H)\twoheadrightarrow\bigwedge^PH^0\big(B,\mathcal E(M)\big)
				\end{equation}
				and $\pi(W)=0$ since $[(\mathcal E,q)]\in\mathbb P^{N'}$. 
				
				The following collection of coherent sheaves on $B$ is bounded 
				\begin{equation}
					\bigcup_{i=1}^2\Big\{\mathcal F: \mathcal F=\mathcal E_i \textrm{ or }\mathcal F=\mathcal E_{\leq i}\textrm{ for } [(\mathcal E,q)]\in\mathrm{Quot}_{V(-m)/B/\Bbbk}^{P(t),\mathcal L}\Big\}
				\end{equation}
				By \cite{HuybrechtsDaniel2010Tgom} Lemma 1.7.6, the set of Hilbert polynomials of the above collection is finite, i.e. the following set is finite 
				\begin{equation}
					\mathcal P:=\bigcup_{i=1}^2\Big\{P(\mathcal F,t):\mathcal F=\mathcal E_i \textrm{ or }\mathcal F=\mathcal E_{\leq i}\textrm{ for } [(\mathcal E,q)]\in\mathrm{Quot}_{V(-m)/B/\Bbbk}^{P(t),\mathcal L}\Big\}. 
				\end{equation}
				Also by \cite{HuybrechtsDaniel2010Tgom} Lemma 1.7.6, regularities of $\mathcal E_i$ and $\mathcal E_{\leq i}$ are bounded. We can enlarge $M$ such that 
				\begin{itemize}
					\item higher cohomology groups of $\mathcal E_{\leq i}(M)$ and $\mathcal E_i(M)$ are zero; 
					\item the following map is order-preserving and bijective 
					\begin{equation}
						\mathcal P\to \mathbb N,\quad P(\mathcal F,t)\mapsto P(\mathcal F,M). 
					\end{equation}
				\end{itemize}

				Let $m_i:=\dim H^0\big(B,\mathcal E_i(M)\big)$ for $1\leq i\leq 2$. When $M$ is large, higher cohomology groups of $\mathcal E_i(M)$ vanish, and then 
				\begin{equation}
					\begin{split}
						&m_i=P(\mathcal E_i,M),\quad i=1,2\\
						&m_1+m_2=P(M). 
					\end{split}
				\end{equation}
				We claim that $m_1\geq P_1$. Otherwise, the following composition is zero
				\begin{equation}
					\begin{split}
						&\bigwedge^{P_1}(V_1\otimes H)\otimes\;\bigwedge^{P_2}(V_2\otimes H)\\
						\subseteq&\bigwedge^P(V\otimes H)\\
						\xrightarrow{\pi}&\bigwedge^{P}H^0\big(B,\mathcal E(M)\big). 
					\end{split}
				\end{equation}
				which contradicts that $[(\mathcal E,q)]\in\big(\mathbb P^{N'}\big)^0_{\min}$. 
				
				Therefore we have either $m_1=P_1$ or $(m_1,m_2)\in\Lambda_0$. We exclude the case when $(m_1,m_2)\in\Lambda_0$. Otherwise, the following composition is non-zero 
				\begin{equation}
					\begin{split}
						&\bigwedge^{m_1}(V_1\otimes H)\otimes\;\bigwedge^{m_2}(V_2\otimes H)\\
						\subseteq&\bigwedge^P(V\otimes H)\\
						\xrightarrow{\pi}&\bigwedge^{P}H^0\big(B,\mathcal E(M)\big)
					\end{split}
				\end{equation}
				which contradicts $\pi(W)=0$. Therefore we have 
				\begin{equation}
					P(\mathcal E_i,M)=P_i(M),\quad i=1,2. 
				\end{equation}
				Since $P(\mathcal E_i,t), P_i(t)\in\mathcal P$, we have $P(\mathcal E_i,t)=P_i(t)$ and this shows $[(\mathcal E,q)]\in Y_\tau$. This proves $\big(\mathbb P^{N'}\big)^0_{\min}\cap\overline{Y_\tau}\subseteq Y_\tau$, thus $(3)$. 
			\end{proof}
		}

		Let $\overline{Y_\tau^{\mathrm{s}}}\hookrightarrow \overline{Y_\tau}$ be the scheme theoretic image of $Y_\tau^{\mathrm{s}}\to \overline{Y_\tau}$. 
		{
			\begin{proposition}\label{proposition: Y_tau^s closure and its Z_min^K-s}
				Assume $\prod_{i=1}^2R_i^{\mathrm{s}}\ne\emptyset$. For $M\gg m\gg0$, we have the following: 
				\begin{itemize}
					\item[(1)] $\big(\overline{Y_\tau^{\mathrm{s}}}\big)^0_{\min}=\overline{Y_\tau^{\mathrm{s}}}\cap Y_\tau$; 

					\item[(2)] Let $Z_{\min}\big(\overline{Y_\tau^{\mathrm{s}}}\big):=\Big(\big(\overline{Y_\tau^{\mathrm{s}}}\big)^0_{\min}\Big)^\lambda$  for $\lambda:\mathrm T_\tau\to \mathrm P_\tau$. Then 
					\begin{equation}
						\begin{split}
							Z_{\min}\big(\overline{Y_\tau^{\mathrm{s}}}\big)=&\overline{Y_\tau^{\mathrm{s}}}\cap \prod_{i=1}^2\mathrm{Quot}_{V_i(-m)/B/\Bbbk}^{P_i(t),\mathcal L}\\
							\prod_{i=1}^2R_i^{\mathrm{s}}\subseteq& Z_{\min}\big(\overline{Y_\tau^{\mathrm{s}}}\big)\hookrightarrow \prod_{i=1}^2\overline{R_i};
						\end{split}
					\end{equation}

					\item[(3)] Let $Z_{\min}^{\mathrm K_\tau,\mathrm{s}}\big(\overline{Y_\tau^{\mathrm{s}}}\big)\subseteq Z_{\min}\big(\overline{Y_\tau^{\mathrm{s}}}\big)$ denote the stable locus for the $\mathrm K_\tau$-linearisation. Then 
					\begin{equation}
						Z_{\min}^{\mathrm K_\tau,\mathrm{s}}\big(\overline{Y_\tau^{\mathrm{s}}}\big)=\prod_{i=1}^2R_i^{\mathrm{s}},\quad p_\tau^{-1}\big(Z_{\min}^{\mathrm K_\tau,\mathrm{s}}\big(\overline{Y_\tau^{\mathrm{s}}}\big)\big)=Y_\tau^{\mathrm{s}}. 
					\end{equation}
				\end{itemize}
			\end{proposition}
			\begin{proof}
				We consider the ample line bundle $\mathcal M|_{\overline{Y_\tau^{\mathrm{s}}}}$ on $\overline{Y_\tau^{\mathrm{s}}}$, which equals the pullback of $\mathcal O(1)$ along $\overline{Y_\tau^{\mathrm{s}}}\hookrightarrow \mathbb P\Big(\big(\bigwedge^P(V\otimes H)\big)^*\Big)$, for $P:=P(M)$. By Lemma \ref{lemma: Y_tau closure and Y_tau}, there is a linear subspace $\mathbb P^{N'}\hookrightarrow\mathbb P\Big(\big(\bigwedge^P(V\otimes H)\big)^*\Big)$ containing $\overline{Y_\tau}$. Then $\mathcal M|_{\overline{Y_\tau^{\mathrm{s}}}}\cong \mathcal O_{\mathbb P^{N'}}(1)|_{\overline{Y_\tau^{\mathrm{s}}}}$. 
				
				For $(1)$, both sides are open in $\overline{Y_\tau^{\mathrm{s}}}$, so we can check that the closed points coincide. We have that $\emptyset\ne\prod_{i=1}^2R_i^{\mathrm{s}}\subseteq\overline{Y_\tau^{\mathrm{s}}}\cap \big(\mathbb P^{N'}\big)^0_{\min}$, and then 
				\begin{equation}
					\big(\overline{Y_\tau^{\mathrm{s}}}\big)^0_{\min}=\overline{Y_\tau^{\mathrm{s}}}\cap \big(\mathbb P^{N'}\big)^0_{\min}. 
				\end{equation}
				By Lemma \ref{lemma: Y_tau closure and Y_tau}, $\overline{Y_\tau}\cap \big(\mathbb P^{N'}\big)^0_{\min}=Y_\tau$ and then 
				\begin{equation}
					\begin{split}
						\big(\overline{Y_\tau^{\mathrm{s}}}\big)^0_{\min}=&\overline{Y_\tau^{\mathrm{s}}}\cap \big(\mathbb P^{N'}\big)^0_{\min}\\
						=&\overline{Y_\tau^{\mathrm{s}}}\cap\overline{Y_\tau}\cap \big(\mathbb P^{N'}\big)^0_{\min}\\
						=&\overline{Y_\tau^{\mathrm{s}}}\cap Y_\tau
					\end{split}
				\end{equation}
				which proves $(1)$. 
				
				For $(2)$, let $Z_{\min}\big(\overline{Y_\tau^{\mathrm{s}}}\big):=\Big(\big(\overline{Y_\tau^{\mathrm{s}}}\big)^0_{\min}\Big)^\lambda$ denote the fixed subscheme with respect to $\lambda$. Recall that $(Y_\tau)^\lambda=\prod_{i=1}^2\mathrm{Quot}_{V_i(-m)/B/\Bbbk}^{P_i(t),\mathcal L}$ and the intersection $\overline{Y_\tau^{\mathrm{s}}}\cap (Y_\tau)^\lambda$ is defined as the fibre product over $\overline{Y_\tau}$. For the first equality, it suffices to show the existence of the following morphisms over $\overline{Y_\tau}$
				\begin{equation}
					\begin{split}
							Z_{\min}\big(\overline{Y_\tau^{\mathrm{s}}}\big)&\to \overline{Y_\tau^{\mathrm{s}}}\cap(Y_\tau)^\lambda\\
							\overline{Y_\tau^{\mathrm{s}}}\cap(Y_\tau)^\lambda&\to Z_{\min}\big(\overline{Y_\tau^{\mathrm{s}}}\big). 
					\end{split}
				\end{equation}
				It is easy that the morphism $Z_{\min}\big(\overline{Y_\tau^{\mathrm{s}}}\big)\to \overline{Y_\tau^{\mathrm{s}}}\cap (Y_\tau)^\lambda$ exists. The existence of $\overline{Y_\tau^{\mathrm{s}}}\cap (Y_\tau)^\lambda\to Z_{\min}\big(\overline{Y_\tau^{\mathrm{s}}}\big)$ is from that the morphism $\overline{Y_\tau^{\mathrm{s}}}\cap (Y_\tau)^\lambda\to \big(Y_\tau^{\mathrm{s}}\big)^0_{\min}$ is $\lambda$-equivariant with the trivial $\lambda$-action on $\overline{Y_\tau^{\mathrm{s}}}\cap (Y_\tau)^\lambda$. This proves $Z_{\min}\big(\overline{Y_\tau^{\mathrm{s}}}\big)=\overline{Y_\tau^{\mathrm{s}}}\cap(Y_\tau)^\lambda$. Note that the equality is between subschemes of $Y_\tau$. 
				
				The morphism $\prod_{i=1}^2R_i^{\mathrm{s}}\to Z_{\min}\big(\overline{Y_\tau^{\mathrm{s}}}\big)$ is from $\prod_{i=1}^2R_i^{\mathrm{s}}\to \big(Y_\tau^{\mathrm{s}}\big)^0_{\min}$, which is $\lambda$-equivariant with trivial action on the source. It is an open immersion, since the composition $\prod_{i=1}^2R_i^{\mathrm{s}}\to Z_{\min}\big(\overline{Y_\tau^{\mathrm{s}}}\big)\to \prod_{i=1}^2\mathrm{Quot}_{V_i(-m)/B/\Bbbk}^{P_i(t),\mathcal L}$ is an open immersion. 
				
				We then construct the morphism $Z_{\min}\big(\overline{Y_\tau^{\mathrm{s}}}\big)\to \prod_{i=1}^2\overline{R_i}$. By \href{https://stacks.math.columbia.edu/tag/01R8}{Lemma 01R8}, $\overline{Y_\tau^{\mathrm{s}}}\cap Y_\tau\hookrightarrow Y_\tau$ is the scheme theoretic image of $Y_\tau^{\mathrm{s}}\to Y_\tau$. The morphism $Y_\tau^{\mathrm{s}}\to Y_\tau$ factors through $p_\tau^{-1}\Big(\prod_{i=1}^2\overline{R_i}\Big)\hookrightarrow Y_\tau$, so the dashed arrow below exists by the universal property of scheme theoretic image
				\begin{equation}
					\begin{tikzcd}
						Y_\tau^{\mathrm{s}}\ar[r,"\subseteq"]\ar[rd]&\overline{Y_\tau^{\mathrm{s}}}\cap Y_\tau\ar[r,hook]\ar[d,dashed]&Y_\tau\\
						&p_\tau^{-1}\Big(\prod_{i=1}^2\overline{R_i}\Big)\ar[ru,hook]. 
					\end{tikzcd}
				\end{equation}
				We obtain a morphism $Z_{\min}\big(\overline{Y_\tau^{\mathrm{s}}}\big)\to \prod_{i=1}^2\overline{R_i}$ by taking fixed point subschemes with respect to $\lambda$ of $\overline{Y_\tau^{\mathrm{s}}}\cap Y_\tau\to p_\tau^{-1}\Big(\prod_{i=1}^2\overline{R_i}\Big)$. It is easy to see $Z_{\min}\big(\overline{Y_\tau^{\mathrm{s}}}\big)\to \prod_{i=1}^2\overline{R_i}$ is a closed immersion over $\overline{Y_\tau}$. This proves $(2)$. 
				
				For $(3)$, the second equality follows immediately from the first. Recall that $\mathrm K_\tau=\prod_{i=1}^2\mathrm{SL}(V_i)$. We have 
				\begin{equation}
					\Big(\prod_{i=1}^2\overline{R_i}\Big)^{\mathrm K_\tau,\mathrm{s}}=\prod_{i=1}^2\overline{R_i}^{\mathrm{SL}(V_i),\mathrm{s}}=\prod_{i=1}^2 R_i^{\mathrm{s}}
				\end{equation}
				where the first equality is a formula on stable loci of product linearisations (See \cite{HoskinsVictoria2012Qous} Lemma 6.6, Proposition 6.7), and the second equality is from Theorem \ref{theorem: simpson's construction of moduli of semistable sheaves} when $M'\gg m\gg0$. By \cite{MumfordDavid1994Git}, Theorem 1.19 and $(2)$ proved above, we have 
				\begin{equation}
					Z_{\min}^{\mathrm K_\tau,\mathrm{s}}\big(\overline{Y_\tau^{\mathrm{s}}}\big)=Z_{\min}\big(\overline{Y_\tau^{\mathrm{s}}}\big)\cap \Big(\prod_{i=1}^2\overline{R_i}\Big)^{\mathrm K_\tau,\mathrm{s}}=\prod_{i=1}^2R_i^{\mathrm{s}}
				\end{equation}
				which proves $(3)$. 
				
			\end{proof}
		}
	}
}

\printbibliography

\end{document}